\documentclass{amsart}
\usepackage{amscd,amsmath,amssymb,amsfonts}
\usepackage[cmtip, all]{xy}

\newtheorem{thm}{Theorem}[subsection]
\newtheorem{prop}[thm]{Proposition}
\newtheorem{lem}[thm]{Lemma}
\newtheorem{cor}[thm]{Corollary}

\newtheorem{claim}{Claim}
\renewcommand{\theclaim}{\kern-3pt}
\newtheorem{apthm}{Theorem}[section]
\newtheorem{approp}[apthm]{Proposition}
\newtheorem{aplem}[apthm]{Lemma}
\newtheorem{apcor}[apthm]{Corollary}

\newtheorem{IntroThm}{Theorem}
\newtheorem{IntroCor}{Corollary}

\theoremstyle{definition}
\newtheorem{Def}[thm]{Definition}

\theoremstyle{remark}
\newtheorem{rem}[thm]{Remark}
\newtheorem{rems}[thm]{Remarks}

\numberwithin{equation}{subsection}

\newcommand{\sA}{{\mathcal A}}
\newcommand{\sB}{{\mathcal B}}
\newcommand{\sC}{{\mathcal C}}
\newcommand{\sD}{{\mathcal D}}
\newcommand{\sE}{{\mathcal E}}
\newcommand{\sF}{{\mathcal F}}
\newcommand{\sG}{{\mathcal G}}
\newcommand{\sH}{{\mathcal H}}

\newcommand{\sK}{{\mathcal K}}
\newcommand{\sL}{{\mathcal L}}
\newcommand{\sM}{{\mathcal M}}

\newcommand{\sO}{{\mathcal O}}
\newcommand{\sP}{{\mathcal P}}

\newcommand{\sS}{{\mathcal S}}

\newcommand{\sZ}{{\mathcal Z}}
\newcommand{\A}{{\mathbb A}}

\newcommand{\G}{{\mathbb G}}
\renewcommand{\H}{{\mathbb H}}

\renewcommand{\P}{{\mathbb P}}
\newcommand{\Q}{{\mathbb Q}}

\newcommand{\W}{{\mathbb W}}

\newcommand{\Z}{{\mathbb Z}}

\renewcommand{\phi}{\varphi}

\renewcommand{\1}{{\mathbf{1}}}

\newcommand{\CH}{{\rm CH}}

\newcommand{\by}[1]{\overset{#1}{\to}}

\newcommand{\iso}{\by{\sim}}
\newcommand{\red}{{\rm red}}
\newcommand{\codim}{{\rm codim}}

\newcommand{\Div}{{\rm Div}}
\newcommand{\Hom}{{\rm Hom}}
\newcommand{\End}{{\rm End}}
\newcommand{\Ext}{{\rm Ext}}

\newcommand{\Spec}{{\rm Spec \,}}

\newcommand{\Gal}{{\rm Gal}}

\newcommand{\0}{\emptyset}
\newcommand{\sHom}{{\mathcal{H}{om}}}

\newcommand{\Nm}{{\operatorname{Nm}}}

\newcommand{\id}{{\operatorname{id}}}
\newcommand{\Zar}{{\text{\rm Zar}}} 
\newcommand{\Ord}{{\mathbf{Ord}}}

\newcommand{\Sch}{{\operatorname{\mathbf{Sch}}}}

\newcommand{\holim}{\mathop{{\rm holim}}}
\newcommand{\op}{{\text{\rm op}}}

\newcommand{\fib}{{\operatorname{\rm fib}}}

\newcommand{\Spt}{{\mathbf{Spt}}}

\newcommand{\Sm}{{\mathbf{Sm}}}
\newcommand{\cofib}{{\operatorname{\rm cofib}}}

\renewcommand{\lim}{\operatornamewithlimits{\varprojlim}}
\newcommand{\colim}{\operatornamewithlimits{\varinjlim}}

\newcommand{\GL}{{\operatorname{\rm GL}}}

 \newcommand{\Ab}{{\mathbf{Ab}}}

\newcommand{\cyc}{{\operatorname{\rm cyc}}}

\newcommand{\SH}{{\operatorname{\sS\sH}}}

\newcommand{\eff}{{\mathop{eff}}}

\newcommand{\DM}{{DM}}

\newcommand{\Nis}{{\operatorname{Nis}}}

\newcommand{\ds}{{/\kern-3pt/}}

\newcommand{\coker}{\operatorname{coker}}

\newcommand{\Nrd}{{\mathop{\rm{Nrd}}}}

\newcommand{\Cor}{{\mathop{\rm{Cor}}}}

\newcommand{\sCH}{\mathcal{C}\mathcal{H}}
\newcommand{\Deg}{{\mathop{\rm{deg}}}}
\newcommand{\reg}{{\mathop{\rm{reg}}}}
\newcommand{\Ralg}{\textbf{$R$-alg}}

\newcommand{\ess}{{\rm ess}}
\newcommand{\PSh}{{\mathop{PSh}}}
\newcommand{\Sh}{{\mathop{Sh}}}
\newcommand{\EM}{{\mathop{EM}}}
\newcommand{\SmCor}{{\mathop{SmCor}}}
\newcommand{\PST}{{\mathop{PST}}}
\newcommand{\mot}{{\mathop{mot}}}
\newcommand{\Mot}{{\mathop{Mot}}}
\newcommand{\gm}{{\mathop{gm}}}
\newcommand{\Sing}{{\operatorname{Sing}}}

\newcommand{\equi}{{\mathop{equi}}}
\newcommand{\SB}{{\operatorname{SB}}}
\newcommand{\et}{{\operatorname{\acute{e}t}}}
\newcommand{\cone}{{\operatorname{cone}}}

\newcommand{\Imm}{{\operatorname{Imm}}}
\newcommand{\qfin}{{\operatorname{q.fin}}}
\newcommand{\Sus}{{\operatorname{Sus}}}
\newcommand{\PGL}{{\operatorname{PGL}}}
\newcommand{\Aut}{{\operatorname{Aut}}}
\newcommand{\Mod}{{\operatorname{Mod}}}

\begin{document}

\title{Motives of Azumaya algebras}
\author{Bruno Kahn}
\address{Institut de Math\'ematiques de Jussieu\\175--179 rue du
Chevaleret\\75013 Paris\\France}
\email{kahn@math.jussieu.fr}
\author{Marc Levine}
\address{
Department of Mathematics\\
Northeastern University\\
Boston, MA 02115\\
USA}
\email{marc@neu.edu}


 
\renewcommand{\abstractname}{Abstract}
\begin{abstract}  We study the slice filtration for the $K$-theory of a sheaf of Azumaya
algebras $A$, and for the motive of a Severi-Brauer variety, the latter in the case of a
central simple algebra of prime degree over a field. Using the Beilinson-Lichtenbaum
conjecture, we apply our results to show the vanishing of $SK_2(A)$ for a central simple algebra
$A$ of square-free index.
\end{abstract}
\date{December 31, 2007}
\maketitle
\tableofcontents

\section*{Introduction} 
Voevodsky \cite{VoevSlice} has defined an analog of  the classical Postnikov tower in the setting of  motivic stable homotopy theory by replacing the classical suspension $\Sigma:=S^1\wedge-$ with $t$-suspension
$\Sigma_t:=\P^1\wedge-$; we call this construction the {\em motivic Postnikov tower}. In this paper, we study the motivic Postnikov tower in the category of $S^1$-spectra, $\SH_{S^1}(k)$, and its analog in the category of effective motives, $\DM^\eff(k)$. We concentrate on objects arising from a central simple algebra $A$ over a field $k$. In the setting of $S^1$-spectra, we look at the presheaf of the $K$-theory spectra $K^A$:
\[
Y\mapsto K^A(Y):=K(Y;A),
\]
where $K(Y;A)$ is the $K$-theory spectrum of the category of $\sO_Y\otimes_kA$-modules which are locally free as $\sO_Y$-modules. For motives, we study the motive $M(X)\in\DM^\eff(k)$, where $X$ is the Severi-Brauer variety of $A$. In the case of the Severi-Brauer variety, we are limited to the case of $A$ having prime degree.

Of course, $K^A$ is a twisted form of the presheaf $K$ of $K$-theory spectra $Y\mapsto K(Y)$
and $X$ is a twisted form of a projective space over $k$, so one would expect the layers in the
respective Postnikov towers of $K^A$ and $M(X)$ to be twisted forms of the layers for $K$ and
$M(\P^n)$. The second author has shown in \cite{LevineHC} that the $n$th layer for $K$ is the
Eilenberg-Maclane spectrum for the Tate motive $\Z(n)[2n]$; similarly, the direct sum
decomposition
\[
M(\P^N)=\oplus_{n=0}^N\Z(n)[2n]
\]
shows that $n$th layer for $M(\P^N)$ is $\Z(n)[2n]$ for $0\le n\le N$, and is 0 for $n$ outside this range. The twisted version of $\Z(n)$ turns out to be $\Z_A(n)$, where $\Z_A\subset \Z$ is the subsheaf of the constant sheaf with transfers having value $\Z_A(Y)\subset\Z(Y)=\Z$ equal to the image of the reduced norm map
\[
\Nrd:K_0(A\otimes_kk(Y))\to K_0(k(Y))=\Z.
\]
Here $Y$ is any smooth irreducible scheme over $k$. Letting $s_n$ and $s_n^\mot$ denote the $n$ layer of the motivic Postnikov tower in $\SH_{S^1}(k)$ and $\DM^\eff(k)$, respectively, and letting $\EM:\DM^\eff(k)\to \SH_{S^1}(k)$ denote the Eilenberg-Maclane functor \cite{OstRond}, our main results are
\begin{IntroThm}\label{Thm:Main1} Let $A$ be a central simple algebra over a field $k$. Then
\[
s_n(K^A)=\EM(\Z_A(n)[2n])
\]
for all $n\ge0$.
\end{IntroThm}

\begin{IntroThm}\label{Thm:Main2} Let $A$ be a central simple algebra over a field $k$ of prime 
degree $\ell$, $X:=\SB(A)$ the associated Severi-Brauer variety. Then 
\[
s_n^\mot(M(X))=\Z_{A^{\otimes n+1}}(n)[2n]
\]
for $0\le n\le \ell-1$, $0$ otherwise.
\end{IntroThm}
See theorems~\ref{thm:Slice} and \ref{thm:SBSlice}, respectively, in the body of the paper.

Since $s_nK^A$ and $s_n^\mot M(X)$ are the layers in the respective motivic Postnikov towers
\begin{gather*}
\ldots\to f_{n+1}K^A\to f_nK^A\to \ldots\to f_0K^A=K^A\\
0=f_\ell^\mot M(X)\to f_{\ell-1}^\mot M(X)\to\ldots\to f_0^\mot M(X)=M(X)
\end{gather*}
our computation of the layers gives us some handle on the spectral sequences
\[
E_2^{p,q}:=\pi_{-p-q}(s_{-q}K^A(Y))\Longrightarrow \pi_{-p-q}K^A(Y)
\]
and
\[
E_2^{p,q}:=\H^{p+q}(Y, s_{-q}^\mot M(X)(n))\Longrightarrow \H^{p+q}(Y,M(X)(n))
\]
arising from the towers.
In fact, we use a version of the first sequence to help compute the layers of $M(X)$. Putting in our computation of the layers into the $K^A$-spectral sequence gives us the spectral sequence
\[
E_2^{p,q}:=H^{p-q}(Y,\Z_A(-q))\Longrightarrow K_{-p-q}(Y;A)
\]
generalizing the Bloch-Lichtenbaum/Friedlander-Suslin spectral sequence from motivic cohomology to $K$-theory \cite{BL, FriedSus}. In particular, taking $Y=\Spec k$, we get
\[
K_1(A)=H^1(k,\Z_A(1))
\]
and for $A$ of square-free index
\[
K_2(A)=H^2(k,\Z_A(2)).
\]
See theorem~\ref{thm:Comp} and theorem~\ref{thm:Codim1Van} .

To go further, we must use the Beilinson-Lichtenbaum conjecture. Recall that it is equivalent to the Milnor-Bloch-Kato conjecture relating Milnor's $K$-theory with Galois cohomology \cite{SuslinVoev}, \cite{GeisserLev}. 
It seems to be now a theorem (see \cite{WeibelOutline}), thanks to work of Rost and Voevodsky; accepted proofs are certainly that of Merkurjev and Suslin in the special case of weight $2$ \cite{K2Ref} and that of Voevodsky at the prime $2$ (in all weights) \cite{Voevmilnor}. Since this seems important to some people, we shall specify in what weights we need the Beilinson-Lichtenbaum (or Milnor-Bloch-Kato) conjecture for our statements.

We thus use our knowledge of the layers of $M(X)$, together with the Beilinson-Lichten\-baum conjecture, to deduce a result comparing $H^p(k,\Z_A(q))$ and $H^p(k,\Z(q))$ via the {\em reduced norm map}
\[
\Nrd:H^p(k,\Z_A(q))\to H^p(k,\Z(q)),
\]
this just being the map induced by the inclusion $\Z_A\subset \Z$. By identifying $\Nrd$ with
the change of topologies map from the Nisnevich to the \'etale topology (using the fact that
$\Z_A(n)^\et=\Z(n)^\et$),  a duality argument leads to

\begin{IntroCor}\label{Cor:Main1} Let $A$ be a central simple algebra of square-free index $e$ over $k$. Let $n\ge 0$ and assume the Beilinson-Lichtenbaum conjecture in weights $\le n+1$ at all primes dividing the index of $A$. Then
\[
\Nrd:H^p(k,\Z_A(n))\to H^p(k,\Z(n))
\]
is an isomorphism for $p<n$, and we have an exact sequence
\begin{multline*}
0\to H^n(k,\Z_A(n))\xrightarrow{\Nrd}H^n(k,\Z(n))\simeq K_n^M(k)\\
\xrightarrow{\cup[A]}H^{n+2}_\et(k,\Z/e(n+1))\to H^{n+2}_\et(k(X),\Z/e(n+1)).
\end{multline*}
\end{IntroCor}
Here $[A]\in H^3_\et(k,\Z(1))=H^2_\et(k,\G_m)$ is the class of $A$ in the Brauer group of $k$, and the map $\cup[A]$ is shorthand for the composition
\[
H^n(k,\Z(n))\iso H^n_\et(k,\Z(n))\xrightarrow{\cup[A]}H^{n+3}_\et(k,\Z(n+1))
\]
(note that this cup-product map lands into $ {}_eH^{n+3}_\et(k,\Z(n+1))\simeq H^{n+2}_\et(k,\Z/e(n+1))$, the latter isomorphism being a consequence of the Beilinson-Lichtenbaum conjecture in weight $n+1$.)

See theorem~\ref{thm:SK0} in the body of the paper for this result.

Combining this result with our identification above of $K_1(A)$ and $K_2(A)$ as ``twisted Milnor $K$-theory" of $k$, we have (see  theorem~\ref{thm:SK0})

\begin{IntroCor}\label{Cor:Main2} Let $A$ be a central simple algebra over $k$  of square-free index. Then the reduced norm maps on  $K_0(A)$, $K_1(A)$ and $K_2(A)$
\[
\Nrd:K_n(A)\to K_n(k);\quad n=0, 1,2
\]
 are  injective; in fact, we have an exact sequence 
 \begin{multline*}
 0\to K_n(A)\xrightarrow{\Nrd}K_n(k)=H^n(k,\Z(n))\\\xrightarrow{\cup[A]}H^{n+2}_\et(k,\Z/e(n+1))
 \to H^{n+2}_\et(k(X),\Z/e(n+1))
 \end{multline*}
 for $n=0,1,2$. (For $n=2$ we need the Beilinson-Lichtenbaum conjecture in weight $3$.)
 \end{IntroCor}

The injectivity of $\Nrd$ on $K_1(A)$ is Wang's theorem \cite{Wang}, and it was proved for $K_2(A)$ and $A$ a quaternion algebra by Rost \cite{rost} and Merkurjev \cite{merk2}. They used it as a step towards the proof of the Milnor conjecture in degree $3$; conversely, the Milnor conjecture in degree $3$ was used in \cite[proof of Th. 9.3]{lens} to give a simple proof of the injectivity in this case. This proof was one of the starting points of the present paper.

For $n=0$, the exact sequence reduces to Amitsur's theorem that $\ker(Br(k)\to Br(k(X))$ is generated by the class of $A$ \cite{amitsur}. For $n=1$,  the exactness at $K_1(k)$ is due to Merkurjev-Suslin \cite[Th. 12.2]{K2Ref} and the exactness at $H^{3}_\et(k,\Z/e(2))$ could be extracted from Suslin \cite{suslin}. For $n=1$ and $A$ a quaternion algebra, the exactness at $H^{3}_\et(k,\Z/2)$ is due to Arason \cite[Satz 5.4]{Arason}. For $n=2$ and a quaternion algebra, it is due to Merkurjev \cite[Prop. 3.15]{merk}. 

The injectivity for
$K_2(A)$ with $A$ of arbitrary square-free index has also been announced recently by A. Merkurjev (joint with A. Suslin); their method also
relies on the Beilinson-Lichtenbaum conjecture, using it to give a computation of the
motivic cohomology of the ``\v{C}ech co-simplicial scheme" $\text{\v{C}}(X)$. 

We begin our paper in section ~\ref{sec:Postnikov} with a quick review of the motivic
Postnikov tower in $\SH_{S^1}(k)$ and $\DM^\eff(k)$, recalling the basic constructions as well
as some of the results from \cite{LevineHC} that we will need. In section~\ref{sec:Birat} we
recall some of the first author's theory of \emph{birational motives}\footnote{jointly with R.
Sujatha} as well as pointing out the role these motives play as the Tate twists of slices of an
arbitrary
$t$-spectrum. We proceed in section~\ref{sec:KA} to define and study the special case of the
birational motive
$\Z_A$ arising from a central simple algebra $A$ over $k$; we actually work in the more general
setting of a sheaf of Azumaya algebras on a scheme. In section~\ref{sec:SpecSeq} we prove our
first main result: we compute the slices of the ``homotopy coniveau tower" for the $G$-theory
spectrum $G(X;\sA)$, where $\sA$ is a sheaf of Azumaya algebras on a scheme $X$. This result
relies on some regularity properties of the functors $K_p(-,A)$  which rely on results due to
Vorst and generalized by van der Kallen; we collect and prove what we need in this direction in
the appendix~\ref{AppSec:Reg}. We also recall some basic results on Azumaya algebras in the
appendix~\ref{AppSec:Azumaya}. Specializing to the case in which $X$ is smooth over a field $k$
and $\sA$ is the pull-back to $X$ of a central simple algebra $A$ over $k$, the results of
\cite{LevineHC} translate our computation of the slices of the homotopy coniveau tower to give
theorem~\ref{Thm:Main1}. 

We also give in Subsection \ref{Sect:4.9} a construction of homomorphisms from $SK_1$ and $SK_2$ of a central simple algebra $A$ to quotients of \'etale cohomology groups of $k$, in the spirit of an idea of Suslin \cite{suslin2, suslin3}, albeit with a very different technique (for $SK_2$ we need the Beilinson-Lichtenbaum conjecture in weight $3$).

 We turn to our study of the motive of a Severi-Brauer variety in
section~\ref{sec:SeveriBrauer}, proving theorem~\ref{Thm:Main2} there. We conclude in
section~\ref{sec:Applications} with a discussion of the reduced norm map and the proofs of
corollaries~\ref{Cor:Main1} and ~\ref{Cor:Main2}.\\

\noindent{\bf Acknowledgements.} The first author would like to thank Philippe
Gille for helpful exchanges about Azumaya algebras and Nicolas Perrin for an
enlightening discussion about the Riemann-Roch theorem. We also thank Wilberd van der Kallen for helpful comments. This work was begun when the second
author was visiting the Institute of Mathematics of Jussieu on a ``Poste rouge CNRS" in  2000,
for which visit the second author expresses his heartfelt gratitude. In addition, the second
author thanks the NSF for support via grants DMS-9876729, DMS-0140445 and DMS-0457195, as well
the Humboldt Foundation for support through the Wolfgang Paul Prize and a Senior Research
Fellowship.

\section{The motivic Postnikov tower in $\SH_{S^1}(k)$ and $\DM^\eff(k)$}\label{sec:Postnikov}
In this section, we assume that $k$ is a perfect field. We review Voevodsky's construction of the motivic Postnikov tower in $\SH(k)$ and $\SH_{S^1}(k)$, as well as the  analog of the tower  in  $\DM^\eff(k)$. We also give the description of these towers in terms of the homotopy coniveau tower, following \cite{LevineHC}, and recall the main results of \cite{LevineHC} on  well-connected theories and on various generalizations of Bloch's higher Chow groups to  higher Chow groups with coefficients.

\subsection{Constructions in $\A^1$ stable homotopy theory}
We start with the unstable $\A^1$ homotopy category over $k$, $\sH(k)$, which is the homotopy category of the category $\PSh(k)$ of pointed presheaves of simplicial sets on $\Sm/k$, with respect to the Nisnevich- and $\A^1$-local model structure (see {\O}stv\ae r-R\"ondigs \cite{OstRond}). We let $T$ denote the pointed presheaf $(\P^1,\infty)$, and $\Sigma_T$ the operation of smash product with $T$.

 $\SH_{S^1}(k)$ is the homotopy category of the model category $\Spt_{S^1}(k)$ of presheaves of spectra on $\Sm/k$, where the model structure is given as in \cite{OstRond}. $\SH(k)$ is the homotopy category of the category $\Spt_T(k)$ of $T$-spectra, where we will take this to mean the category of $T$-spectra in $\Spt_{S^1}(k)$, that is, objects are sequences
\[
\sE:=(E_0, E_1,\ldots)
\]
with the $E_n\in \Spt_{S^1}(k)$, together with bonding maps
\[
\epsilon_n:\Sigma_TE_n\to E_{n+1}
\]
We give $\Spt_T(k)$ the model structure of spectrum objects as defined by Hovey \cite{Hovey}. 

 $\Sigma_T$ (on both $\Spt_{S^1}(k)$ and $\Spt_T(k)$) has as a right adjoint the $T$-loops functor $\Omega_T$, which passes to the homotopy categories by applying $\Omega_T$ to fibrant models. Concretely, for $E\in \Spt_{S^1}(k)$ we have
\[
\Omega_TE(X):=\fib(E(X\times\P^1)\to E(X\times\infty)).
\]
Note that $\Sigma_T$ is an equivalence on $\SH(k)$, with inverse $\Omega_T$, but this is not the case on $\SH_{S^1}(k)$.

\subsection{Postnikov towers for $T$-spectra and $S^1$-spectra}
Voevodsky \cite{VoevSlice} has defined a canonical tower on the motivic stable homotopy category $\SH(k)$, which we call the \emph{motivic Postnikov tower}. This is defined as follows: Let $\SH^\eff(k)\subset \SH(k)$ be the localizing subcategory generated by the $T$-suspension spectra of smooth $k$-schemes, $\Sigma^\infty_TX_+$, $X\in\Sm/k$. Acting by the equivalence $\Sigma_T$ of $\SH(k)$ gives us the tower of localizing subcategories
\[
\cdots\subset\Sigma^{n+1}_T\SH^\eff(k)\subset \Sigma^n_T\SH^\eff(k)\subset \cdots\subset \SH(k)
\]
for $n\in\Z$. Voevodsky notes that the inclusion $i_n:\Sigma^n_T\SH^\eff(k)\to \SH(k)$ admits a right adjoint $r_n: \SH(k)\to\Sigma^n_T\SH^\eff(k)$; let $f_n:=i_n\circ r_n$ with counit $f_n\to \id$. Thus, for each $\sE\in \SH(k)$, there is a canonical tower in $\SH(k)$
\begin{equation}\label{eqn:MotTTower}
\ldots\to f_{n+1}\sE\to f_n\sE\to\ldots\to \sE
\end{equation}
which we call the \emph{motivic Postnikov tower}. The cofiber $s_n\sE$ of $f_{n+1}\sE\to f_n\sE$ is the $n$th slice of $\sE$.

Voevodsky also defined the motivic Postnikov tower on the category of $S^1$-spectra, using the $T$-suspension operator  as above to define the localizing subcategories $\Sigma_T^n\SH_{S^1}(k)$, $n\ge0$, 
\[
\cdots\subset \Sigma_T^{n+1}\SH_{S^1}(k)\subset \Sigma_T^n\SH_{S^1}(k)\subset\cdots\subset \SH_{S^1}(k),
\]
with right adjoint $r_n:\SH_{S^1}(k)\to \Sigma_T^n\SH_{S^1}(k)$ to the inclusion. This gives the $T$-Postnikov tower for $S^1$-spectra
\begin{equation}\label{eqn:S1TTower}
\ldots\to f_{n+1}E\to f_nE\to\ldots\to f_0E=E.
\end{equation}
One difference from the tower for $T$-spectra is that the $S^1$ tower terminates at $f_0E=E$, whereas the tower for $T$ spectra is in general infinite in both directions.

We write $f_{n/n+r}E$ for the cofiber of $f_{n+r}E\to f_nE$; for $r=1$, we use the notation $s_n:=f_{n/n+1}$ to denote the $n$th  slice in the Postnikov tower. We will mainly be using the $S^1$ tower, so the overlap in notations with the tower in $\SH(k)$ should not cause any confusion.

We have the pair of adjoint functors
\[
\xymatrix{
\SH_{S^1}(k)\ar@<3pt>[r]^{\Sigma_T^\infty}&\SH(k)\ar@<3pt>[l]^{\Omega_T^\infty}
}
\]
respecting the two Postnikov towers. In addition, we have
\begin{align*}
&\Sigma_T\circ f_n=f_{n+1}\circ \Sigma_T\\
&\Omega_T\circ f_{n+1}=f_n\circ \Omega_T
\end{align*}
for $T$-spectra (and all $n$) and
\begin{equation}\label{eqn:S1LoopIso}
\Omega_T\circ f_{n+1}=f_n\circ \Omega_T
\end{equation}
 for $S^1$-spectra (for $n\ge0$). The identities for $T$ spectra follow from the fact that $\Sigma_T$ is an equivalence, while the result for $S^1$-spectra is more difficult, and is proved in \cite[theorem 7.4.2]{LevineHC}.
 
 \begin{rem} \label{rem:SliceSusp}The identity $f_n\circ \Omega_TE=\Omega_T\circ f_{n+1}E$ for an $S^1$ spectrum $E$ gives by adjointness a map (in $\SH_{S^1}(k)$)
\[
\phi_{E,n}:\Sigma_Tf_n\Omega_TE\to f_{n+1}E
\]
which is not in general an isomorphism. If however $E=\Omega^\infty_T\sE$ for a $T$-spectrum $\sE\in\SH(k)$, then the isomorphism
\[
\Sigma_Tf_n\Omega_T\sE\cong f_{n+1}\Sigma_T\Omega_T\sE\cong f_{n+1}\sE
\]
passes to $E$, showing that $\phi_{E,n}$ is an isomorphism for all $n$.
\end{rem}

\subsection{Postnikov towers for motives}
There is an analogous picture for motives; we will describe the situation analogous to the $S^1$-spectra. The corresponding category of motives is  the enlargement $\DM^\eff(k)$ of $\DM^\eff_-(k)$. We recall briefly the construction.

The starting point is the category $\SmCor(k)$, with objects the smooth quasi-projective $k$-schemes $\Sm/k$, and morphisms given by the {\em finite correspondences} $\Cor_k(X,Y)$, this latter being the group of cycles on $X\times_kY$ generated by the integral closed subschemes $W\subset X\times_kY$ such that $W\to X$ is finite and surjective over some component of $X$. Composition is by the usual formula for composition of correspondences:
\[
W'\circ W:=p_{XZ*}(p_{XY}^*(W)\cdot p_{YZ}^*(W')).
\]
Sending $f:X\to Y$ to the graph $\Gamma_f\subset X\times_kY$ defines a functor $m:\Sm/k\to \SmCor(k)$.

Next, one has the category $\PST(k)$ of {\em presheaves with transfer}, this being simply the category of presheaves of abelian groups on $\SmCor(k)$. Restriction to $\Sm/k$ gives the functor 
\[
m^*:\PST(k)\to \PSh(\Sm/k);
\]
we let $\Sh^{tr}_\Nis(k)\subset \PST(k)$ be the full subcategory of $P$ such that $m^*(P)$ is a Nisnevich sheaf on $\Sm/k$. Such a $P$ is a {\em Nisnevich sheaf with transfers} on $\Sm/k$. 

The inclusion $\Sh^{tr}_\Nis(k)\to \PST(k)$ has as left adjoint the {\em sheafification functor}. $\PST(k)$ is an abelian category with kernel and cokernel defined pointwise;  as usual, $\Sh^{tr}_\Nis(k)$ is an abelian category with kernel the presheaf kernel and cokernel the sheafification of the presheaf cokernel.

Recall that the category $\DM^\eff_-(k)$ is the full subcategory of the bounded above derived category $D^-(\Sh^{tr}_\Nis((k))$ with objects the complexes $C^*$ for which the hypercohomology presheaves
\[
X\mapsto \H^n_\Nis(X,C^*)
\]
are $\A^1$ homotopy invariant for all $n$, i.e., the pull-back map
\[
 \H^n_\Nis(X,C^*)\to  \H^n_\Nis(X\times\A^1,C^*)
 \]
 is an isomorphism for all $n$ and for all $X$ in $\Sm/k$.
 
\begin{Def}  $\DM^\eff(k)$ is the full subcategory of $D(\Sh^{tr}_\Nis(k))$ consisting of complexes $C$ such that the hypercohomology presheaves  $\H^n_\Nis(-,C^*)$ are $\A^1$ homotopy invariant for all $n$.
\end{Def}

We note that there is a model structure on $C(\Sh^{tr}_\Nis(k))$ for which $\DM^\eff(k)$ is equivalent to the homotopy category of $C(\Sh^{tr}_\Nis(k))$ (see \cite{CisDeg, OstRond}); we will often use this result to lift constructions in $\DM^\eff(k)$ by taking fibrant models. Also, we have the localization functor
\[
RC_*^\Sus:D(\Sh^{tr}_\Nis(k))\to \DM^\eff(k)
\]
extending Voevodsky's localization functor
\[
RC_*^\Sus:D^-(\Sh^{tr}_\Nis(k))\to \DM_-^\eff(k)
\]
The equivalence of the subcategory of homotopy invariant objects in $D(\Sh^{tr}_\Nis(k))$ with the localization of $D^-(\Sh^{tr}_\Nis(k))$ is proved in \cite[\S 3.10]{CisDeg}.

We write $L(X)$ for the presheaf with transfers represented by $X\in\Sm/k$; this is in fact a Nisnevich sheaf with transfers. We have the Tate object $\Z(1)$ defined as the image in $\DM^\eff(k)$ of the complex
\[
L(\P^1)\to L(k)
\]
with $L(\P^1)$ in degree 2.

\begin{rem} The inclusion $D^-(\Sh^{tr}_\Nis(k))\to D(\Sh^{tr}_\Nis(k))$ is a full embedding, so Voevodsky's embedding theorem \cite[V, theorem 3.2.6]{FSV}, that the Suslin complex functor
\[
RC^\Sus_*\circ L:\DM^\eff_\gm(k)\to \DM^\eff_-(k)
\]
is a full embedding, yields the full embedding
\[
RC^\Sus_*\circ L:\DM^\eff_\gm(k)\to \DM^\eff(k)
\]
For $X\in\Sm/k$ we write $M(X)$ for the image of $L(X)$ in $\DM^\eff(k)$.
\end{rem}

The operation of the functor   $\otimes\Z(1)[2]$ on $\DM^\eff(k)$ gives the tower of localizing subcategories (for $n\ge0$)
\[
\cdots\subset\DM^\eff(k)(n+1)\subset \DM^\eff(k)(n)\subset \cdots\subset \DM^\eff(k)
\]
where $\DM^\eff(k)(n)$ is the localizing subcategory of $\DM^\eff(k)$ generated by objects $M(X)(n)[2n]$, $X\in\Sm/k$. Just as for $\SH_{S^1}(k)$, we have the right adjoint $r^\mot_n: \DM^\eff(k)\to \DM^\eff(k)(n)$ to the inclusion $i^\mot_n$. Thus, for $E$ in $\DM^\eff(k)$,  we  the \emph{motivic Postnikov tower} in $\DM^\eff(k)$
\begin{equation}\label{eqn:MotPostnikovTower}
\ldots\to f^\mot_{n+1}E\to f^\mot_nE\to\ldots\to f_0^\mot E=E
\end{equation}
with $f^\mot_n:=i^\mot_n\circ r^\mot_n$.

One has the pair of adjoint functors
\[
\xymatrix{
\SH_{S^1}(k)\ar@<3pt>[r]^{\Mot}&\ar@<3pt>[l]^{\EM}\DM^\eff(k) 
}
\]
where $\EM$ is the Eilenberg-Maclane functor, associating to a  presheaf of abelian groups the corresponding  presheaf of Eilenberg-Maclane spectra, and $\Mot$ associates to a presheaf of spectra $E:=(E_0, E_1,\ldots)$, first of all, the   presheaf of   singular chain complexes
\[
\Sing E:=\colim_n\Sing E_n[n], 
\]
where the   maps $\Sing E_n[n]\to \Sing E_{n+1}[n+1]$ are  induced by the bonding maps for $E$ and the natural map
\[
\Sigma \Sing E_n\to \Sing\Sigma E_n.
\]
 One then takes  the associated  freely generated complexes of presheaves with transfer,
 i.e.
 \[
 \Mot:=\Z^{tr}\circ\Sing.
 \]
 These functors respect the two towers of subcategories, and hence commute with the two truncation functors $f_n$ and $f^\mot_n$. It follows from work of {\O}stv{\ae}r-R\"ondigs \cite{OstRond} that the Eilenberg-Maclane functor $\EM$ is faithful and conservative.
 
 In particular, the identity \eqref{eqn:S1LoopIso} implies the identity
  \begin{equation}\label{eqn:MotLoopIso}
\Omega_T\circ f^\mot_{n+1}=f^\mot_n\circ \Omega_T
\end{equation}
for the motivic truncation functors.
 
Let $E$ be in $\SH_{S^1}(k)$, $Y\in \Sm/k$ and $W\subset Y$ a closed subset. We let $E^W(Y)$ denote the homotopy fiber of
\[
\tilde{E}(Y)\to \tilde{E}(Y\setminus W)
\]
where $\tilde{E}$ is a fibrant model of $E$ in $\Spt_{S^1}(k)$. We make a similar definition for $\sF\in C(\PST(k))$. If $E$ is homotopy invariant and satisfies Nisnevich excision, then the map of the homotopy fiber of  $E(Y)\to E(Y\setminus W)$ to $E^W(Y)$ is a weak equivalence; we will sometimes use this latter spectrum for $E^W(Y)$ without explicit mention.

Similarly, we lift the functors $s_n$, $f_n$ to operations on $\Spt_{S^1}(k)$ by taking the fibrant model of the corresponding object in $\SH_{S^1}(k)$; we make a similar lifting to $C(\Sh^{tr}_\Nis(k))$ for the functors $f^\mot_n$, $s^\mot_n$.

\subsection{Purity}
Let $i:W\to Y$ be a closed immersion in $\Sm/k$ such that the normal bundle $\nu:=N_{W/Y}$ admits a trivialization $\phi:\sO_W^q\to \nu$. This gives us the Morel-Voevodsky purity isomorphism  \cite[Theorem 2.23]{MorelVoev}  in $\SH$
\begin{equation}\label{eqn:purity}
\theta_{\phi,E}:E^W(Y)\to (\Omega^q_TE)(W)
\end{equation}
and the isomorphism on homotopy groups
\begin{equation}\label{eqn:purityn}
\theta_{\phi,n,E}:\pi_n(E^W(Y))\to \pi_n((\Omega^q_TE)(W)).
\end{equation}
In general, the $\theta_{\phi,n,E}$ depend on the choice of $\phi$.

For later use, we record the following result:

\begin{lem}\label{lem:Comp} Let $W\subset Y$ be a closed subset, $Y\in\Sm/k$, such  that $\codim_YW\ge q$ for some integer $q\ge0$. \\
\\
1. For $E\in\SH_{S^1}(k)$, the canonical map $f_qE\to E$ induces a weak equivalence
\[
(f_qE)^W(Y)\to E^W(Y)
\]
2. For $\sF\in\DM^\eff(k)$, the canonical map $f_q^\mot\sF\to \sF$ induces a weak equivalence
\[
(f_q^\mot\sF)^W(Y)\to \sF^W(Y)
\]
\end{lem}

\begin{proof} We prove (1), the proof of (2) is the parallel. Note that
\[
\pi_n(E^W(Y))\cong \Hom_{\SH_{S^1}(k)}(\Sigma^\infty(Y/Y\setminus W),\Sigma^{-n}E).
\]
and similarly for $\pi_n((f_qE)^W(Y))$.

Suppose at first that $W$ is smooth and has trivial normal bundle in $Y$, $\nu\cong\sO_W^p$, $p\ge q$. Then
\[
\Sigma^\infty(Y/Y\setminus W)\cong \Sigma_T^p(\Sigma^\infty W_+)
\]
hence $\Sigma^\infty(Y/Y\setminus W)$ is in $\Sigma_T^q\SH_{S^1}(k)$. In general, since $k$ is perfect, $W$ admits a filtration by closed subsets
\[
\0=W_{-1}\subset W_0\subset\ldots\subset W_N=W
\]
such that $W_n\setminus W_{n-1}$ is smooth and has trivial normal bundle in $Y\setminus W_{n+1}$. We have the homotopy cofiber sequence
\[
(Y\setminus W_{n+1})/(Y\setminus W_n)\to Y/(Y\setminus W_n)\to Y/(Y\setminus W_{n+1})
\]
so by induction, $\Sigma^\infty(Y/Y\setminus W)$ is in $\Sigma_T^q\SH_{S^1}(k)$. Thus, the universal property of $f_nE\to E$ implies
\[
\Hom_{\SH_{S^1}(k)}(\Sigma^\infty(Y/Y\setminus W),\Sigma^{-n}f_qE)\to
\Hom_{\SH_{S^1}(k)}(\Sigma^\infty(Y/Y\setminus W),\Sigma^{-n}E)
\]
is an isomorphism, as desired.
\end{proof}

\subsection{The homotopy coniveau tower} 
\begin{Def} 1. For $X\in\Sm/k$, and $q,n\ge0$ integers, set
\begin{align*}
\sS_X^{(q)}(n):=\{W\subset X\times\Delta^n\ | &W\text{ is closed}\\&\text{and }\codim_{X\times F}W\cap X\times F\ge q\\&\text{for all faces }F\subset\Delta^n\}
\end{align*}
Set
\begin{align*}
X^{(q)}(n):=\{w\in X\times\Delta^n\ | w&\text{ is the generic point of}\\&\text{some irreducible }W\in \sS_X^{(q)}(n) \}
\end{align*}
2. For $E\in\Spt_{S^1}(k)$, $X\in\Sm/k$ and integer $q\ge0$, define
\[
f^q(X,n;E)=\colim_{W\in\sS_X^{(q)}(n)}E^{W}(X\times\Delta^n)
\]
3. For $E\in\Spt_{S^1}(k)$, $X\in\Sm/k$ and integer $q\ge0$, define
\[
s^q(X,n;E)=\colim_{\substack{W\in\sS_X^{(q)}(n)\\W'\in\sS_X^{(q+1)}(n)}}E^{W\setminus W'}(X\times\Delta^n\setminus W')
\]
\end{Def}
For fixed $q$, $n\mapsto \sS^{(q)}_X(n)$ forms a simplicial set, and $n\mapsto f^q(X,n;E)$, $n\mapsto s^q(X,n;E)$ form simplicial spectra. We let $f^q(X,-;E)$ and $s^q(X,-;E)$ denote the respective total spectra.

For $\sF\in C(\PST(k))$, we make the analogous definition yielding the simplicial complexes $n\mapsto  f_\mot^q(X,n;\sF)$ and $n\mapsto  s_\mot^q(X,n;\sF)$; we let $f_\mot^q(X,*;\sF)$ and $s_\mot^q(X,*;\sF)$ be the associated total complexes.

\begin{prop}[\hbox{\cite[theorem 7.1.1]{LevineHC}}]\label{prop:HC}  Take $X\in \Sm/k$ and $q\ge0$ an integer. Let  $E\in\Spt_{S^1}(k)$ be homotopy invariant and satisfy Nisnevich excision. Then there are natural isomorphisms in $\SH$
\begin{align*}
&\alpha_{X,q;E}:f^q(X,-;E)\xrightarrow{\sim} f_q(E)(X)\\
&\beta_{X,q;E}:s^q(X,-;E)\xrightarrow{\sim} s_q(E)(X)
\end{align*}
\end{prop}

\begin{cor} \label{cor:HC}  Take $X\in \Sm/k$ and $q\ge0$ an integer. Let  $\sF\in C(\PST(k))$ be homotopy invariant and satisfy Nisnevich excision. Then there are natural isomorphisms in $D(\Ab)$
\begin{align*}
&\alpha_{X,q;\sF}:f_\mot^q(X,*;\sF)\xrightarrow{\sim} f^\mot_q(\sF)(X)\\
&\beta_{X,q;\sF}:s_\mot^q(X,*;\sF)\xrightarrow{\sim} s^\mot_q(\sF)(X)
\end{align*}
\end{cor}
Indeed, the corollary follows directly from proposition~\ref{prop:HC}  by using the Eilen\-berg-Maclane functor.

\begin{rems}\label{rem:HC} 1. The isomorphisms in proposition~\ref{prop:HC}(1) are natural   with respect to flat morphisms $X\to X'$ in $\Sm/k$ and with respect to maps  $E\to E'$ in $\Spt_{S^1}(k)$, for $E$, $E'$ which are homotopy invariant and satisfy Nisnevich excision. \\
\\
2. As each $W'\in\sS_X^{(q+1)}(n)$ is in $\sS_X^{(q)}(n)$, we have the natural maps
$f_{q+1}(X,n;E)\to f_q(X,n;E)$, compatible with the simplicial structure. Similarly, we have the natural restriction maps $f_q(X,n;E)\to s_q(X,n;E)$.  Since $E$ satisfies Nisnevich excision, the sequence
\[
f_{q+1}(X,n;E)\to f_q(X,n;E)\to s_q(X,n;E)
\]
is a homotopy cofiber sequence, giving the homotopy  cofiber sequence 
\[
f_{q+1}(X,-;E)\to f_q(X,-;E)\to s_q(X,-;E)
\]
on the total spectra. In addition, the diagram
\[
\xymatrix{
f_{q+1}(X,-;E)\ar[r]\ar[d]_{\alpha_{X,q+1;E}}& f_q(X,-;E)\ar[r]\ar[d]_{\alpha_{X,q;E}}& s_q(X,-;E)\ar[d]^{\beta_{X,q;E}}\\
f_{q+1}(E)(X)\ar[r]& f_q(E)(X)\ar[r]& s_q(E)(X)}
\]
commutes.\\
\\
4. The analogous statements hold for $\sF$ in $C(\PST(k))$ as in corollary~\ref{cor:HC}.
\end{rems}

For $E\in \SH_{S^1}(k)$, we have the diagram
\[
E\xleftarrow{\tau_q}f_qE\xrightarrow{\pi_q}s_qE
\]

\begin{lem} \label{lem:SliceComp} Take $E\in\SH_S{^1}(k)$, $X\in\Sm/k$ and  integers $q,n\ge0$.  For all $p\ge q$ the map $\tau_q:f_qE\to E$ induces  weak equivalences 
\begin{align*}
&f^p(X,n;f_qE)\xrightarrow{\tau_q}f^p(X,n;E)\\
&s^p(X,n;f_qE)\xrightarrow{\tau_q}s^p(X,n;E)
\end{align*}
\end{lem}

\begin{proof} That $\tau_q:f^p(X,n;f_qE)\to f^p(X,n;E)$ is a weak equivalence follows from lemma~\ref{lem:Comp}. We have the map of distinguished triangles
\[
\xymatrix{
f^{p+1}(X,n;f_qE)\ar[r]\ar[d]^{\tau_q}&f^p(X,n;f_qE)\ar[r]\ar[d]^{\tau_q}&s^p(X,n;f_qE)\ar[d]^{\tau_q}\\
f^{p+1}(X,n;E)\ar[r]&f^p(X,n;E)\ar[r]&s^p(X,n;E)}
\]
hence $\tau_q:s^p(X,n;f_qE)\to s^p(X,n;E)$ is also a weak equivalence.
\end{proof}

\begin{prop}\label{prop:SliceComp} Take $E\in\SH_S{^1}(k)$, $X\in\Sm/k$ and  integer $q\ge0$.\\
\\
1. For all $p\ge q$, the map $\tau_q:f_qE\to E$ induces  weak equivalences 
\begin{align*}
&f^p(X,-;f_qE)\xrightarrow{\tau_q}f^p(X,-;E)\\
&s^p(X,-;f_qE)\xrightarrow{\tau_q}s^p(X,-;E)
\end{align*}
2. The map $\pi_q:f_q\to s_q$ induces a weak equivalence
\[
s^q(X,-;f_qE)\xrightarrow{\pi_q}s^q(X,-;s_qE)
\]
\end{prop}

\begin{proof} (1) follows from lemma~\ref{lem:SliceComp}. For (2), we have the commutative diagram   in $\SH$
\[
\xymatrix{
s^q(X,-;f_qE)\ar[d]_{\beta_{X,q;f_qE}}\ar[r]^{\pi_q}&s^q(X,-;s_qE)\ar[d]^{\beta_{X,q;s_qE}}\\
s_q(f_qE)(X)\ar[r]_{s_q(\pi_q)}&s_q(s_qE)(X)}
\]
with vertical arrows isomorphisms. The bottom horizontal diagram extends to the distinguished triangle
\[
s_q(f_{q+1}E)\to s_q(f_qE)\xrightarrow{s_q(\pi_q)} s_q(s_qE)\to s_q(f_{q+1}E)[1]
\]
and we have the defining distinguished triangle for $s_q$:
\[
f_{q+1}(f_{q+1}E)\to f_q(f_{q+1}E)\to s_q(f_{q+1}E)\to f_{q+1}(f_{q+1}E)[1]
\]
Since $f_{q+1}E$ is in $\Sigma^{q+1}_T\SH_{S^1}(k)\subset \Sigma^{q}_T\SH_{S^1}(k)$, the canonical maps
\[
f_{q+1}(f_{q+1}E)\to f_{q+1}E,\ f_{q}(f_{q+1}E)\to f_{q+1}E
\]
are isomorphisms, hence $s_q(f_{q+1}E)\cong0$ and $s_q(\pi_q)$ is an isomorphism.
\end{proof}

\begin{rem} Making the evident changes, the analogs of lemma~\ref{lem:SliceComp} and proposition~\ref{prop:SliceComp} hold for $\sF\in \DM^\eff(k)$.
\end{rem}

\subsection{The 0th slice} \label{subsec:0thSlice}

Let $F$ be a presheaf of spectra on $\Sm/k$ which is $\A^1$-homotopy invariant and satisfies Nisnevich excision. Then $F$ is pointwise weakly equivalent to its fibrant model. In addition, these properties pass to $\sHom(X, F)$ for $X\in\Sm/k$. 

Furthermore, $(s_0F)(Y)$ can be described using the cosimplicial scheme of \emph{semi-local $\ell$-simplices} $\hat\Delta^\ell$ (denoted $\Delta^\ell_0$ in \cite{LevineHC}). In fact, for $Y\in\Sm/k$, let $\sO(\ell)_{k(Y),v}$ be the semi-local ring of the set $v$ of vertices of $\Delta^\ell_{k(Y)}$ and set
\[
\hat\Delta^\ell_{k(Y)}:=\Spec \sO(\ell)_{k(Y),v}.
\]
Clearly $\ell\mapsto \hat\Delta^\ell_{k(Y)}$ forms a cosimplicial subscheme of $\Delta^*_{k(Y)}$.  It follows from proposition~\ref{prop:HC} below that $(s_0F)(Y)$  weakly equivalent to total spectrum of the simplicial spectrum  
\[
\ell\mapsto  F(\hat\Delta^\ell_{k(Y)}),
\]
which we denote by $F(\hat\Delta^*_{k(Y)})$.

We have an analogous description of $s_0^\mot \sF(Y)$ for $\sF\in C(\PST(k))$ which is  $\A^1$-homotopy invariant and satisfies Nisnevich excision. Using the Eilenberg-Maclane functor, it follows from the case of $S^1$ spectra that $s_0^\mot \sF(Y)$ is represented by the total complex associated to the simplicial object of $C(\Ab)$ 
\[
\ell\mapsto \sF(\hat\Delta^\ell_{k(Y)}), 
\]
which we denote by  $\sF(\hat\Delta^*_{k(Y)})$.

\begin{rem}\label{rem:0SliceSupp} The 0th slice computes the $q$th slice with supports in codimension $\ge q$, as follows. Let $W\subset Y$ be a closed subset, $Y\in\Sm/k$. We have shown in \cite{LevineHC}
\begin{enumerate}
\item Let $W^0\subset W$ be an open subset of $W$ such that $W\setminus W^0$ has codimension $>q$ on $Y$, let $F:=W\setminus W^0$. Then the restriction 
\[
s_q(E)^W(Y)\to s_q(E)^{W^0}(Y\setminus F)
\]
is a weak equivalence. This follows directly from lemma~\ref{lem:Comp}.
\item Suppose that $W$ is smooth with trivial normal bundle $\nu$ in $Y$ and that $\codim_YW=q$.  A choice of trivialization $\psi:\nu\to \sO_W^q$ together with the  purity isomorphism \eqref{eqn:purity} gives an isomorphism
\[
s_q(E)^W(Y)\cong \Omega^q_T(s_q(E))(W)
\]
in $\SH$. Combining with the de-looping isomorphism \eqref{eqn:S1LoopIso} gives us the isomorphism
\[
\hat\theta_{E,W,Y,q}:s_q(E)^W(Y)\to s_0(\Omega^q_TE)(W)
\]
It is shown in \cite[corollary 4.2.4]{LevineHC} (essentially a consequence of proposition~\ref{prop:HC})  that $\hat\theta$ is in fact independent of the choice of trivialization $\psi$. 
\end{enumerate}
Let $Y_W^{(q)}$ be the set of generic points of $W$ of codimension exactly $q$ on $Y$. Combining (1) and (2) we have, for each $W\subset Y$ of codimension $\ge q$, a natural isomorphism
\[
\rho_{E,W,Y,q}:s_q(E)^W(Y)\to \amalg_{w\in Y_W^{(q)}}s_0(\Omega^q_TE)(k(w)).
\]
\end{rem}

\subsection{Connected spectra}\label{subsec:Conn} We continue to assume the field $k$ is perfect.

\begin{Def} Call $E\in\SH_{S^1}(k)$ {\em connected} if  for each  $X\in\Sm/k$, the spectrum  $\tilde{E}(X)$ is -1 connected, where $\tilde{E}\in\Spt_{S^1}(k)$ is a fibrant model for $E$.
\end{Def}

\begin{lem}\label{lem:ConnProps} Let $E\in \SH_{S^1}(k)$ be connected. Then \\
\\
1. For each $q\ge0$, $\Omega^q_TE$ is connected.\\
\\
2. For $X\in\Sm/k$ and $W\subset X$ a closed subset, the spectrum with supports $E^W(X)$ is -1 connected.\\
\\
3. Let $j:U\to X$ be an open immersion in $\Sm/k$, $W\subset X$ a closed subset . Then
\[
j^*:\pi_0(E^W(X))\to \pi_0(E^{W\cap U}(U))
\]
is surjective.
\end{lem}

\begin{proof} For (1) it suffices to prove the case $q=1$. Take $X\in\Sm/k$. Since $\infty\hookrightarrow \P^1$ is split by $\P^1\to\Spec k$,  $(\Omega_TE)(X)$ is a retract of $E(X\times\P^1)$. Since $E(X\times\P^1)$ is -1 connected by assumption, it follows that $(\Omega_TE)(X)$ is also -1 connected, hence $\Omega_TE$ is connected.

For (2), suppose first that $i:W\to X$ is a closed immersion in $\Sm/k$ and that the normal bundle $\nu$ of $W$ in $X$ admits a trivialization, $\nu\cong \sO_W^q$. We have the Morel-Voevodsky purity isomorphism \eqref{eqn:purity}
\[
E^W(X)\cong (\Omega^q_TE)(W).
\]
By (1) $(\Omega^q_TE)(W)$ is -1 connected, verifying (2) in this case.

In general, we proceed by descending induction on $\codim_XW$, starting with the trivial case $\codim_XW=\dim_kX+1$, i.e. $W=\0$ In general, suppose that $\codim_XW\ge q$ for some integer $q\le\dim_kX$. Then there is a closed subset $W'\subset W$ with $\codim_XW'>q$ such that $W\setminus W'$ is smooth and has trivial normal bundle in $X\setminus W'$. We have the homotopy fiber sequence
\[
E^{W'}(X)\to E^W(X)\to E^{W\setminus W'}(X\setminus W')
\]
thus the induction hypothesis, and the -1 connectedness of $E^{W\setminus W'}(X\setminus W')$ implies that $E^W(X)$ is -1 connected.

(3) follows from the homotopy fiber sequence
\[
E^{W\setminus U}(X)\to E^W(X)\to E^{W\cap U}(U)
\]
and the -1 connectedness of $E^{W\setminus U}(X)$.
\end{proof}

\begin{lem} \label{lem:SliceConn} Suppose $E\in\SH_{S^1}(k)$ is connected. Then for $X\in\Sm/k$ and every $q,n\ge0$, $f^q(X,n;E)$ and $s^q(X,n;E)$ are -1 connected.
\end{lem}

\begin{proof} This follows from lemma~\ref{lem:ConnProps}(2), noting that $f^q(X,n;E)$ and $s^q(X,n;E)$ are both colimits over spectra with supports $E^W(X\times\Delta^n)$, $E^{W\setminus W'}(X\times\Delta^n\setminus W')$.
\end{proof}

\begin{prop}  \label{prop:SliceConn} Suppose $E\in\SH_{S^1}(k)$ is connected. Then for every $q\ge0$, $f_qE$ and $s_qE$ are  connected.
\end{prop}

\begin{proof} Take $X\in \Sm/k$. We have isomorphism in $\SH$:
\[
f_qE(X)\cong f^q(X,-;E),\ s_qE(X)\cong s^q(X,-;E)
\]
By lemma~\ref{lem:SliceConn}, the total spectra $ f^q(X,-;E)$ and $s^q(X,-;E)$ are -1 connected, whence the result.
\end{proof}

\begin{Def} Fix an integer $q\ge0$ and let $W\subset Y$ be a closed subset with $Y\in\Sm/k$ and $\codim_YW\ge q$. For $E\in\SH_{S^1}(k)$, define the {\em comparison map}
\[
\psi^E_W(Y):\pi_0(E^W(Y))\to \pi_0(s_q(E)^W(Y))
\]
as the composition
\[
\pi_0(E^W(Y))\xleftarrow{\sim}\pi_0((f_qE)^W(Y))\to \pi_0(s_q(E)^W(Y))
\]
noting that $\pi_0((f_qE)^W(Y))\to \pi_0(E^W(Y))$ is an isomorphism by lemma~\ref{lem:Comp}.
\end{Def}

\begin{lem}\label{lem:pi0} Let $w\in Y^{(q)}$ be a codimension $q$ point of $Y\in\Sm/k$ and let $Y_w:=\Spec\sO_{Y,w}$. 
Take $E\in\SH_{S^1}(k)$ and suppose that $E$ is connected. Then the comparison map
\[
\psi^E_w(Y_w):\pi_0(E^w(Y_w))\to \pi_0(s_q(E)^w(Y_w))
\]
is an isomorphism.
\end{lem}

\begin{proof} Since $\hat\Delta^0_{k(Y)}=\Spec k(Y)$, we have the natural map
\[
\pi_0((\Omega^q_TE)(k(Y)))\to \pi_0((\Omega^q_TE)(\hat\Delta^*_{k(Y)}))
\]
which is an isomorphism. Indeed,  by lemma~\ref{lem:ConnProps}(1), $\Omega^q_TE$ is connected for all $q\ge0$. In particular,  $(\Omega^q_TE)(\hat\Delta^n_{k(Y)})$ is -1 connected for all $Y$ and all $n\ge0$. Thus  we have the presentation of $ \pi_0((\Omega^q_TE)(\hat\Delta^*_{k(Y)}))$:
\[
\pi_0((\Omega^q_TE)(\hat\Delta^1_{k(Y)}))\xrightarrow{i_0^*-i_1^*}\pi_0((\Omega^q_TE)(k(Y))\to \pi_0((\Omega^q_TE)(\hat\Delta^*_{k(Y)})).
\]
By lemma~\ref{lem:ConnProps}(3) and a limit argument, the map 
\[
\pi_0((\Omega^q_TE)(\Delta^1_{k(Y)}))\to \pi_0((\Omega^q_TE)(\hat\Delta^1_{k(Y)}))
\]
is surjective; since $\Delta^1_{k(Y)}=\A^1_{k(Y)}$ and $\Omega^q_TE$ is homotopy invariant, the map $i_0^*-i_1^*$ is the zero map.

Choose a trivialization of the normal bundle $\nu$ of $w\in Y_w$,  $k(w)^q\cong \nu$. This gives us the purity isomorphisms $E^w(Y_w)\cong (\Omega^q_TE)(w)$, $(s_qE)^w(Y_w)\cong s_0(\Omega^q_TE)(w)\cong (\Omega^q_TE)(\hat\Delta^*_{k(w)})$, giving the commutative diagram
\[
\xymatrix{
\pi_0(E^w(Y_w))\ar[r]^{\psi^E_w(Y_w)}\ar[d]&\pi_0(s_q(E)^w(Y_w))\ar[d]\\
\pi_0(\Omega^q_TE(w))\ar[r]&\pi_0((\Omega^q_TE)(\hat\Delta^*_{k(w)}))
}
\]
with the two vertical arrows and the bottom horizontal arrow isomorphisms. Thus $\psi^E_w(Y_w)$ is an isomorphism.
\end{proof}

\begin{lem} Suppose $E\in\SH_{S^1}(k)$ is connected. Fix an integer $q\ge0$ and let $W\subset Y$ be a closed subset, with $Y\in\Sm/k$ and $\codim_YW\ge q$. Then the comparison map 
\[
\psi^E_W(Y):\pi_0(E^W(Y))\to \pi_0(s_q(E)^W(Y))
\]
is surjective.
\end{lem}

\begin{proof}  Recall that $Y_W^{(q)}$ denotes the set of generic points $w$ of $W$ with $\codim_Y w=q$. Let $Y_W:=\Spec \sO_{Y,Y_W^{(q)}}$. By remark~\ref{rem:0SliceSupp}, the restriction map
\[
s_q(E)^W(Y)\to \amalg_{w\in Y_W^{(q)}}s_q(E)^w(Y_W)
\]
is a weak equivalence. By lemma~\ref{lem:pi0},  
\[
\psi^E_w(Y_W):\pi_0(E^w(Y_W))\to \pi_0(s_q(E)^w(Y_W))
\]
is an isomorphism for all $w\in Y_W^{(q)}$. Thus we have the commutative diagram
\[
\xymatrixcolsep{50pt}
\xymatrix{
\pi_0(E^W(Y))\ar[r]^{\psi^E_W(Y)}\ar[d]&\pi_0(s_q(E)^W(Y))\ar[d]\\
\oplus_{w\in Y_W^{(q)}}\pi_0(E^w(Y_W))\ar[r]_{\Sigma_w\psi^E_w(Y_W)}&\oplus_{w\in Y_W^{(q)}}
\pi_0(s_q(E)^w(Y_W)).}
\]
By  remark~\ref{rem:0SliceSupp}, the right hand vertical arrow is an isomorphism; the bottom horizontal arrow is an isomorphism by lemma~\ref{lem:pi0}. It follows from lemma~\ref{lem:ConnProps}(3) that the left hand vertical arrow is surjective, hence 
$\psi^E_W(Y)$ is surjective as well.
\end{proof}

\begin{lem}  Suppose that $E\in\Spt_{S^1}(k)$ is connected. Take $Y\in\Sm/k$, $w\in Y^{(q)}$ and let $Y_w:=\Spec\sO_{Y,w}$. Then the purity isomorphism
\[
\theta_{\phi,0,E}:\pi_0(E^w(Y_w))\to \pi_0(\Omega^q_TE(w))
\]
is independent of the choice of trivialization $\phi$.
\end{lem}

\begin{proof} We have the commutative diagram of isomorphisms
\[
\xymatrix{
\pi_0(E^w(Y_w))\ar[r]^{\psi^E_w(Y_w)}\ar[d]_{\theta_{\phi,0E}}&\pi_0(s_q(E)^w(Y_w))\ar[d]^{\theta_{\phi,0,s_qE}}\\
\pi_0(\Omega^q_TE(w))\ar[r]&\pi_0((\Omega^q_TE)(\hat\Delta^*_{k(w)}))
}
\]
By \cite[corollary 4.2.4]{LevineHC}, $\theta_{\phi,0,s_qE}$ is independent of the choice of $\phi$, whence the result.
\end{proof}

Take $E\in \SH_{S^1}(k)$ connected. For each closed subset $W\subset Y$, $Y\in\Sm/k$,  we have the canonical map
\[
\rho_{E,Y,W}:E^W(Y)\to \EM(\pi_0(E^W(Y))).
\]

\begin{Def}  Let $E\in\SH_{S^1}(k)$ be connected. Let $Y$ be in $\Sm/k$ and let $W\subset Y$ be a closed subset of codimension $\ge q$. The {\em cycle map}
\[
\cyc_E^W(Y):E^W(Y)\to \EM(\oplus_{w\in Y_W^{(q)}}\pi_0((\Omega^q_TE)(w)))
\]
is the composition
\begin{multline*}
E^W(Y)\xrightarrow{\rho_{E,Y,W}} \EM(\pi_0(E^W(Y)))\\
\xrightarrow{res}
 \EM(\oplus_{w\in Y_W^{(q)}}\pi_0(E^w(Y_W)))\\
 \xrightarrow{\theta_{\phi,0E}}
  \EM(\oplus_{w\in Y_W^{(q)}}\pi_0((\Omega^q_TE)(w))).
 \end{multline*}
 We let
 \[
\pi_0(\cyc_E^W(Y)):\pi_0(E^W(Y))\to \oplus_{w\in Y_W^{(q)}}\pi_0((\Omega^q_TE)(w))
\]
be the  the composition
\[
\pi_0(E^W(Y))
\xrightarrow{res}
\oplus_{w\in Y_W^{(q)}}\pi_0(E^w(Y_W))
 \xrightarrow{\theta_{\phi,0,E}}
\oplus_{w\in Y_W^{(q)}}\pi_0((\Omega^q_TE)(w)).
 \]
 \end{Def}
In other words, $\pi_0(\cyc_E^W(Y))$ is the map on $\pi_0$ induced by $\cyc_E^W(Y)$.

\begin{Def} Let $E\in\SH_{S^1}(k)$ be connected. For $X\in \Sm/k$ and integers $q,n\ge0$ define
\[
z^q(X,n;E):=\oplus_{w\in X^{(q)}(n)}\pi_0((\Omega^q_TE)(w)).
\]
\end{Def}
Taking the limit of the maps $\cyc^{W\setminus W'}_E(X\times\Delta^n\setminus W')$ for for $E\in \SH_{S^1}(k)$ connected, $W\in \sS_X^{(q)}(n)$, $W'\in  \sS_X^{(q+1)}(n)$ we have the maps of spectra
\[
\cyc_E(X,n):s^q(X,n;E)\to \EM(z^q(X,n;E))
\]
and the maps of abelian groups
\[
\pi_0(\cyc_E(X,n)):\pi_0(s^q(X,n;E))\to z^q(X,n;E)
\]

\begin{lem}\label{lem:CycSlice} Let $E\in\SH_{S^1}(k)$ be connected and let $X$ be in $\Sm/k$. Then
\[
\pi_0(\cyc_{s_qE}(X,n)):\pi_0(s^q(X,n;s_qE))\to z^q(X,n;s_qE)
\] 
is an isomorphism.
\end{lem}

\begin{proof} First note that, by proposition~\ref{prop:SliceConn}, $s_qE$ is connected, hence all terms in the statement are defined. By remark~\ref{rem:0SliceSupp}, the restriction map
\[
\pi_0((s_qE)^W(Y))\to \oplus_{w\in Y_W^{(q)}}\pi_0((s_qE)^w(Y_W))
\]
is an isomorphism; since $\pi_0(\cyc_{s_qE}(X,n))$ is constructed by composing restriction maps with purity isomorphisms, this proves the result.
\end{proof}

\begin{lem} \label{lem:CycSimp} Let $E\in\SH_{S^1}(k)$ be connected and let $X$ be in $\Sm/k$. There is a unique structure of a simplicial abelian group
\[
n\mapsto z^q(X,n;E)
\]
such that the maps $\pi_0(\cyc_E(X,n))$ define a map of simplicial abelian groups
\[
[n\mapsto \pi_0(s^q(X,n;E))]\xrightarrow{\pi_0(\cyc_E(X,-))}[n\mapsto z^q(X,n;E)].
\]
\end{lem}

\begin{proof} Since $E$ is connected, the cycle maps
\[
\pi_0(E^W(Y))\xrightarrow{res} \oplus_{w\in Y_W^{(q)}}\pi_0(E^w(Y_W))\cong
 \oplus_{w\in Y_W^{(q)}}\pi_0((\Omega^q_TE(w))
 \]
 are surjective. Thus $\pi_0(\cyc_E(X,n))$ is surjective, which proves the uniqueness.
 
 For existence,  the map $\pi_0(\cyc_E(X,n))$ is natural with respect to  $E$. In addition, by proposition~\ref{prop:SliceConn}, both $f_qE$ and $s_qE$ are connected;  applying $\pi_0(\cyc_?(X,n))$ to the diagram
 \[
 E\leftarrow f_qE\to s_qE
 \]
 gives the commutative diagram
 \[
 \xymatrix{
  \pi_0(s^q(X,n;E))\ar[d]_{\pi_0(\cyc_E(X,n))}&  \pi_0(s^q(X,n;f_qE))\ar[d]_{\pi_0(\cyc_{f_qE}(X,n))}\ar[r]\ar[l]&
  \pi_0(s^q(X,n;s_qE))\ar[d]_{\pi_0(\cyc_{s_qE}(X,n))}\\
  z^q(X,n;E)&z^q(X,n;f_qE)\ar[l]\ar[r]&z^q(X,n;s_qE)}
  \]
  By lemma~\ref{lem:SliceComp}, the left hand map in the top row is an isomorphism. The maps in the bottom rows are induced by maps
  \[
  \pi_0((\Omega^q_TE)(w))\leftarrow \pi_0((\Omega^q_Tf_qE)(w))\to  \pi_0((\Omega^q_Ts_qE)(w))
  \]
  By \eqref{eqn:S1LoopIso}, $\Omega^q_Tf_qE=f_0(\Omega^q_Tf_qE)=\Omega^q_TE$ and similarly 
  $\Omega^q_Ts_qE=s_0(\Omega^q_TE)$. Thus the bottom row is a sum of isomorphisms
  \[
   \pi_0((\Omega^q_TE)(w))\to  \pi_0(s_0(\Omega^q_TE)(w)).
   \]
   Finally, the right hand vertical map is an isomorphism by lemma~\ref{lem:CycSlice}. As the top row is the degree $n$ part of a diagram of maps of simplicial abelian groups, the isomorphisms 
   \[
    \pi_0(s^q(X,n;s_qE))\to z^q(X,n;s_qE)\leftarrow z^q(X,n;E)
    \]
 induce the structure of a simplicial abelian group from $[n\mapsto \pi_0(s^q(X,n;s_qE))]$ to 
 $[n\mapsto  z^q(X,n;E)]$, so that the maps $\pi_0(\cyc_E(X,n))$ define a map of simplicial abelian groups.
 \end{proof}

We use the above results to give a generalization of the higher cycle complexes of Bloch:
\begin{Def} Let $E\in\SH_{S^1}(k)$ be connected. For $X\in\Sm/k$, and $q,n\ge0$ integers, let $z^q(X,*;E) $ be the complex associated to the simplicial abelian group $n\mapsto z^q(X,n;E)$. 
Similarly, for $\sF\in C(\PST(\Sm/k))$ which is homotopy invariant and satisfies Nisnevich excision, we set
\[
z^q(X,n;\sF)= \oplus_{w\in X^{(q)}(n)}H^0((\Omega^q_T\sF)(w)),
\]
giving the simplicial abelian group $n\mapsto z^q(X,n;\sF)$. We denote the associated complex by  $z^q(X,*;\sF)$.

For integers $q,n\ge0$, set
\[
\CH^q(X,n;E):=H_n(z^q(X,*;E))
\]
and
\[
\CH^q(X,n;\sF):=H_n(z^q(X,*;\sF))
\]
\end{Def}

\subsection{Well-connected spectra}\label{subsec:WellConn}

Following \cite{LevineHC} we have
\begin{Def}\label{Def:WC} $E\in\SH_{S^1}(k)$  is {\em well-connected} if 
\begin{enumerate}
\item $E$ is connected.
\item For each $Y\in\Sm/k$, and each $q\ge0$, the total spectrum $(\Omega^q_TE)(\hat\Delta^*_{k(Y)})$ has
\[
\pi_n((\Omega^q_TE)(\hat\Delta^*_{k(Y)}))=0
\]
for $n\neq0$.
\end{enumerate}
\end{Def}

\begin{rem} Under the Eilenberg-Maclane map, the corresponding notion in $\DM^\eff(k)$ is: Let $\sF\in C(\PST(k))$ be $\A^1$ homotopy invariant and satisfy Nisnevich excision. Call $\sF$ well-connected if 
\begin{enumerate}
\item  $\sF$ is connected
\item For each $Y\in\Sm/k$, the total  complex $(\Omega^q_T\sF)(\hat\Delta^*_{k(Y)})$ satisfies
\[
H^n((\Omega^q_T\sF)(\hat\Delta^*_{k(Y)})))=0
\]
for $n\neq0$.
\end{enumerate}
\end{rem}

\begin{rem} We gave a slightly different definition of well-connectedness in \cite[definition 6.1.1]{LevineHC}, replacing the connectedness condition (1) with: $E^W(Y)$ is -1 connected for all closed subsets $W\subset Y$, $Y\in\Sm/k$. By lemma~\ref{lem:ConnProps}, this condition is equivalent with the connectedness of $E$.
\end{rem}

The main result on well-connected spectra is:
\begin{thm}\label{thm:WellConn} 1. Suppose $E\in\SH_{S^1}(k)$ is well-connected. Then
\[
\cyc_E(X):s^q(X,-;E)\to \EM(z^q(X,-;E))
\]
is a weak equivalence. In particular, there is a natural isomorphism
\[
\CH^q(X,n;E)\cong \pi_n((s_qE)(X))\cong \Hom_{\SH_{S^1}(k)}(\Sigma^\infty_TX_+,\Sigma^{-n}_ss_q(E)).
\]
2. Suppose $\sF\in C(\PST(k))$ is well-connected. Then 
\[
\cyc_E(X):s^q(X,*;\sF)\to z^q(X,*;\sF).
\]
is a quasi-isomorphism. In particular, there is a natural isomorphism
\[
\CH^q(X,n;\sF)\cong \H^{-n}_\Nis(X,s^\mot_q\sF)\cong \Hom_{\DM^\eff(k)}(M(X), s^\mot_q(\sF)[-n]).
\]
\end{thm}

\begin{proof} We prove (1), the proof of (2) is the same. We have  commutative diagram (in $\SH$)
\[
\xymatrix{
s^q(X,-;E)\ar[d]_{\cyc_E(X)}&\ar[d]_{\cyc_{f_qE}(X)}\ar[l]_{\tau_q} s^q(X,-;f_qE)\ar[r]^{\pi_q}&s^q(X,-;s_qE)\ar[d]^{\cyc_{s_qE}(X)}\\
\EM(z^q(X,-;E))&\ar[l]^{\tau_q}\ar[r]^{\pi_q}\EM(z^q(X,-;f_qE))&\EM(z^q(X,-;s_qE))}
\]
By proposition~\ref{prop:SliceComp}, the arrows in the top row are isomorphisms. As we have seen in the proof of lemma~\ref{lem:CycSimp}  the arrows in the bottom row are also  isomorphisms. Thus, it suffices to prove the result with $E$ replaced by $s_qE$.

The map $\cyc_{s_qE}(X)$ is just the map on total spectra induced by the map on $n$-simplices
\[
\cyc_{s_qE}(X,n):s^q(X,n;s_qE)\to \EM(z^q(X,n;s_qE))
\]
By lemma~\ref{lem:CycSlice}, the map on $\pi_0$, 
\[
\pi_0(\cyc_{s_qE}(X,n)):\pi_0(s^q(X,n;s_qE))\to z^q(X,n;s_qE),
\] 
is an isomorphism. However, since $E$ is well-connected, and since
\[
s^q(X,n;s_qE)\cong \amalg_{w\in X^{(q)}(n)}s_0(\Omega^q_TE)(k(w)),
\]
it follows that 
\[
s^q(X,n;s_qE)=\EM(\pi_0(s^q(X,n;s_qE))), 
\]
and  $\cyc_{s_qE}(X,n)$ is the map  induced by $\pi_0(\cyc_{s_qE}(X,n))$. Thus $\cyc_{s_qE}(X,n)$ is a weak equivalence for every $n$, hence $\cyc_{s_qE}(X)$ is an isomorphism in $\SH$, as desired.
\end{proof}

\section{Birational motives and higher Chow groups}\label{sec:Birat} Birational motives have been introduced and studied by Kahn-Sujatha \cite{KahnSujatha} and Huber-Kahn \cite{HuberKahn}. In this section we re-examine their theory, emphasizing the relation to the slices in the motivic Postnikov tower. We also  extend Bloch's construction of cycle complexes and higher Chow groups: Bloch's construction may be considered as the case of the cycle complex with constant coefficients $\Z$ whereas our generalization allows the coefficients to be in a  {\em birational motivic sheaf}. Finally, we extend the identification of Bloch's higher Chow groups with motivic cohomology \cite{FSV, VoevChow}  to the setting of birational motivic sheaves.

\subsection{Birational motives} 
\begin{Def} A motive $\sF\in \DM^\eff(k)$ is called {\em birational} if for every dense open immersion $j:U\to X$ in $\Sm/k$ and every integer $n$, the map
\[
j^*:\Hom_{\DM^\eff(k)}(M(X),\sF[n])\to \Hom_{\DM^\eff(k)}(M(U),\sF[n])
\]
is an isomorphism. If $\sF$ is a sheaf, i.e., $\sF\cong \sH^0(\sF)$ in $D(\Sh^{tr}_\Nis(k))$, we call $\sF$ a {\em birational motivic sheaf}. 
\end{Def}

\begin{rems}  1. For $X\in\Sm/k$ and $\sF\in \DM^\eff(k)\subset D(\Sh^{tr}_\Nis(k))$, there is a natural isomorphism
\[
\Hom_{\DM^\eff(k)}(M(X),\sF[n])\cong\H^n_\Nis(X,\sF)
\]
Thus a motive $\sF\in  \DM^\eff(k)$ is birational if an only if the hypercohomology presheaf
\[
U\mapsto \H^n_\Nis(U,\sF)
\]
on $X_\Zar$ is the constant presheaf on each connected component of $X$.\\
\\
2. Let $\sF$ be a Nisnevich sheaf with transfers that is birational and homotopy invariant. Then $\sF$ is a birational motivic sheaf. Indeed, since $\sF$ is birational, the restriction of $\sF$ to $X_\Zar$ is a locally constant sheaf. We have
\[
\Hom_{D(\Sh^{tr}_\Nis(k))}(L(X),\sF[n])=H^n_\Nis(X,\sF)=H^n_\Zar(X,\sF);
\]
the Zariski cohomology $H^n_\Zar(X,\sF)$ is zero for $n>0$ since a constant Zariski sheaf is flasque. In particular, $\sF$ is strictly homotopy invariant and thus an object of $\DM^\eff_-(k)\subset\DM^\eff(k)$. Finally
\[
\Hom_{\DM^\eff(k)}(M(X),\sF[n])=\Hom_{D(\Sh^{tr}_\Nis(k))}(L(X),\sF[n])=\begin{cases}\sF(X)&\text{ for }n=0\\0&\text{ for }n\neq0\end{cases}
\]
hence $\sF$ is a birational motive.

Of course, this result also follows from Voevodsky's theorem \cite{FSV} that a Nisnevich sheaf with transfers that is homotopy invariant is also strictly homotopy invariant, but the above argument avoids having to use this deep result.
\end{rems}

\subsection{The Postnikov tower for birational motives} In this section, we give a treatment
of the slices of a birational motive. These results are obtained in \cite{KahnSujatha};
here we develop part of the theory of  \cite{KahnSujatha} in a slightly different and
independent way.

Let $\sF$ be in $\DM^\eff(k)$. Since $f^\mot_0\sF\to \sF$ is an isomorphism, we have the canonical map
\[
\pi_0:\sF\to s_0^\mot\sF.
\]

The following result is taken from \cite{KahnSujatha} in slightly modified form:
\begin{thm}\label{thm:BiratS0}  For $\sF$ in $\DM^\eff(k)$, $\pi_0:\sF\to s_0^\mot\sF$ is an isomorphism if and only if $\sF$ is a birational motive.  In particular, since $s_0^\mot\sF=s^\mot_0(s_0^\mot\sF)$,   $s_0^\mot\sF$ is a birational motive.
\end{thm}

\begin{proof} Since we have the distinguished triangle
\[
f^\mot_1\sF\to \sF\xrightarrow{\pi_0} s^\mot_0\sF\to f^\mot_1\sF[1]
\]
$\pi_0$ is an isomorphism if and only if $f^\mot_1\sF\cong0$.

Suppose that $\pi_0$ is an isomorphism. Let $j:U\to X$ be a dense open immersion in $\Sm/k$ and let $W=X\setminus U$. We show that 
\[
\Hom_{\DM^\eff(k)}(M(X),\sF[n])\xrightarrow{j^*} \Hom_{\DM^\eff(k)}(M(U),\sF[n])
\]
is an isomorphism by induction on $\codim_XW$, starting with $\codim_XW=\dim_kX+1$, i.e., $W=\0$. We may assume that $X$ is irreducible.

By induction we may assume that $W$ is smooth of codimension $d\ge1$, giving us the Gysin distinguished triangle
\[
M(U)\xrightarrow{j}M(X)\to M(W)(d)[2d]\to M(U)[1].
\]
But as $d\ge1$, we have
\[
0=\Hom_{\DM^\eff(k)}(M(W)(d)[2d],f^\mot_1\sF[n])\cong \Hom_{\DM^\eff(k)}(M(W)(d)[2d],\sF[n])
\]
hence $j^*$ is an isomorphism.

Now suppose that $\sF$ is birational. We may assume that $\sF$ is fibrant as a complex of Nisnevich sheaves, so that
\[
\Hom_{\DM^\eff(k)}(M(X),\sF[n])=H^n(\sF(X))
\]
for all $X\in\Sm/k$.

 Take an  irreducible $X\in \Sm/k$. By remark~\ref{rem:0SliceSupp}(2) (applied with $Y=W=X$), we have a natural isomorphism
\[
\Hom_{\DM^\eff(k)}(M(X),s^\mot_0\sF[n])\cong H^n(\sF(\hat\Delta^*_{k(X)}))
\]
Also, as $\sF$ is birational, the restriction to the generic point gives an isomorphism
\[
\Hom_{\DM^\eff(k)}(M(X),\sF[n])\cong H^n(\sF(k(X))),
\]
and the map
\[
\Hom_{\DM^\eff(k)}(M(X),\sF[n])\xrightarrow{\pi_0}\Hom_{\DM^\eff(k)}(M(X),s^\mot_0\sF[n])
\]
is given by the map on $H^n$ induced by the canonical map
\[
\sF(k(X))=\sF(\hat\Delta^0_{k(X)})\to \sF(\hat\Delta^*_{k(X)}).
\]

On the other hand, since $\sF$ is birational, the map
\[
\sF(\Delta^n_{k(X)})\to \sF(\hat\Delta^n_{k(X)})
\]
is a quasi-isomorphism for all $n$, and hence the map of total complexes
\[
\sF(\Delta^*_{k(X)})\to \sF(\hat\Delta^*_{k(X)})
\]
is a quasi-isomorphism. Since $\sF$ is homotopy invariant, the map
\[
\sF(k(X))=\sF(\Delta^0_{k(X)})\to \sF(\Delta^*_{k(X)})
\]
is a quasi-isomorphism; thus the composition
\[
\sF(k(X))\to \sF(\Delta^*_{k(X)})\to \sF(\hat\Delta^*_{k(X)})
\]
is a quasi-isomorphism as well. Taking $H^n$, we see that 
\[
\Hom_{\DM^\eff(k)}(M(X),\sF[n])\xrightarrow{\pi_0}\Hom_{\DM^\eff(k)}(M(X),s^\mot_0\sF[n])
\]
is an isomorphism for all $X\in\Sm/k$. Since the localizing subcategory of $\DM^\eff(k)$ generated by the  $M(X)$ for $X\in\Sm/k$ is all of  $\DM^\eff(k)$, it follows that $\pi_0$ is an isomorphism.
\end{proof}

\begin{cor}\label{cor:BiratSlice} Let $\sF$ be a birational motive. Then
\[
f^\mot_m(\sF(n))=\begin{cases}0&\text{ for }m>n\\\sF(n)&\text{ for }m\le n.\end{cases}
\]
\end{cor}

\begin{proof}  Suppose $n\ge m\ge0$. As $\sF(n)$ is in $\DM^\eff(k)(m)$ , we have $f^\mot_m(\sF(n))=\sF(n)$.

Now take $m>n$. As a localizing subcategory of $\DM^\eff(k)$, $\DM^\eff(k)(m)$ is generated by objects $M(X)(m)$, $X\in\Sm/k$. Thus it suffices to show that 
\[
\Hom_{\DM^\eff(k)}(M(X)(m),\sF(n)[p])=0
\]
for all $X\in\Sm/k$ and all $p$.   By Voevodsky's  cancellation theorem \cite{voecan}, we have
\[
\Hom_{\DM^\eff(k)}(M(X)(m),\sF(n)[p])=\Hom_{\DM^\eff(k)}(M(X)(m-n),\sF[p])
\]
But since $m-n\ge1$, we have
\[
\Hom_{\DM^\eff(k)}(M(X)(m-n),\sF[p])\cong\Hom_{\DM^\eff(k)}(M(X)(m-n),f^\mot_1\sF[p])
\]
which is zero by theorem~\ref{thm:BiratS0}.
\end{proof}

\begin{rem}\label{rem:BiratSlice} Let $\sF$ be a birational motive. Then $\sF(n)=s^\mot_n(\sF(n))$ for all $n\ge0$. Indeed, $f^\mot_n(\sF(n))=\sF(n)$ and $f^\mot_{n+1}(\sF(n))=0$.
\end{rem}

\begin{rem}\label{rem:BiratVan} Let $\sF$ be a birational motive. Then for all $\sG$ in $\DM^\eff(k)$ and all integers $m>n\ge0$, we have
\[
\Hom_{\DM^\eff(k)}(\sG(m),\sF(n))=0
\]
Indeed, the universal property of $f^\mot_m\sF\to \sF$ gives the isomorphism
\[
\Hom_{\DM^\eff(k)}(\sG(m),f^\mot_m(\sF(n)))\cong \Hom_{\DM^\eff(k)}(\sG(m),\sF(n))
\]
but $f^\mot_m(\sF(n))=0$ by corollary~\ref{cor:BiratSlice}.
\end{rem}

\subsection{Cycles and slices} If $F/k$ is a finitely generated field extension, we define the motive $M(F)$ in $\DM^\eff(k)$ as the homotopy limit of the motives $M(Y)$ as $Y\in\Sm/k$ runs over all smooth models of $F$. Since we will really only be using the functor $\Hom_{\DM^\eff(k)}(M(F),-)$, the reader can, if she prefers, view  this as a notational short-hand for the functor on $\DM^\eff(k)$
\[
M\mapsto \colim_{\substack{Y\\k(Y)=F}}\Hom_{\DM^\eff(k)}(M(Y),M)
\]
This limit is just
\[
\lim_{\substack{Y\\k(Y)=F}}\H^0_\Zar(Y, M)
\]
in other words, just the stalk of the 0th  hypercohomology sheaf of $M$ at the generic point of $Y$.

\begin{lem}\label{lem:BiratVan2} Let $\sF$ be a  homotopy invariant Nisnevich sheaf with transfers. Then
\[
\Hom_{\DM^\eff(k)}(M(k(Y)),\sF(n)[2n+r]))=0
\]
for $r>0$ and for all $Y\in\Sm/k$.
\end{lem}

\begin{proof} Let $F=k(Y)$. $\sF(n)[2n]$ is a summand of $\sF\otimes M(\P^n)$, so it suffices to show that 
\[
\Hom_{\DM^\eff(k)}(M(F),\sF\otimes M(\P^n)[r]))=0
\]
for $r>0$. We can represent $\sF\otimes M(\P^n)$ by $C_*(\sF\otimes^{tr}L(\P^n))$. We have the canonical left resolution
\[
\sL(\sF)\to \sF
\]
 of $\sF$ (as a Nisnevich sheaf with transfers), where the terms in $\sL(\sF)$ are direct sums of representable sheaves, so we can replace $C_*(\sF\otimes^{tr}L(\P^n))$ with the total complex of
\[
\ldots\to C_*(\sL(\sF)_n\otimes L(\P^n))\to\ldots\to C_*(\sL(\sF)_0\otimes^{tr}L(\P^n))
\]
This in turn is a complex supported in degrees $\le 0$ with all terms direct sums of representable sheaves $L(Y)$, $Y\in\Sm/k$. But for any $X\in\Sm/k$, we have
\[
\Hom_{\DM^\eff(k)}(M(X), M(Y)[r])\cong \H^r_\Zar(X,C_*(Y))
\]
Thus 
\[
\Hom_{\DM^\eff(k)}(M(F), M(Y)[r])\cong H^r(C_*(Y)(F))
\]
which is zero for $r>0$, and thus 
\[
\Hom_{\DM^\eff(k)}(M(F),\sF(n)[2n+r]))\subset H^r(C_*(\sL(\sF)\otimes^{tr}L(\P^n)))=0
\]
for $r>0$.
\end{proof}

\begin{prop} Let $\sF$ be a birational motivic sheaf. Then $\sF(n)[2n]$ is well-connected. 
\end{prop}

\begin{proof} We first show that $\sF$ is connected, i.e., that
\[
\H^r_\Zar(X,\sF(n)[2n])=\Hom_{\DM^\eff(k)}(M(X), \sF(n)[2n+r])=0
\]
for all $r>0$ and all $X\in\Sm/k$. We have the Gersten-Quillen spectral sequence
\begin{multline*}
E_1^{p,q}=\oplus_{x\in X^{(p)}} \Hom_{\DM^\eff(k)}(M(k(x))(p)[2p], \sF(n)[2n+p+q])\\
\Longrightarrow
\Hom_{\DM^\eff(k)}(M(X), \sF(n)[2n+p+q]).
\end{multline*}
 For $p>n$, 
$E_1^{p,q}=0$ by remark~\ref{rem:BiratVan}. Using lemma~\ref{lem:BiratVan2} and Voevodsky's cancellation theorem \cite{voecan}, we see that $E_1^{p,q}=0$ for $p+q>0$, $p\le n$, whence the claim. 

Next, note that
\[
\Omega^m_T(\sF(n)[2n])=\begin{cases}\sF(n-m)[2n-2m]&\text{ for }0\le m\le n\\0&\text{ for }m>n.\end{cases}
\]
Indeed, note that, for $\sG\in \DM^\eff(k)$, 
\[
\Hom_{\DM^\eff(k)}(\sG,\Omega^m_T(\sF(n)[2n]))\cong\Hom_{\DM^\eff(k)}(\sG(m)[2m],\sF(n)[2n]).
\]
For $m\le n$, we have the canonical evaluation map $ev:\sF(n-m)[2n-2m]\to \Omega^m_T(\sF(n)[2n])$; the above identity says that $ev$ induces the Tate twist map 
\[
\Hom_{\DM^\eff(k)}(\sG,\sF(n-m)[2n-2m]))\to\Hom_{\DM^\eff(k)}(\sG(m)[2m],\sF(n)[2n]).
\]
Voevodsky's  cancellation theorem \cite{voecan} implies that the Tate twist map is an isomorphism; as $\sG$ was arbitrary, it follows that $ev$ is an isomorphism.  For the case $m>n$, the right-hand side
$\Hom_{\DM^\eff(k)}(\sG(m)[2m],\sF(n)[2n])$ is zero by remark~\ref{rem:BiratVan}.

Thus
\[
s_0^\mot(\Omega^m_T(\sF(n)[2n]))=\begin{cases}0&\text{ for }m\ge0, m\neq n\\\sF&\text{ for }m=n.
\end{cases}
\]
In fact, we need only check for $0\le m\le n$. If $0\le m<n$, then $\Omega^m_T(\sF(n)[2n])$ is in $\DM^\eff(k)(1)$, hence the $s^\mot_0(\Omega^m_T(\sF(n)[2n]))=0$. Finally, $\Omega^n_T(\sF(n)[2n])=\sF$, and thus $s_0\Omega^n_T(\sF(n)[2n])=s_0^\mot(\sF)=\sF$ by remark~\ref{rem:BiratSlice}.

In particular, $s_0^\mot(\Omega^m_T(\sF(n)[2n]))$ is concentrated in cohomological degree 0 for all $m$, which shows that $\sF(n)[2n]$ is well-connected.
\end{proof}

\begin{thm}\label{thm:BiratMotCoh} Let $\sF$ be a birational motivic sheaf. Then there is a natural isomorphism
\[
H^{2q-p}(X,\sF(q)):=\Hom_{\DM^\eff(k)}(M(X),\sF(q)[2q-p])\cong \CH^q(X,p;\sF(q)[2q])
\]
\end{thm}

\begin{proof} Since $\sF(q)[2q]$ is well-connected, it follows from theorem~\ref{thm:WellConn} that the slices $s^\mot_q(\sF(q)[2q])$ are computed by the cycle complexes, i.e., there is a natural isomorphism
\[
\Hom_{\DM^\eff(k)}(M(X),s^\mot_q(\sF(q)[2q])[-p])\cong \CH^q(X,p;\sF(q)[2q]).
\]
But $s^\mot_q(\sF(q)[2q])=\sF(q)[2q]$ by remark~\ref{rem:BiratSlice}.
\end{proof}

\begin{rem}\label{rem:coef} Let $\sF$ be a birational sheaf. For $Y\in\Sm/k$, we can define the group of codimension $q$ cycles on $Y$ with values in $\sF$ as
\[
z^q(Y)_\sF:=\oplus_{w\in Y^{(q)}}\sF(k(w)),
\]
that is, an $\sF$-valued cycle on $Y$ is a formal finite sum $\sum_i a_iW_i$ with each $W_i$ a codimension $q$ integral closed subscheme of $Y$ and $a_i\in\sF(k(W_i))$. The canonical identification 
\[
\sF(k(w))\cong H^0((\sF(q)[2q])^W(Y))
\]
for $W\subset Y$ a codimension $q$ integral closed subscheme gives the $\sF$-valued cycle groups the usual properties of algebraic cycles, including proper pushforward, and partially defined pull-back. In particular, for $\sF=\Z$, we have the identification 
\[
z^q(Y)_\Z=z^q(Y);
\]
we will show in the next section that this identification is compatible with the operations of proper pushforward, and pull-back (when defined).

In addition, we have
\begin{align*}
s_0^\mot(\Omega^q_T(\sF(q)[2q]))&\cong s_0^\mot(\sF)\\
&\cong\sF
\end{align*}
hence
\[
z^q(X,n;\sF(q)[2q])=\oplus_{w\in X^{(q)}(n)}\sF(k(w)).
\]
Thus we can think of  $z^q(X,*;\sF(q)[2q])$ as the cycle complex of codimension $q$ $\sF$-valued cycles in good position on $X\times\Delta^*$.
\end{rem}

\subsection{The sheaf $\Z$} The most basic example of a birational motivic sheaf is the constant sheaf $\Z$. Here we show that the constructions of the previous section are compatible with the classical operations on algebraic cycles.

Let $W\subset Y$ be a closed subset with $Y\in\Sm/k$. We let $z^q_W(Y)$ be the subgroup of $z^q(Y)$ consisting of cycles with support contained in $W$.

\begin{Def}\label{Def:Imm} The category of closed immersions $\Imm_k$ has objects $(Y,W)$ with $Y\in\Sm/k$ and $W\subset Y$ a closed subset. A morphism $f:(Y,W)\to (Y',W')$ is a morphism $f:Y\to Y'$  in $\Sm/k$ such that  $f^{-1}(W')_\red\subset W$. Let $\Imm_k(q)\subset \Imm_k$ be the full subcategory of closed subsets $W\subset Y$ such that each component of $W$ has codimension $\ge q$.
\end{Def}

 Note that for each morphism $f:(W\subset Y)\to (W'\subset Y')$, the pull-back of cycles gives a well-defined map $f^*:z^q_{W'}(Y')\to z^q_W(Y)$.

\begin{Def}Let $f:Y'\to Y$ be a morphism in $\Sch_k$, with $Y$ and $Y'$ equi-dimensional over $k$. We let $z^q(Y,*)_f \subset z^q(Y,*)$ be the subcomplex defined by letting $z^q(Y,n)_f$ be the subgroup of $z^q(Y,n)$ generated by irreducible $W\subset Y\times\Delta^n$, $W\in z^q(Y,n)$, such that for each face $F\subset \Delta^n$, each irreducible component of $(f\times\id_F)^{-1}(W\cap X\times F)$ has codimension $q$ on $Y'\times F$. 
\end{Def}

Assuming that $f(Y')$ is contained in the smooth locus of $Y$, the maps $(f\times\id_{\Delta^n})^*$ thus define the morphism of complexes 
\[
f^*:z^q(Y,*)_f\to z^q(Y',*)
\]

We recall Chow's moving lemma in the following form:
\begin{thm}[Bloch  \hbox{\cite{BlochWeb}}]\label{thm:CML} Suppose that $Y$ is a quasi-projective $k$-scheme, and that $f:Y'\to Y$ has image contained in the smooth locus of $Y$. Then the inclusion $:z^q(Y,*)_f\to :z^q(Y,*)$ is a quasi-isomorphism.
\end{thm}

One proves this by first using ``moving by translation" to prove the result for $Y=\P^n$, then using the method of the projecting cone to prove the result for $Y$ projective, and finally using Bloch's localization theorem, applied to a projective completion $Y\to \bar{Y}$, to prove the general case.

\begin{lem}\label{lem:CycComp} Take $Y\in\Sm/k$, $W\subset Y$ a closed subset. Suppose that each irreducible component of $W\subset Y$ has codimension $\ge q$. Then there is an isomorphism
\[
\rho_{Y,W,q}:H^{2q}_W(Y,\Z(q))\to z^q_W(Y)
\]
such that the $\rho_{Y,W,q}$ define a natural isomorphism of functors from $\Imm_{k,q}^\op$ to $\Ab$. In addition, the maps $\rho_{Y,W,q}$ are natural with respect to proper push-forward.
\end{lem}

\begin{proof} By definition, $\Z(1)[2]$ is the reduced motive of $\P^1$, 
\[
\Z(1)[2]=\tilde{M}(\P^1)\cong \cone(M(k)\xrightarrow{i_{\infty*}}M(\P^1)),
\]
and $\Z(q)[2q]$ is the $q$th tensor power of $\Z(1)[2]$. Via the localization functor
\[
RC^\Sus_*:D^-(\Sh^{tr}_\Nis(\Sm/k))\to \DM^\eff_-(k)
\]
and using \cite[V, corollary 4.1.8]{FSV}, we have the isomorphism
\[
\Z(q)[2q]\cong C^\Sus_*(z_\qfin(\A^q))
\]
and the natural identification
\[
H^{2q+p}(Y,\Z(q))\cong \H^p_\Nis(Y,C^\Sus_*(z_\qfin(\A^q)))\cong H^p(C^\Sus_*(z_\qfin(\A^q))(Y)).
\]
In particular, we have the natural identification of the motivic cohomology with supports
\[
H^{2q}_W(Y,\Z(q))\cong H_0(\cone(C^\Sus_*(z_\qfin(\A^q))(Y)\to C^\Sus_*(z_\qfin(\A^q))(Y\setminus W))[-1]).
\]
Set
\[
C^\Sus_*(z_\qfin(\A^q))(Y)_W:=\cone(C^\Sus_*(z_\qfin(\A^q))(Y)\to C^\Sus_*(z_\qfin(\A^q))(Y\setminus W))[-1].
\]

In addition, from the definition of the Suslin complex, we have the evident inclusion of complexes
\[
C^\Sus_*(z_\qfin(\A^q))(Y)\subset z^q(Y\times\A^q,*)_{f\times\id}\subset z^q(Y\times\A^q,*).
\]
It follows from  \cite[VI, theorem 3.2, V, theorem 4.2.2]{FSV} that the inclusion 
\[
C^\Sus_*(z_\qfin(\A^q))(Y)\subset z^q(Y\times\A^q,*)
\]
is a quasi-isomorphism; by theorem~\ref{thm:CML}, the inclusion 
\[
C^\Sus_*(z_\qfin(\A^q))(Y)\subset z^q(Y\times\A^q,*)_{f\times\id}
\]
is a quasi-isomorphism as well. 

Let $U=Y\setminus W$, $U':=Y'\setminus W'$ and let $f_U:U'\to U$ be the restriction of $f$. Setting
\[
z^q(Y,*)_{W,f}=\cone(z^q(Y,*)_{f}\to
z^q(U,*)_{f_U})[-1],
\]
we thus have the quasi-isomorphism 
\[
C^\Sus_*(z_\qfin(\A^q))(Y)_W\to  z^q(Y\times\A^q,*)_{W\times\A^q,f\times\id}.
\]

We have the commutative diagram
\[
\xymatrix{
C^\Sus_*(z_\qfin(\A^q))(Y)_W\ar[r]\ar[d]_{(f^*,f^*_U)}&  z^q(Y\times\A^q,*)_{W\times\A^q,f\times\id}\ar[d]^{(f\times\id^*, f_U\times\id^*)}\\
C^\Sus_*(z_\qfin(\A^q))(Y')_{W'}\ar[r]&  z^q(Y'\times\A^q,*)_{W'\times\A^q}}
\]
Since the horizontal maps are quasi-isomorphisms, we can use the right-hand side to compute $f^*:H^{2q}_W(Y,\Z(q))\to H^{2q}_{W'}(Y',\Z(q))$.

By the homotopy property for the higher Chow groups, and using the moving lemma again, the pull-back maps
\begin{align*}
&p_1^*:z^q(Y,*)_{W,f}\to z^q(Y\times\A^q,*)_{W\times\A^q,f\times\id}\\
&p_1^*:z^q(Y',*)_{W'}\to z^q(Y\times\A^q,*)_{W'\times\A^q}
\end{align*}
are quasi-isomorphisms. Thus we can use
\[
f^*:z^q(Y,*)_{W,f}\to z^q(Y',*)_{W'}
\]
to compute $f^*:H^{2q}_W(Y,\Z(q))\to H^{2q}_{W'}(Y',\Z(q))$.

Let $d=\dim_kY$. Chow's moving lemma together with the localization distinguished triangle
\[
z_{d-q}(W,*)\to z_{d-q}(Y,*)\to z_{d-q}(U,*)
\]
shows that the inclusion $z_{d-q}(W,*)\subset z_{d-q}(Y,*)_f$ induces a quasi-isomorphism
\[
z_{d-q}(W,*)\to z^q(Y,*)_{W,f}.
\]
Similarly, the inclusion $z_{d'-q}(W',*)\subset z_{d'-q}(Y',*)$, $d':=\dim_kY'$,  induces a quasi-isomorphism
\[
z_{d'-q}(W',*)\to z^q(Y',*)_{W'}.
\]
Since each component of $W$ has codimension $\ge q$ on $Y$, it follows that the inclusion
\[
z_{d-q}(W)=z_{d-q}(W,0)\to z_{d-q}(W,*)
\]
is a quasi-isomorphism. As $z_{d-q}(W)=z^q_W(Y)$, we thus have the isomorphism
\[
\rho_{Y,W,q}:z^q_W(Y)\to H^{2q}_W(Y,\Z(q))
\]
In addition, the diagram
\[
\xymatrix{
z_{d-q}(W)\ar@{=}[r]&z^q_W(Y)\ar[r]\ar[d]_{f^*}&z^q(Y,*)_{W,f}\ar[d]^{f^*}\\
z_{d'-q}(W')\ar@{=}[r]&z^q_{W'}(Y')\ar[r]&z^q(Y',*)_{W'}}
\]
commutes. Combining this with our previous identification of $H^{2q}_W(Y,\Z(q))$ with 
$H_0(z^q(Y,*)_{W,f})$ and $H^{2q}_{W'}(Y',\Z(q))$ with $H_0(z^q(Y',*)_{W'})$ shows that 
the isomorphisms $\rho_{Y,W,q}$ are natural with respect to pull-back.

The compatibility of the $\rho_{Y,W,q}$ with proper push-forward is similar, but easier, as one does not need to introduce the complexes $z^q(Y\times\A^q,*)_{f\times\id}$, etc., or use Chow's moving lemma. We leave the details to the reader.
\end{proof}

Now take $X\in\Sm/k$, $W\in\sS^{(q)}_X(n)$. We thus have the isomorphism
\[
\rho_{X\times\Delta^n,W,q}:H^{2q}_W(X\times\Delta^n,\Z(q))\to z^q_W(X\times\Delta^n) 
\]
In addition, if $W'\subset W$ is a closed subset of codimension $>q$ on $X\times\Delta^n$, then the restriction map
\[
H^{2q}_W(X\times\Delta^n,\Z(q))\to H^{2q}_{W\setminus W'}(X\times\Delta^n\setminus W',\Z(q))
\]
is an isomorphism.  Noting that
\[
H^0((\Z(q)[2q])^W(X\times\Delta^n))=H^{2q}_W(X\times\Delta^n,\Z(q))
\]
it follows from the  definition of $z^q(X,n;\Z(q)[2q])$ that we  have 
\[
z^q(X,n;\Z(q)[2q])=\colim_{\substack{W\subset X\times\Delta^n\\ 
W\in\sS^{(q)}_X(n)}}H^{2q}_W(X\times\Delta^n,\Z(q)).
\]
Thus taking the limit of the isomorphisms $\rho_{X\times\Delta^n,W,q}$ over $W\in \sS^{(q)}_X(n)$ gives the isomorphism
\[
\rho_{X,n}:z^q(X,n;\Z(q)[2q])\to z^q(X,n).
\]

\begin{prop}\label{prop:CycComp} For $X\in\Sm/k$, the maps $\rho_{X,n}$ define an isomorphism of complexes
\[
z^q(X,*;\Z(q)[2q])\xrightarrow{\rho_X}z^q(X,*)
\]
natural with respect to flat pull-back.
\end{prop}

\begin{proof}   It follows from lemma~\ref{lem:CycComp} that the isomorphisms $\rho_{X,W,n}$ are  natural with respect to the pull-back maps in $\Imm_k(q)$; in particular, with respect to flat pull-back and with respect to the face maps $X\times\Delta^{n-1}\to X\times\Delta^n$. Passing to the limit over $W\in\sS^{(q)}_X(n)$ proves the result.
\end{proof}

\section{The sheaves $\sK_0^\sA$ and $\Z_\sA$} \label{sec:KA} We apply the results of the previous sections to the $K$-theory of Azumaya algebras. The basic construction will be valid for a sheaf of Azumaya algebras over a fairly general base-scheme $X$; as the general theory we have already discussed is only available for $X=\Spec k$, we are forced to repeat some of the constructions in the more general setting before reducing the proof of the main result to the case $X=\Spec k$.

\subsection{$\sK_0^\sA$: definition and first properties}
Fix  a sheaf of Azumaya algebras $\sA$ on
an $R$-scheme of finite type $X$. For $p:Y\to X\in\Sch_X$, we have the sheaf
$p^*\sA$ of Azumaya algebras on $Y$. We may sheafify the $K$-groups of $p^*\sA$ for the Zariski topology on
$Y$, giving us the Zariski sheaves $\sK_n^\sA$ on $\Sch_X$.

\begin{lem}\label{l2.1} Suppose that $X$ is regular. Then 
\begin{enumerate}
\item  $\sK_0^\sA$ is an  $\A^1$ homotopy invariant presheaf on $\Sm/X$.
\item $\sK_0^\sA$ is a
\emph{birational} presheaf on $\Sm/X$, i.e., for $Y\in\Sm/X$, $j:U\to Y$ a dense open subscheme, the
restriction map
 \[
 j^*:\sK_0^\sA(Y)\to \sK_0^\sA(U)
 \]
 is an isomorphism.  Equivalently, $\sK_0^\sA$ is locally constant for the Zariski topology on $\Sm/X$,
hence is a sheaf for the Nisnevich topology on $\Sm/X$. 
\end{enumerate}
 \end{lem}
 
\begin{proof}  The homotopy invariance follows from the fact that $Y\mapsto K_0(Y;\sA)$ is homotopy invariant, and that the restriction map $K_0(Y,\sA)|to K_0(U,\sA)$ is surjective for each open immersion $U\to Y$ in $\Sm/X$.

For the birationality property, we may assume that $Y$ is irreducible. By corollary~\ref{cA1}, any object in the category $\sP_{X;\sA}$ is locally $\sA$-projective, hence  it suffices to show that for
each $y\in Y$, the map
 \[
 K_0(\sA\otimes_{\sO_B}\sO_{Y,y})\to
 K_0(\sA\otimes_{\sO_B}k(Y))
 \]
 is an isomorphism.
 
Since $Y$ is regular, surjectivity follows easily from corollary~\ref{cA2}. On the other hand, since $\sO_{Y,y}$ is local, the category of finitely generated projective
$\sA\otimes_{\sO_B}\sO_{Y,y}$ modules has a unique indecomposable generator (\cite{deMeyer}, \cite[III.5.2.2]{Knus}) and
similarly, the category of  finitely generated projective $\sA\otimes_{\sO_B}k(Y)$
modules has a unique simple generator. Thus the map is also injective, completing the proof that $\sK_0^\sA$ is birational.

To see that $\sK^\sA_0$ is a sheaf for the Nisnevich topology, it suffices to check the sheaf condition on elementary Nisnevich squares; this follows directly from the birationality property.
\end{proof}
 
\subsection{The reduced norm map}\label{subsec:RedNorm1} Let $\Spec F\to X$ be a point. We define a map 
\[
\Nrd_F:\Z\simeq K_0(\sA_F)\to K_0(F)=\Z
\]
by mapping the positive generator of $K_0(\sA_F)$ to $e_F[F]$, where $e_F$ is the \emph{index} of $\sA_F$. Recall that, by definition, $e_F^2=[D:F]$ where $D$ is the unique division $F$-algebra similar to $\sA_F$.

\begin{lem}\label{l2.2} The assignment $F\mapsto \Nrd_F$ defines a morphism of sheaves
 \[
 \Nrd:\sK_0^\sA \to \Z
 \]
which realizes $\sK_0^\sA$ as a subsheaf of the constant sheaf $\Z$ on $Y$. This is the \emph{reduced norm map} attached to $\sA$.
 \end{lem}
 
\begin{proof} In view of lemma~\ref{l2.1}, it suffices to check that if $L$ is a separable extension of $F$, the diagram
\[\begin{CD}
K_0(\sA_L)@>\Nrd_L>> K_0(L)\\
@AAA @AAA\\
K_0(\sA_F)@>\Nrd_K>> K_0(F)
\end{CD}\]
commutes. This is classical: by Morita invariance, we may replace $\sA_F$ by a similar division algebra $D$. Choose a maximal commutative subfield $E\subset D$ which is separable over $F$. First assume that $L=E$: then $D_L$ is split and  $\Nrd_L$ is an isomorphism by Morita invariance; on the other hand, the generator $[D]$ of $K_0(D)$ maps to $e$ times the generator of $K_0(D_L)$, which proves the claim in this special case. The general case reduces to the special case by considering a commutative cube involving the extension $LE$.
\end{proof}

 \subsection{The presheaf with transfers $\Z_\sA$}   For a scheme $X$ we let $\sM_X$ denote the category of coherent sheaves (of $\sO_X$ modules) on $X$. Given a sheaf of Azumaya algebras $\sA$ on $X$, we let $\sM_X(\sA)$ denote the category of sheaves of $\sA$-modules $\sF$ which are coherent as $\sO_X$-modules, using the structure map $\sO_X\to \sA$ to define the $\sO_X$-module structure on $\sF$. We let $G(X;\sA)$ denote the $K$-theory spectrum of the abelian category $\sM_X(\sA)$. If $f:Y\to X$ is a morphism, we often write $G(Y;\sA)$ for $G(Y;f^*\sA)$.
 
Suppose that $X$ is regular. Let $f:Z\to Y$ be a finite morphism in $\Sch_X$ with $Y$ in $\Sm/X$. Restriction of scalars defines a map of sheaves
 \[
 f_*:f_*\sK_0(Z;\sA)\to \sG_0(Y;\sA).
 \]
Using corollary~\ref{cA2}, we see that the natural map
 \[
 \sK_0(Y;\sA)\to \sG_0(Y;\sA)
 \]
 is an isomorphism, giving us the pushforward map
 \[
 f_*:\sK_0^\sA(Z)\to \sK_0^\sA(Y)
 \]

 Now take $Y, Y'\in\Sm/X$ and let $Z\subset Y\times_XY'$ be an integral subscheme which is finite over $Y$ and surjective onto a component of $Y$; let $p:Z\to Y$, $p':Z\to Y'$ be the maps induced by the projections. Define
 \[
 Z^*:\sK_0^\sA(Y')\to \sK_0^\sA(Y)
 \]
 by $Z^*:=p_*\circ p^{\prime*}$. For $X$ regular, this operation extends to $\Cor_X(Y,Y')$ by linearity.
 
 \begin{lem}\label{l2.3} Suppose $X$ regular. For $Z_1\in \Cor_X(Y,Y')$, $Z_2\in \Cor_X(Y',Y'')$ we have
 \[
 (Z_2\circ Z_1)^*=Z_1^*\circ Z_2^*
 \]
 \end{lem}
 
 \begin{proof} We already have a canonical operation of $\Cor_X(-,-)$ on the constant sheaf $\Z$ making $\Z$ a sheaf with transfers; one easily checks that this action agrees with the action we have defined above for $\sA=\sO_X$. It is similarly easy to check that, for $Z$ integral and $f:Z\to Y$ finite and surjective with $Y$ smooth, $f_*$ commutes with $\Nrd$. Since $\Nrd$ is injective, this implies that $\sK_0^\sA$ is also a sheaf with transfers, as desired.
 \end{proof}
 
 \begin{Def} Let $X$ be a regular $R$-scheme of finite type, $\sA$ a sheaf of Azumaya algebras on $X$. We let $\Z_\sA$ denote the Nisnevich sheaf with transfers on $\Sm/X$ defined by $\sK_0^\sA$.
 \end{Def}
 
 \begin{rem}\label{rem:Nrd} The reduced norm map $\Nrd:\sK_0^\sA\to\Z$ defines a map of Nisnevich sheaves with transfers $\Nrd:\Z_\sA\to \Z$.
 \end{rem}

\begin{lem} The subsheaf with transfers $(\Z_\sA,\Nrd)$ of the constant sheaf (with transfers) $\Z$  only depends on the subgroup of $Br(X)$ generated by $\sA$. In particular, it is Morita-invariant.
\end{lem}

\begin{proof} Indeed, if $\sB$ generates the same subgroup of $Br(X)$ as $\sA$, there exist integers $r,s$ such that $\sA^{\otimes_X s}$ is similar to $\sB$ and $\sB^{\otimes_X s}$ is similar to $\sA$. This implies readily that $\sA$ and $\sB$ have the same splitting fields (say, over a point $\Spec F$ of $X$), hence have the same index (say, over any extension of $F$).
\end{proof}

\begin{rem} The maps $K_0(F)\to K_0(\sA_F)$ given by extension of scalars also define a morphism of sheaves $\Z\to\Z_\sA$. But this morphism is not Morita-invariant.
\end{rem}

In case $X$ is the spectrum of a field, lemma~\ref{l2.1} yields

\begin{prop}\label{prop:AzumayaBiratSheaf} Take $X=\Spec k$, $k$ a field, and  let $A$ be a central simple algebra over $k$. Then the sheaf with transfers $\Z_A$ on $\Sm/k$ is a birational motivic sheaf.
\end{prop} 

\subsection{Severi-Brauer varieties}  Let $p:SB(\sA)\to X$ be the Severi-Brauer variety associated to $\sA$. 

\begin{lem} Suppose $X=\Spec k$, $k$ a field. Then the subgroup $\Nrd(K_0(\sA))\subset
K_0(k)=\Z$ is the same as the image 
 \[
 p_*(\CH_0(SB(\sA)))\subset \CH_0(B)=\Z.
 \]
 Moreover, $p_*:\CH_0(SB(\sA))\to\Z$ is injective.
 \end{lem} 
 
 \begin{proof} This is a theorem of Panin \cite{Panin}. We recall the proof of the first
statement. Let $x=\Spec K$ be a closed point of $SB(\sA)$. Then $K$ is a finite extension of
$F$ which is a splitting field of $\sA$. It is classical that $K$ is a maximal commutative
subfield of some algebra similar to $\sA$; in particular,  $[K:F]$ is divisible by the index of
$A$. Conversely, replacing $A$ by a similar division algebra $D$, for any maximal commutative
subfield $L\subset D$, $[L:F]$ equals the index of $A$. 
 \end{proof}
 
 Now let us come back to the case where $X$ is regular. Let us denote by $\sCH_0(SB(\sA)/X)$
the sheafification (for the Zariski topology) of the presheaf on $\Sm/X$
 \[
 U\mapsto \CH_{\dim_kU}(SB(\sA)\times_XU).
 \]
 The push-forward
 \[
 p_{U*}:\CH_{\dim_kU}(SB(\sA)\times_XU)\to \CH_{\dim_kU}(U)=\Z
 \]
 defines the map
 \[
\Deg:\sCH_0(SB(\sA)/X)\to \Z
\]
where $\Z$ is viewed as a constant sheaf on $(\Sm/X)_\Zar$. 

\begin{lem} The map $\Deg$ identifies $\sCH_0(SB(\sA)/X))$  with the locally constant subsheaf   $\Nrd(\Z_\sA)\subset \Z$. In other words, there is a canonical isomorphism of subsheaves of $\Z$
\[(\Z_\sA,\Nrd)\simeq (\sCH_0(SB(\sA)/X)),\Deg).\]
\end{lem}

\begin{proof} As we have already remarked, the result is true at $\Spec F$, $F$ a field. For $Y$ local, the restriction map
\[
j^*:\CH_{\dim X}(SB(\sA)\times_X Y)\to \CH_0(SB(\sA\otimes_{\sO_B}k(Y))
\]
($\dim X:=$ the Krull dimension) is surjective, from which the result easily follows.
\end{proof}

\begin{rem} It is evident that the transfer structure of lemma~\ref{l2.3} on $\Z_\sA$ coincides with the natural transfer structure on $\sCH_0(SB(\sA)/X))$ .
\end{rem}

\subsection{$\sK_0^\sA$ for embedded schemes} Let $k$ be a field. We fix a sheaf of Azumaya algebras $\sA$ on some finite type $k$-scheme $X$; we do not assume that $X$ is regular.

As a technical tool, we extend the definition of the category $\Imm_k$ (definition~\ref{Def:Imm}) as follows:

\begin{Def} The category of closed immersions $\Imm_{X,k}$ has objects $(Y,W)$ with $Y\in\Sm/k$ and $W\subset X\times_kY$ a closed subset. A morphism $f:(Y,W)\to (Y',W')$ is a morphism $f:Y\to Y'$  in $\Sm/k$ such that  $(\id\times f)^{-1}(W')_\red\subset W$.
\end{Def}

Let $Y$ be a smooth $k$-scheme,  let $i:W\to X\times_kY$ be a reduced closed subscheme of pure codimension. Letting $W_\reg\subset W$ be the regular locus, we have the (constant) Zariski sheaf $\sK_0^\sA$ defined on $W_\reg$. We describe how to extend $\sK_0^\sA$ to $W\subset X\times_kY$ so that 
\[
(Y,W)\mapsto \sK^\sA_0(W\subset X\times_kY)
\]
defines a presheaf $\sK_0^\sA$ on $\Imm_{X,k}$.

For this, we define $\sK_0^\sA$ on $i:W\to X\times_kY$ to be $\sK_0^\sA(W_\reg)$, where $j:W_\reg\to W$ is the regular locus of $W$. The trick is to define the pull-back maps.

We let $G^W(X\times_kY;\sA)$ denote the homotopy fiber of the restriction map
\[
G(X\times_kY;\sA)\to G(X\times_kY\setminus W;\sA)
\]

\begin{lem}\label{lem:QuillenGersten} Suppose that $X$ is local, with closed point $x$.  Let $i:Y'\to Y$ be a closed embedding in $\Sm^\ess/k$, with $Y$ local having closed point $y$. Let $W\subset X\times Y$ be a closed subset  such that $X\times Y'\cap W=(x,y)$ (as a closed subset). If $\codim_{X\times Y} W>\codim_{X\times Y'}(x,y)$, then the restriction map
\[
i^*:G_0^W(X\times Y;\sA)\to G_0^{(x,y)}(X\times Y';\sA)
\]
is the zero map.
\end{lem}

\begin{proof} The proof is a modification of Quillen's proof of Gersten's conjecture. Making a base-change to $k(x,y)$, and noting that $G_0^{(x,y)}(X\times Y;\sA))=G_0((x,y);\sA)$, we may assume that $k(y)=k(x)=k$. Since $K$-theory commutes with direct limits (of rings) we may replace $Y$ and $Y'$ with finite type, smooth affine $k$-schemes, and we are free to shrink to a smaller neighborhood of $y$ in $Y$ as needed. 

Let $\bar W\subset Y$ be the closure of $p_2(W)$. Note that the condition  $\codim_{X\times Y} W>\codim_{X\times Y'}(x,y)$ implies that $\dim_k W< \dim_kY$, hence $\bar W$ is a proper closed subset of $Y$. Take a  divisor $D\subset Y$ containing $\bar W$. Then there is a morphism
\[
\pi:Y\to \A^n,
\]
$n=\dim_kY-1$, such that $\pi$ is smooth in a neighborhood of $y$ and $\pi:D\to \A^n$ is finite.  Let 
\[
W':= \pi^{-1}(\pi(W)).
\]
Choosing $\pi$ general enough, and noting that
\[
\codim_{X\times Y}W'=\codim_{X\times Y}W-1\ge \codim_{X\times Y'}(x,y)=\dim_k X\times Y',
\]
we may assume that  $W'  \cap X\times Y'$ is a finite set of closed points, say $T$. Let $S\subset D$ be the finite set of closed points $\pi^{-1}(\pi(y))\cap D$.

The inclusion $D\to Y$ induces a section $s:D\to Y\times_{\A^n}D$ to $p_2:Y\times_{\A^n}D\to D$; since $\pi$ is smooth at $y'$, $s(D)$ is contained in the regular locus of $Y\times_{\A^n}D$ and is hence a Cartier divisor on $Y\times_{\A^n}D$. Noting that $p_1:Y\times_{\A^n}D\to Y$ is finite, there is an open neighborhood $U$ of $S$ in $Y$ such that $s(D)\cap Y\times_{\A^n}U$ is a principal divisor; let $t$ be a defining equation. Let  $D_U:=D\cap U$.

This  gives us the commutative diagram
 \[
 \xymatrix{
Y\times_{\A^n}U\ar[r]^-q\ar@<3pt>[d]^p& U\\
D_U\ar@<3pt>[u]^{s}\ar[ur]_-{i}
 }
 \]
 with $q$ finite. Thus we have, for $M\in \sM_{D_U;\sA}$, the exact sequence
\[
0\to q_*(p^*M)\xrightarrow{q_*(\times t)}q_*(p^*M))\to i_*M\to 0
\]
natural in $M$. 

Note that, if $M$ is supported in $W$, then $q_*(p^*M)$ is supported in $W'$. Letting $i':W\to W'$ be the inclusion, our exact sequence gives us the identity
\[
[i'_*M]=0\text{ in }G^{W'}_0(Y;\sA),
\]
hence 
\[
i^*([i'_*M])=0\text{ in }   G^{W'\cap Y'}(Y';\sA).
\]
Let $\bar{i}:(x,y)\to T$ be the inclusion. We have the commutative diagram
\[
\xymatrix{
G_0^{W}(X\times Y;\sA)\ar[r]^-{i'_*}\ar[d]_{i^*}& G_0^{W'}(X\times Y';\sA)\ar[d]^{i^*}\\
G_0^{(x,y)}(X\times Y';\sA)\ar[r]^-{\bar{i}_*},&G_0^{T}(X\times Y';\sA).
}
\]
Since $T$ is a finite set of points containing $(x,y)$,
\[
G_0^{T}(X\times Y';\sA)=G_0^{(x,y)}(X\times Y';\sA)\oplus G_0^{T\setminus \{(x,y)\}}(X\times Y';\sA),
\]
with $\bar{i}_*$ the inclusion of the summand $G_0^{(x,y)}(X\times Y';\sA)$,  from which the result follows directly.
\end{proof}

For a closed immersion $i:W\to X\times Y$, restricting to the generic points of $W$ and using the canonical weak equivalence
\[
G(W;\sA)\to G^W(X\times Y;\sA)
\]
gives the map
\[
\phi_W:G_0^W(X\times Y;\sA)\to  \sK_0^\sA(W).
\]

Each map of pairs $f:(i':W'\to X\times Y')\to (i:W\to X\times Y)$ induces a commutative diagram of inclusions
\[
\xymatrix{
X\times Y'\setminus W'\ar[r]\ar[d]&X\times Y'\ar[d]\\
X\times Y\setminus W\ar[r]&X\times Y;}
\]
Noting that $\id\times f:X\times Y'\to X\times Y$ is an lci morphism, we may apply $G(-)$ to this diagram, giving us the induced map on the homotopy fibers
\[
f^*:G_0^W(X\times Y;\sA)\to G_0^{W'}(X';\sA).
\]

Thus, we have the diagram
\[
\xymatrix{
G_0^W(X\times Y;\sA)\ar[r]^{f^*}\ar[d]_{\phi_W}&G_0^{W'}(X\times Y';\sA)\ar[d]^{\phi_{W'}}\\
 \sK_0^{\sA}(W)& \sK_0^{\sA}(W')
}
\]
In order that $f^*$ descend to  a map 
\[
f^*: \sK_0^\sA(W)\to \sK_0^\sA(W'), 
\]
it therefore suffices to prove:
\begin{lem} \label{lem:restriction} (1) For each $i:W\to X\times Y$, the map $\phi_W$ is surjective. \\
\\
(2) $\phi_{W'}(f^*(\ker\phi_W))=0$.
\end{lem}

\begin{proof} The surjectivity of $\phi_W$ follows from Quillen's localization theorem, which first of all identifies $K_0^W(X\times Y;\sA)$ with $G_0(W;\sA)$ and secondly implies that the restriction map
\[
j^*:G_0(W;\sA)\to G_0(k(W);\sA)=K_0(k(W);\sA)
\]
is surjective. 

For (2), we can factor $f$ as a composition of a closed immersion followed by a smooth morphism. In the second case,  $f^{-1}(W\setminus \Spec k(W))$ is a proper closed subset of $W'$, hence classes supported in $W\setminus \Spec k(W)$ die when pulled back by $f$ and restricted to $k(W')$. Thus we may assume $f$ is a closed immersion. 

Fix a generic point $w'=(x,y)$ of $W'$. We may replace $X$ with $\Spec \sO_{X,x}$ and replace $Y$ with  $\Spec \sO_{Y,y}$. Making a base-change, we may assume that $k(x,y)$ is finite over $k$. Since $X\times_k Y$ is smooth, it follows that 
\[
\codim_{X\times Y}W\ge \codim_{X\times Y'}(x,y).
\]

 Let  $W''\subset W$ is a closed subset of $W$ containing no generic point of $W$. Then 
\[
\codim_{X\times Y}W''> \codim_{X\times Y'}(x,y),
\]
hence by lemma~\ref{lem:QuillenGersten} the restriction map
\[
G_0^{W''}(X\times Y;\sA)\to G_0{(x,y)}(X\times Y';\sA)
\]
is the zero map. By Quillen's localization theorem we have
\[
\ker\phi_W=\colim G_0^{W''}(X\times Y;\sA)
\]
over such $W''$, which proves the lemma.
\end{proof}

\subsection{The cycle complex}
Let $T$ be a finite type $k$-scheme. We let $\dim_kT$ denote the Krull dimension of $T$; we sometimes write $d_T$ for $\dim_kT$.

We fix as above a finite type $k$-scheme $X$ and a sheaf of Azumaya algebras $\sA$ on $X$. We have the cosimplicial scheme $\Delta^*$ with
\[
\Delta^n:=\Spec k[t_0,\ldots, t_n]/\sum_it_i-1
\]
and the standard coface and codegeneracy maps. A \emph{face} $F$ of $\Delta^n$ is a closed subscheme of the form $t_{i_1}=\ldots, =t_{i_s}=0$. We let $\sS^X_r(n)$ be the set of closed subsets $W\subset X\times\Delta^n$ with 
\[
\dim_kW\cap X\times F\le r+\dim_kF
\]
for all faces $F\subset \Delta^n$. We order $\sS^X_r(n)$ by inclusion. If $g:\Delta^m\to \Delta^n$ is the map corresponding to a map $g:[m]\to[n]$ in $\Ord$, and $W$ is in $\sS^X_r(n)$, then $g^{-1}(W)$ is in $\sS_r^X(m)$, so $n\mapsto \sS^X_r(n)$ defines a simplicial set. We let $X_r(n)\subset \sS_r^X(n)$ denote the set of irreducible $W\in \sS^X_r(n)$ with $\dim_kW=r+n$.

\begin{Def} 
\[
z_r(X,n;\sA):= \oplus_{W\in X_r(n)} K_0(k(W);\sA).
\]
\end{Def} 

\begin{rem} Let $W\subset X\times\Delta^n$ be a closed subset. Then restriction to the generic points of $W$ gives the isomorphism
\[
\sK_0^\sA(W\subset X\times\Delta^n)\cong \oplus_{w\in W^{(0)}} K_0(k(w);\sA).
\]
Thus, we can identify $z_r(X,n;\sA)$ with the quotient:
\[
z_r(X,n;\sA)\cong \frac{\colim_{W\in \sS^X_r(n)}\sK_0^\sA(W\subset X\times\Delta^n)}{
\colim_{W'\in \sS^X_{r-1}(n)}\sK_0^\sA(W'\subset X\times\Delta^n)}
\]

Suppose each irreducible 
$W'\in  \sS^X_{r-1}(n)$ is contained in some irreducible $W\in  \sS^X_{r}(n)$ with $\dim_kW=r+n$; as the map
\[
\sK_0^\sA(W')\to \sK_0^\sA(W)
\]
is in this case the zero-map, it follows that 
\[
z_r(X,n;\sA)\cong \colim_{W\in \sS^X_r(n)}\sK_0^\sA(W\subset X\times\Delta^n)
\]
if this condition is satisfied, e.g., for $X$ quasi-projective over $k$.
\end{rem}

 Let $g:\Delta^m\to \Delta^n$ be a map corresponding to $g:[m]\to [n]$ in $\Ord$.  By lemma~\ref{lem:restriction} and the above remark,  we have a well-defined pullback map
 \[
 \id\times g^*: z_r(X,n;\sA)\to z_r(X,m;\sA),
 \]
 giving us the simplicial abelian group $n\mapsto z_r(X,n;\sA)$. We let $(z_r(X,*;\sA),d)$ denote the associated complex, i.e.,
 \[
 d_n:=\sum_{i=0}^n(-1)^i(\id\times\delta^{n-1}_i)^*:z_r(X,n;\sA)\to z_r(X,n-1;\sA).
 \]
 
 \begin{Def} We define the higher Chow groups of dimension $r$ with coefficients in $\sA$ as 
 \[
 \CH_r(X,n;\sA):=H_n(z_r(X,*;\sA)).
 \]
 \end{Def}
 
 \subsection{Elementary properties} The standard elementary properties of the cycle complexes  are also valid with coefficients in $\sA$, if properly interpreted.\\
 \\
  \emph{Projective pushforward}. Let $f:X'\to X$ be a proper morphism. For $Y\in\Sm/k$ and $W\subset X'\times Y$, we have the pushfoward map
 \[
f\times\id_*: G_0^W(X'\times Y, f^*\sA)\to G_0^{f\times\id(W)}(X\times Y;\sA)
\]
Also, if $g:Y'\to Y$ is a morphism in $\Sm/k$, then the diagram
\[
\xymatrix{
X'\times Y'\ar[r]^{f\times\id}\ar[d]_{\id\times g}&X\times Y'\ar[d]^{\id\times g}\\
X'\times Y\ar[r]_{f\times\id}&X\times Y
}
\]
is cartesian and Tor-independent, and the vertical maps are lci morphisms, from which it follows that the diagram (with $W'=(\id\times g)^{-1}(W)$)
\[
\xymatrix{
G_0^{W'}(X'\times Y')\ar[r]^{f\times\id}&G_0^{f\times\id(W')}(X\times Y')\ar[d]^{\id\times g}\\
G_0^W(X'\times Y)\ar[u]^{\id\times g^*}\ar[r]_{f\times\id}&G_0^{f\times\id(W)}(X\times Y)\ar[u]^{\id\times g^*}
}
\]
is commutative.

Thus, then maps $(f\times\id_{\Delta^n})_*$ induce a map of complexes
\[
f_*:z_r(X',*;f^*\sA)\to z_r(X,*;\sA)
\]
with the evident functoriality.\\
\\
\emph{Flat pullback}. Let $f:X'\to X$ be a flat  morphism. For $Y\in\Sm/k$ and $W\subset X\times_kY$, we have the pull-back map
\[
f\times\id^*: G_0^W(X\times Y, \sA)\to G_0^{(f\times\id)^{-1}(W)}(X'\times Y,f^*\sA)
\]
commuting with the pull-back maps $\id\times g^*$ for $g:Y'\to Y$ a map in $\Sm/k$.  Since $f$ is flat, the codimension of $W$ is preserved, hence the pullback maps $f\times\id_{\Delta^n}^*$ induce a map of complexes
\[
f^*:z_r(X,*;\sA)\to z_r(X',*;f^*\sA)
\]
functorially in $f$. \\
\\
\emph{Elementary moving lemmas and homotopy property}. 
\begin{Def} Fix a $Y\in\Sm/k$ and let $\sC$ be a finite set of locally closed subsets of $Y$. Let $X\times Y_r^\sC(n)$ be the set of irreducible dimension $r+n$ closed subsets $W$ of $X\times Y\times\Delta^n$ such that $W$ is in $X\times Y_r(n)$ and for each $C\in \sC$
\[
W\cap X\times C\times\Delta^n \text{ is in } \sS_r^{X\times C}(n).
\]
We have the subcomplex $z_r(X\times Y,*;\sF)_\sC$ of $z_r(X\times Y,*;\sF)$, with
\[
z_r(X\times Y,n;\sF)_\sC=\oplus_{W\in X\times Y_r^\sC(n)}\sK^\sA_0(W).
\]
\end{Def}

Exactly the same proof as for \cite[lemma 2.2]{AlgCyc}, using translation by $\GL_n$, gives the following:
\begin{lem}\label{lem:Translation} Let $\sC$ be a finite set of locally closed subsets of $Y$, with $Y=\A^n$ or $Y=\P^{n-1}$.   Then the inclusion
\[
z_r(X\times Y,*;\sA)_\sC\to z_r(X\times Y,*;\sA)
\]
is a quasi-isomorphism.
\end{lem}

Similarly, we have 
\begin{lem}\label{lem: homotopy} The pull-back map
\[
z_r(X,*;\sA)\to z_{r+1}(X\times\A^1;\sA)
\]
is a quasi-isomorphism.
\end{lem}

\subsection{Localization}  Let $j:U\to X$ be an open immersion with closed complement $i:Z\to X$.  Let $Y$ be in $\Sm/k$. If $W\subset X\times Y$ is an irreducible closed subset supported in $Z\times Y$, then $i\times\id$ induces an isomorphism
\[
i\times\id_*:G_0^W(X\times Y,i^*\sA)\to G_0^W(X\times Y;\sA),
\]
which in turn induces the isomorphism
\[
i_*:\sK_0^{i^*\sA}(W)\to \sK_0^\sA(W)
\]
Similarly, if the generic point of $W$ lives over $U\times Y$, then we have the surjection
\[
j\times\id^*:G_0^W(X\times Y;\sA)\to G_0^{W\cap U\times Y}(U\times Y,j^*\sA)
\]
inducing an isomorphism
\[
j^*:\sK_0^\sA(W)\to \sK_0^{j^*\sA}(W\cap U\times Y)
\]
This yields the termwise exact sequence of complexes
\begin{equation}\label{eqn:loc}
0\to z_r(Z,*,i^*\sA)\xrightarrow{i_*}z_r(X,*;\sA)\xrightarrow{j^*}z_r(U,*,j^*\sA)
\end{equation}
It follows from the main result of \cite{Loc} that
\begin{lem} The inclusion
\[
j^*(z_r(X,*\sA))\subset z_r(U,*,j^*\sA)
\]
is a quasi-isomorphism
\end{lem}
hence
\begin{cor}\label{cor:ZLoc}
The sequence \eqref{eqn:loc} thus determines a canonical distinguished triangle in $D^-(\Ab)$, and we have the long exact \emph{localization sequence}
\begin{multline*}
\ldots\to \CH_r(Z,n,i^*\sA)\xrightarrow{i_*}\CH_r(X,n;\sA)\\\xrightarrow{j^*}\CH_r(U,n,j^*\sA)\to
\CH_r(Z,n-1,i^*\sA)\to\ldots
\end{multline*}
\end{cor}

This in turn yields the \emph{Mayer-Vietoris} distinguished triangle for $X=U\cup V$, $U, V\subset X$ Zariski open subschemes
\begin{multline}\label{eqn:MV}
z_r(X,*;\sA)\to z_r(U,*;\sA_U)\oplus z_r(V,*;\sA_V)\\\to z_r(U\cap V,*;\sA_{U\cap V})\to
z_r(X,*-1;\sA)
\end{multline}

\subsection{Extended functoriality}  Assume now that $X$ is equi-dimensional over $k$ (but not necessarily smooth). A modification of the method of \cite{LevineChowMov}, derived from Chow's moving lemma,  yields a functorial model for the assignment $Y\mapsto z^r(X\times Y,*;\sA)$; as we will not need the functoriality in this paper, we omit a further discussion of this topic.

\subsection{Reduced norm} \label{subsec:RedNorm2} For $X\in\Sch_k$,  $\sA=k$, the complex $z_r(X,*;k)$ is just Bloch's cycle complex $z_r(X,*)$. Indeed, for a field $F$, we have the canonical identification of $K_0(F)$ with $\Z$ by the dimension function, giving the isomorphism  
\[
z_r(X,n;k)=\oplus_{w\in X_{(r)}(n)}K_0(k(w))\cong \oplus_{w\in X_{(r)}(n)}\Z=z_r(X,n).
\]
In addition, if $W\subset X\times\Delta^n$ is an integral closed subscheme of dimension $d$, $i:\Delta^{n-1}\to\Delta^n$ is a codimension one face and if $W$ is not contained in $X\times i(\Delta^{n-1})$, then it follows directly from Serre's intersection multiplicity formula that the image of $(\id\times i)^*([\sO_W])$ in $\oplus_{w\in (X\times\Delta^{n-1})_{(d-1)}}K_0(k(w))$  goes to the pull-back cycle $(\id\times i)^*([W])$ under the isomorphism
\[
\oplus_{w\in (X\times\Delta^{n-1})_{(d-1)}}K_0(k(w))\cong z_{d-1}(X\times\Delta^{n-1}).
\]

Now take $\sA$ to be a sheaf of Azumaya algebras on $X$. The collection of reduced norm maps
\[
\Nrd_{\sA_{k(w)}}:K_0(k(w);\sA)\to K_0(k(w))
\]
thus defines the homomorphism
\[
\Nrd_{X,n ;\sA}:z_r(X,n;\sA)\to z_r(X,n).
\]

\begin{lem}\label{lem:NrdSimp} The maps $\Nrd_{X,n;\sA}$ define a map of simplicial abelian groups
\[
n\mapsto[\Nrd_{X,n;\sA}:z_r(X,n;\sA)\to z_r(X,n)].
\]
\end{lem}

\begin{proof} We note that the maps $\Nrd_{X',n;\sA}$ for $X'\to X$ \'etale define a map of presheaves on $X_\et$. Both $z_r(X,n;\sA)$ and $z_r(X,n)$ are sheaves for the Zariski topology on $X$ and $\Nrd_{X,n;\sA}$ defines a map of sheaves, so we may assume that $X$ is local.  If $X'\to X$ is an \'etale cover, then $z_r(X,n;\sA)\to z_r(X',n;\sA)$ and $z_r(X,n)\to z_r(X',n)$ are injective, so we may replace $X$ with any \'etale cover. Since  $\sA$ is locally a sheaf of matrix algebras on $X_\et$, we may assume that $\sA=M_n(\sO_X)$. In this case, $\Nrd_{X,n;\sA}$ is just the Morita isomorphism;  we thus may extend 
$\Nrd_{X,n;\sA}$ to the Morita isomorphism
\[
\Nrd^{W}:G_0^W(X\times\Delta^n;\sA)\to G_0^W(X\times\Delta^n)
\]
for every $W\in \sS_{(r)}^X(n)$. But the pull-back maps $g^*:z_r(X,n;\sA)\to z_r(X,m;\sA)$ and
$g^*:z_r(X,n)\to z_r(X,m)$ for $g:[m]\to [n]$ in $\Ord$ are defined by lifting elements in $z_r(X,n;\sA)$ (resp. $z_r(X,n)$) to $G_0^W(X\times\Delta^n;\sA)$ (resp.  $G_0^W(X\times\Delta^n)$) for some $W$, applying $(\id\times g)^*$ and mapping to $z_r(X,m;\sA)$ (resp. $z_r(X,m)$). Thus the maps $\Nrd_{X,n;\sA}$ define an isomorphism of simplicial abelian groups, completing the proof.
\end{proof}

Thus we have maps
\begin{align*}
&\Nrd_{X,\sA}:z_r(X,*;\sA)\to z_r(X,*)\\
&\Nrd_{X;\sA}:\CH_r(X,n;\sA)\to \CH_r(X,n)
\end{align*}
The naturality properties of $\Nrd$ show that the maps $\Nrd_{X,\sA}$ are natural with respect to flat pull-back and proper push forward (on the level of complexes).

\section{The spectral sequence}\label{sec:SpecSeq} We are now ready for the first of our main
constructions and results. We begin by discussing the {\em homotopy coniveau tower} associated
to the $G$-theory of sheaf of Azumaya algebras $\sA$ on a scheme $X$. Our main result
(theorem~\ref{thm:main1}) is the identification of the layers in the homotopy coniveau tower
with the Eilenberg-Maclane spectra associated to the twist cycle complex $z_p(X,*;\sA)$. The
proof is exactly the same as for standard $K$-theory $K(X)$  (see \cite{LevineKThyMotCoh,
LevineHC}), except that at one point we need to use an extension of some regularity results
from $K(-)$ to $K(-;\sA)$; this extension is given in Appendix~\ref{AppSec:Reg}.

We then turn to the case $X=\Spec k$, where we have the motivic Postnikov tower for the presheaf $K^\sA$. We show how our computation of the layers in the homotopy coniveau tower for $K^A(X)=K(X;A\otimes_k\sO_X)$, for each $X\in\Sm/k$, lead to a computation of of the layers in the motivic Postnikov tower for $K^A$. This completes the proof of our first main theorem~\ref{Thm:Main1} (see theorem~\ref{thm:Slice}). We conclude this section with a comparison of the reduced norm maps in motivic cohomology and $K$-theory, and some computations of the Atiyah-Hirzebruch spectral sequence in low degrees.

\subsection{The homotopy coniveau filtration} Following \cite{LevineHC} we define
\[
G_{(p)}(X,n;\sA):=\colim_{W\in\sS^X_p(n)}G^W(X\times\Delta^n;\sA)
\]
giving the simplicial spectrum $n\mapsto G_{(p)}(X,n;\sA)$, denoted $G_{(p)}(X,-;\sA)$. Note that, for all $p\ge d_X$, the evident map
\[
G_{(p)}(X,-;\sA)\to G(X\times\Delta^*;\sA)
\]
is an isomorphism.

\begin{rem} In order that $n\mapsto G_{(p)}(X,n;\sA)$ form a simplicial spectrum, one needs to make the $G$-theory with support strictly functorial. This is done by first replacing the categories $\sM_{X\times \Delta^n}(\sA)$ with the full subcategory $\sM_{X\times \Delta^n}(\sA)'$ of $\sA$ modules which are coherent sheaves on $X\times\Delta^n$ and are {\em flat} with respect to all inclusions $X\times F\to X\times\Delta^n$, $F\subset \Delta^n$ a face. Quillen's resolution theorem shows that
\[
K(\sM_{X\times \Delta^n}(\sA)')\to K(\sM_{X\times \Delta^n}(\sA))
\]
is a weak equivalence. One then uses the usual trick of replacing $\sM_{X\times \Delta^n}(\sA)'$ with sequences of objects together with isomorphisms (indexed by the morphisms in $\Ord$) to make the pull-backs strictly functorial. 

A similar construction makes $Y\mapsto G(X\times_kY,\sA)$ strictly functorial on $\Sm/k$; we will use this modification from now on without further mention.
\end{rem}

Since $G(X\times-;\sA)$ is homotopy invariant, the canonical map 
\[
G(X;\sA)\to G_{(d_X)}(X,-;\sA)
\]
 is a weak equivalence (on the total spectrum). This gives us the \emph{homotopy coniveau tower}
\begin{equation}\label{eqn:HCTower}
\ldots\to G_{(p-1)}(X,-;\sA)\to G_{(p)}(X,-;\sA)\to\ldots\to G_{(d_X)}(X,-;\sA)\sim G(X;\sA).
\end{equation}

Setting $G_{(p/p-r)}(X,-;\sA)$ equal to the homotopy cofiber of $G_{(p-r)}(X,-;\sA)\to G_{(p)}(X,-;\sA)$, the tower \eqref{eqn:HCTower} yields the spectral sequence
\begin{equation}\label{eqn:HCSS}
E^{p,q}_2=\pi_{-p-q}(G_{(q/q-1)}(X,-;\sA))\Longrightarrow G_{-p-q}(X;\sA)
\end{equation}

\begin{rems}\label{rem:Connect} 1. Let $T$ be a finite type $k$-scheme, $W\subset T$ a closed subscheme with open complement $j:U\to T$ and $\sA$ a sheaf of Azumaya algebras on $T$.  We have the homotopy fiber sequence
\[
G^W(T;\sA)\to G(T;\sA)\to G(U;j^*\sA)
\]
In addition, the spectra $G(T;\sA)$ and $G(U;j^*\sA)$ are -1 connected, and the restriction map 
\[
j^*:G_0(T;\sA)\to G_0(U;j^*\sA)
\]
is surjective. Thus $G^W(T;\sA)$ is -1 connected, hence the spectra $G_{(p)}(X,n;\sA)$ are -1 connected for all $n$ and $p$. \\
\\
2. Noting that $\sS^X_p(n)=\0$ for $p+n<0$,  the  -1 connectedness of $G_{(p)}(X,n;\sA)$  implies that
\[
\pi_N(G_{(p)}(X,-;\sA))=0
\]
for $N<-p$, i.e., that $G_{(p)}(X,-;\sA)$ is $-p-1$ connected. This in turn implies that $G_{(p/p-r)}(X,-;\sA)$
is $-p-1$ connected for all $r\ge0$, that
 the natural map
\[
G(X;\sA)\to \holim_nG_{(d_X/-n)}(X;\sA)
\]
is a weak equivalence and that the spectral sequence \eqref{eqn:HCSS} is strongly convergent.
\end{rems}

Our main result in this section is
\begin{thm}\label{thm:main1} There is a natural isomorphism
\[
\pi_n(G_{(p/p-1)}(X,-;\sA))\cong\CH_p(X,n;\sA).
\]
\end{thm}

\begin{cor} \label{cor:SS}There is a strongly convergent spectral sequence
\[
E^{p,q}_2=\CH_q(X,-p-q;\sA))\Longrightarrow G_{-p-q}(X;\sA)
\]
\end{cor}

The proof is in three steps: we first define a natural ``cycle map" 
\[
\cyc: \pi_n(G_{(p/p+1)}(X,-;\sA))\to\CH_p(X,n;\sA).
\]
which will define the isomorphism. We then use the localization properties of $G_{(p/p+1)}(X,-;\sA)$ and $\CH_p(X,*;\sA)$ to reduce to the case $X=\Spec F$, $F$ a field, and finally we apply  theorem~\ref{thm:WellConn}  to complete the proof.

\subsection{The cycle map} Let $T$ be a finite type $k$-scheme, $W\subset T$ a closed subscheme with open complement $j:U\to T$ and $\sA$ a sheaf of Azumaya algebras on $T$.  We have the homotopy fiber sequence
\[
G^W(T;\sA)\to G(T;\sA)\to G(U;j^*\sA)
\]
In addition, the spectra $G(T;\sA)$ and $G(U;j^*\sA)$ are -1 connected, and the restriction map 
\[
j^*:G_0(T;\sA)\to G_0(U;j^*\sA)
\]
is surjective. Thus $G^W(T;\sA)$ is -1 connected.  In particular, this implies that the spectra $G_{(p)}(X,n;\sA)$ are all -1 connected. A similar argument shows that the spectra $G_{(p/p-r)}(X,n;\sA)$ are all -1 connected.

As we have seen in remark~\ref{rem:Connect}(1), the  the spectra $G_{(p/p-1)}(X,n;\sA)$ are all -1 connected. Let $\EM(\pi_0G_{(p/p-1)}(X,n\sA))$ denote the Eilenberg-Maclane spectrum with $\pi_0=\pi_0G_{(p/p-1)}(X,n\sA)$ and all other homotopy groups 0. Since $G_{(p/p-1)}(X,n;\sA)$ is -1 connected, we have the  map of spectra
\[
\phi_n:G_{(p/p-1)}(X,n;\sA)\to \EM(\pi_0G_{(p/p-1)}(X,n;\sA))
\]
natural in $n$. Letting $\EM(\pi_0G_{(p/p-1)}(X,-;\sA))$ denote the simplicial spectrum $n\mapsto \EM(\pi_0G_{(p/p-1)}(X,n;\sA))$, this gives us the natural map of simplicial spectra
\[
\phi:G_{(p/p-1)}(X,-;\sA)\to \EM(\pi_0G_{(p/p-1)}(X,-;\sA)).
\]

\begin{lem} There is a natural map
\[
\psi_n:\pi_0(G_{(p/p-1)}(X,n;\sA))\to z_p(X,n;\sA),
\]
which is an isomorphism if $X=\Spec F$, $F$ a field.
\end{lem}

\begin{proof} Let $W\subset X\times\Delta^n$ be a closed subset with generic points $w_1,\ldots, w_r$. We have the evident restriction map
\[
G^W_0(X\times\Delta^n;\sA)=G_0(W;\sA)\to \oplus_iG_0(k(w_i);\sA).
\]
Since $\Z_\sA(W)= \oplus_iG_0(k(w_i);\sA)$, we may define
\[
\psi_n:\pi_0(G_{(p/p-1)}(X,n;\sA))\to z_p(X,n;\sA)
\]
by projecting $\oplus_iG_0(k(w_i);\sA)$ on the factors coming from the generic points of $W\in\sS_p^X(n)$ having dimension $n+r$ over $k$. By lemma~\ref{lem:restriction}, $\psi_n$ is natural in $n$.

Suppose now that $X=\Spec F$, $F$ a field; making a base-change and replacing $p$ with $p-\dim_kX$, we may assume that $F=k$ (note that in this case we may assume $p\le 0$). Thus implies that $X\times\Delta^n\cong\A^n_k$. It is easy to see that, for each $W\in \sS_p^X(n)$, the intersection of $-p$ hypersurfaces of sufficiently high degree, containing $W$, is in $\sS_p^X(n)$ and has pure dimension $p+n$. Thus the closed subsets $W\in 
\sS_p^X(n)$ of pure dimension $p+n$ are cofinal in $\sS_p^X(n)$. 

Identify $ z_p(X,n;\sA)$ with the direction sum $\oplus_wG_0(k(w);\sA)$ as $w$ runs over the generic points of dimension $r+n$ $W\in\sS^X_p(n)$. From the localization sequence, we see that the map
\[
\colim_{W\in\sS_p^X(n)} G_0(W;\sA)\to \oplus_wG_0(k(w);\sA)
\]
is surjective, with kernel the subgroup generated by the image of groups $G_0(W';\sA)$ with $\dim W'<p+n$ and $W'\subset W$ for some $W\in\sS_p^X(n)$. The result thus follows from lemma~\ref{lem:Sherman} below.
\end{proof}

\begin{lem}\label{lem:Sherman}  Suppose that $X=\Spec k$. Let $W'\subset \Delta^n_k$ be a closed subset with $W'\in\sS^q_X(n)$ and $\codim_{\Delta^n}W'>q$. Then  the natural map
\[
G_0(W';\sA)\to \colim_{W\in\sS^q_X(n)} G_0(W;\sA)
\]
is the zero-map.
\end{lem}

\begin{proof} This is a modification of the proof of Sherman \cite{Sherman} that the Gersten complex for $\A^n$ is exact. We may assume that $k$ is infinite. Take a general linear linear projection 
\[
\pi:\Delta^n_k=\A^n_k\to \A^{n-1}_k
\]
and let $W=\pi^{-1}(\pi(W'))$.  Then 
\[
\pi:W'\to \A^{n-1}_k
\]
is finite and $W$ is in $\sS^q_X(n)$. In addition, $\pi$ makes $\A^n$ into  a trivial $\A^1$-bundle over $\A^{n-1}$. Thus the canonical section $s:W'\to W'\times_{\A^{n-1}}\A^n$ makes $W'\times_{\A^{n-1}}\A^n\to W'$ into a trivial line bundle over $W$, hence $s(W')\subset W'\times_{\A^{n-1}}\A^n$ is a principal Cartier divisor. Letting $t$ be a defining equation, we have the functorial exact sequence
\[
0\to p_{2*}p_1^*(M)\xrightarrow{\times t}p_{2*}p_1^*(M)\to i_*(M)\to 0
\]
where $p_1:W'\times_{\A^{n-1}}\A^n\to W'$, $p_2:W'\times_{\A^{n-1}}\A^n\to W\subset \A^n$ are the projections and $i:W'\to W$ is the inclusion. Thus 
\[
i_*:G_0(W';\sA)\to G_0(W;\sA)
\]
is the zero-map, completing the proof.
\end{proof}

We denote the composition $\EM(\psi_n)\circ\phi_n$ by
\[
\cyc_n:G_{(p)}(X,n;\sA)\to \EM(z_p(X,n;\sA))
\]
and the map on the associated simplicial objects by
\[
\cyc:G_{(p)}(X,-;\sA)\to \EM(z_p(X,-;\sA))
\]

\subsection{Localization} Consider an open subscheme $j:U\to X$ with closed complement $i:Z\to X$. We let $\sS_r^{U_X}(n)\subset \sS_r^U(n)$ denote the set of closed subsets $W\subset U\times\Delta^n$ such that
\begin{enumerate}
\item $W$ is in $ \sS_r^U(n)$
\item The closure $\bar W$ of $W$ in $X\times\Delta^n$ is in $\sS^X_r(n)$.
\end{enumerate}
Define the spectrum $G_{(r)}(U_X,n;\sA)$ by
\[
G_{(r)}(U_X,n;j^*\sA):=\colim_{W\in\sS^{U_X}_r(n)}G^W(U\times\Delta^n;j^*\sA)
\]
giving us the simplicial spectrum  $G_{(r)}(U_X,-;\sA)$. The restriction map 
\[
j^*: G_{(r)}(U_X,n;\sA)\to  G_{(r)}(U_X,n;\sA)
\]
 factors through $G_{(r)}(U_X,n;\sA)$, giving us the commutative diagram
\[
\xymatrix{
G_{(r)}(X,-;\sA)\ar[r]^{\hat{j}^*}\ar[rd]_{j^*}&G_{(r)}(U_X,-;\sA)\ar[d]^{\phi}\\
&G_{(r)}(U,-;\sA)
}
\]

\begin{lem} The sequence
\[
G_{(r)}(Z,-;i^*\sA)\xrightarrow{i_*}G_{(r)}(X,-;\sA)\xrightarrow{\hat{j}^*}G_{(r)}(U_X,-;\sA)
\]
is a homotopy fiber sequence.
\end{lem}

\begin{proof} In fact, it follows from the localization theorem for $G(-;\sA)$ that, for each $n$, the sequence
\[
G_{(r)}(Z,n;i^*\sA)\xrightarrow{i_*}G_{(r)}(X,n;\sA)\xrightarrow{\hat{j}^*}G_{(r)}(U_X,n;\sA)
\]
whence the result.
\end{proof}

The localization results of \cite{Loc} yield the following result:
\begin{thm} The map
\[
\phi:G_{(r)}(U_X,-;\sA)\to G_{(r)}(U,-;\sA)
\]
is a weak equivalence.
\end{thm}
Thus, we have
\begin{cor} \label{cor:GLoc} The sequences
\[
G_{(r)}(Z,-;i^*\sA)\xrightarrow{i_*}G_{(r)}(X,-;\sA)\xrightarrow{j^*}G_{(r)}(U,-;\sA)
\]
and
\[
G_{(r/r-s)}(Z,-;i^*\sA)\xrightarrow{i_*}G_{(r/r-s)}(X,-;\sA)\xrightarrow{j^*}G_{(r/r-s)}(U,-;\sA)
\]
are  homotopy fiber sequences.
\end{cor}

In addition, we have
\begin{lem}\label{lem:CycFunct} The diagram
\[
\xymatrix{
G_{(r/r-1)}(Z,-;i^*\sA)\ar[r]^{i_*}\ar[d]_\cyc&G_{(r/r-1)}(X,-;\sA)\ar[r]^{j^*}\ar[d]_\cyc&G_{(r/r-1)}(U,-;\sA)\ar[d]_\cyc\\
\EM(z_r(Z,-;i^*\sA))\ar[r]_{i_*}&\EM(z_r(X,-;\sA))\ar[r]_{j^*}&\EM(z_r(U,-;j^*\sA))
}
\]
defines a map of distinguished triangles in $\SH$.
\end{lem}

\begin{proof} It is clear the maps $\cyc_n$ are functorial with respect to the closed immersion $i$ and the open immersion $j$, hence the diagram
\[
\xymatrix{
\pi_0G_{(r/r-1)}(Z,n;i^*\sA)\ar[r]^{i_*}\ar[d]_{\cyc_n}&\pi_0G_{(r/r-1)}(X,n;\sA)\ar[r]^{\hat{j}^*}\ar[d]_{\cyc_n}&\pi_0G_{(r/r-1)}(U_X,n;\sA)\ar[d]_{\cyc_n}\\
z_r(Z,n;i^*\sA)\ar[r]_{i_*}&z_r(X,n;\sA)\ar[r]_{j^*}&z_r(U_X,n;j^*\sA)
}
\]
commutes for each $n$. Similarly, the diagram
\[
\xymatrix{
\pi_0G_{(r/r-1)}(U_X,n;\sA)\ar[r]^{\phi}\ar[d]_{\cyc_n}&\pi_0G_{(r/r-1)}(U,n;\sA)\ar[d]_{\cyc_n}\\
z_r(U_X,n;\sA)\ar[r]_{\phi}&z_r(U,n;j^*\sA)
}
\]
commutes for each $n$. The result follows directly from this.
\end{proof}

\begin{prop}\label{prop:Reduction} Suppose that the map
\[
\cyc(X):G_{(r/r-1)}(X,-;\sA)\to \EM(z_r(X,-;\sA))
\]
is a weak equivalence for $X=\Spec F$, $F$ a finitely generated field extension of $k$. Then $\cyc(X)$ is a weak equivalence for all $X$ essentially of finite type over $k$.
\end{prop}

\begin{proof} This follows from corollary~\ref{cor:ZLoc}, corollary~\ref{cor:GLoc}, lemma~\ref{lem:CycFunct} and noetherian induction.
\end{proof}

\subsection{The case of a field} We have reduced to the case $X=\Spec k$, where we may apply the method of \cite[\S 6.4]{LevineHC}, as explained in section~\ref{subsec:WellConn}.

Let $K^\sA\in\Spt_{S^1}(k)$ be the presheaf of spectra $X\mapsto K(X;\sA)$. We note that 
\begin{lem}\label{lem:KAProps}\ \\
\begin{enumerate}
\item $K^\sA$ is homotopy invariant and satisfies Nisnevich excision.
\item $K^\sA$ is connected
\item $K^\sA\cong \Omega_T(K^\sA)$.
\end{enumerate}
\end{lem}
We have already seen (1); (2) follows from the weak equivalence $K(-;\sA)\to G(-;\sA)$ on $\Sm/k$ and (3) follows from the projective bundle formula.

In particular, for $Y$ in $\Sm/k$ and integer $q\ge0$, we have the simplicial abelian group $z^q(Y,-;K^\sA)$ and  the cycle map
\[
\cyc_{K^\sA}:s^q(Y,-;K^\sA)\to \EM(z^q(Y,-;K^\sA)).
\]

\begin{lem} Let $Y$ be in $\Sm/k$, $d=\dim_kY$. Fix an integer $q\ge0$ and let $r=d-q$. There is a weak equivalence of simplicial spectra
\[
n\mapsto \phi_n: s^q(Y,n;K^\sA)\to G_{(r/r-1)}(Y,n;\sA)
\]
and an isomorphism of simplicial abelian groups
\[
n\mapsto \psi_n:z^q(Y,n;K^\sA)\to z_r(Y,n;\sA) 
\]
such that the diagram
\[
\xymatrix{
s^q(Y,-;K^\sA)\ar[r]^\phi\ar[d]_{\cyc_{K^\sA}(Y)}& G_{(r/r-1)}(Y,-;\sA)\ar[d]^{\cyc(Y)}\\
\EM(z^q(Y,n;K^\sA))\ar[r]_{\EM(\psi)}&\EM(z_r(Y,-;\sA))}
\]
commutes in $\SH$.
\end{lem}

\begin{proof} We have the natural transformation of functors on $\Sm/k$
\[
K(-;\sA)\to G(-;\sA)
\]
In particular, for $T\in \Sm/k$ and $W\subset T$ a closed subset, we have the map
\[
\phi_{T,W}:K^W(T;\sA)\to G^W(T;\sA)
\]
defining a natural transformation of presheaves of spectra on $\Imm_k$. Applying $\phi_{-,-}$ to the colimit of spectra with supports forming the definition of $s^q(Y,n;K^\sA)$ and $G_{(r/r-1)}(Y,n;\sA)$ gives $ \phi_n$. The map $\psi_n$ is defined similarly, using the maps $\pi_0(\phi_{T,W})$. The compatibility with the cycle maps follows directly from the definitions.
\end{proof}

Thus, to prove that $\cyc(Y):G_{(r/r-1)}(Y,-;\sA)\to \EM(z_r(Y,-;\sA))$ is an isomorphism in $\SH$ for all  $r$ and all $Y\in\Sm/k$ (in particular, for $Y=\Spec k$), it suffices to show that
\[
\cyc_{K^\sA}(Y):s^q(Y,-;K^\sA)\to \EM(z^q(Y,n;K^\sA))
\]
is an isomorphism in $\SH$ for all $q$ and all $Y$. For this, it suffices by theorem~\ref{thm:WellConn} to show that $K^\sA$ is well-connected. 

We have already seen that $K^\sA$ is connected (lemma~\ref{lem:KAProps}(2)). By lemma~\ref{lem:KAProps}(3) we need only show that 
\[
\pi_n(K(\hat\Delta^*_{k(Y)};\sA))=0
\]
for $n\neq0$.

We have shown in \cite[theorem 6.4.1]{LevineHC} that the theory $Y\mapsto K(Y)$ is well-connected, in particular, that $\pi_n(K(\hat\Delta^*_{k(Y)};\sA))=0$
for $n\neq0$ and for $\sA=k$. Using the results of appendix~\ref{AppSec:Reg}, especially proposition~\ref{prop:descent}, the same argument shows  $\pi_n(K(\hat\Delta^*_{k(Y)};\sA))=0$ for $n\neq0$ for arbitrary $\sA$. This completes the proof of theorem~\ref{thm:main1}.

\begin{rem} \label{rem:WC} For later use, we record the fact we have just proved above, namely, that the presheaf $Y\mapsto K(Y;\sA)$ is well-connected.
\end{rem}

\subsection{The slice filtration for an Azumaya algebra} 
By proposition~\ref{prop:AzumayaBiratSheaf}, $\Z_\sA$ is a birational motivic sheaf, hence the cycle complex $z^q(X,*;\sZ_A(q)[2q])$ is defined. 

\begin{thm}\label{thm:HigherChowMotCoh} Let $\sA$ be a central simple algebra over a field $k$. For $X\in\Sm/k$, there is an isomorphism of complexes
\[
z^q(X,*;\sA)\xrightarrow{\phi_{X,\sA}} z^q(X,*;\Z_\sA(q)[2q]),
\]
natural with respect to proper push-forward and flat pull-back.
\end{thm}

\begin{proof} We first define for each $n\ge0$ an isomorphism
\[
\phi_{X,\sA,n}:z^q(X,n;\sA)\cong z^q(X,n;\Z_\sA(q)[2q])
\]
Indeed, by definition
\[
z^q(X,n;\sA)=\oplus_{w\in X^{(q)}(n)}\sK^\sA_0(k(w)).
\]
By remark~\ref{rem:coef}, we have
\[
z^q(X,n;\Z_\sA(q)[2q])=\oplus_{w\in X^{(q)}(n)}\Z_\sA(k(w)).
\]
But $\Z_\sA$ is just $\sK^\sA_0$ considered as a sheaf with transfers, giving us the desired isomorphism. 

This isomorphism $\phi_{X,\sA,n}$ is clearly compatible with proper push-forward and flat pull-back. It thus suffices to show that  the $\phi_{X,\sA,n}$ are compatible with the face maps $X\times\Delta^{n-1}\to X\times\Delta^n$.

Let $k\to k'$ be an extension of fields. 
\begin{align*}
&z^q(X,n;\sA)\to z^q(X_{k'},n;\sA(q)[2q])\\
&z^q(X,n;\Z_\sA(q)[2q])\to z^q(X_{k'},n;\Z_\sA(q)[2q])
\end{align*}
are injective, it suffices to check in case $\sA$ is a matrix algebra. By Morita equivalence, it suffices to check for $\sA=k$.

Recall from proposition~\ref{prop:CycComp} the isomorphism of simplicial abelian groups
\[
n\mapsto[\rho_{X,n}:z^q(X,n;\Z(q)[2q])\to z^q(X,n)
\]
and from \S\ref{subsec:RedNorm2} and lemma~\ref{lem:NrdSimp} the reduced norm map (of simplicial abelian groups)
\[
n\mapsto[\Nrd_{X,n ;\sA}:z^q(X,n;\sA)\to z^q(X,n)].
\]
In case $\sA=k$, the maps $\Nrd_{X,n ;\sA}$ are isomorphisms.   It is easy to check that (for $\sA=k$) the  diagram of isomorphisms
\[
\xymatrix{
z^q(X,n;\sA)\ar[rr]^{\phi_{X,\sA,n}}\ar[rd]_{\Nrd_{X,n ;\sA}}&&z^q(X,n;\Z(q)[2q])\ar[ld]^{\rho_{X,n}}\\
&z^q(X,n)}
\]
commutes. Since both the $\Nrd_{X,n ;\sA}$ and $\rho_{X,n}$ define maps of simplicial abelian groups, it follows that  the  $\phi_{X,\sA,n}$ are  define maps  of  of simplicial abelian groups as well.
\end{proof}

\begin{rem} We have the reduced norm map $\Nrd_\sA:\Z_\sA\to \Z$ (as a map of Nisnevich sheaves with transfers) inducing a reduced norm map $\Nrd_\sA(q):\Z_\sA(q)[2q]\to \Z(q)[2q]$ and thus a map of complexes
\[
\Nrd_\sA(q)_X:z^q(X,*;\Z_\sA(q)[2q])\to z^q(X,*;\Z(q)[2q]).
\]
We have as well the reduced norm map of \S\ref{subsec:RedNorm2}
\[
\Nrd_{X;\sA}:z^q(X,*;\sA)\to z^q(X,*).
\]
We claim that the diagram
\[
\xymatrixcolsep{40pt}
\xymatrix{
z^q(X,*;\sA)\ar[r]^{\phi_{X,\sA}}\ar[d]_{\Nrd_{X;\sA}}&z^q(X,*)\ar[d]^{\phi_{X,k}}\\
z^q(X,*;\Z_\sA(q)[2q])\ar[r]_{\Nrd_\sA(q)_X}&z^q(X,*;\Z(q)[2q])}
\]
commutes. Indeed, on $z^q(X,n;\sA)=\oplus_w\sK^\sA_0(k(w))$, both compositions are just sums of the reduced norm maps
\[
\Nrd:K_0(\sA_{k(w)})\to K_0(k(w))=\Z.
\]
\end{rem}

\begin{cor}\label{cor:HigherChowMotCoh}  Let $A$ be a central simple algebra over a perfect field $k$, $Y\in\Sm/k$. Then there is an isomorphism
\[
\psi_{p,q;\sA}:\CH^q(Y,2q-p;\sA)\to H^p(Y,\Z_\sA(q))
\]
natural with respect to flat pull-back and proper push-forward, and compatible with the respective reduced norm maps.
\end{cor}

\begin{proof} This follows from theorem~\ref{thm:BiratMotCoh} and  theorem~\ref{thm:HigherChowMotCoh}.
\end{proof}

\begin{cor}\label{cor:SSMotCoh} Let $A$ be a central simple algebra over a perfect field $k$, $Y\in\Sm/k$. Then there is a strongly convergent spectral sequence
\[
E_2^{p,q}=H^{p-q}(Y,\Z_\sA(-q))\Longrightarrow K_{-p-q}(Y;\sA).
\]
\end{cor}

\begin{proof} By corollary~\ref{cor:SS}, we have the strongly convergent $E_2$ spectral sequence
\[
E_2^{p,q}=\CH^{-q}(Y,-p-q;\sA)\Longrightarrow K_{-p-q}(Y;\sA).
\]
By corollary~\ref{cor:HigherChowMotCoh} we have the isomorphism
\[
\CH^{-q}(Y,-p-q;\sA)\cong H^{p-q}(Y,\Z_\sA(-q))
\]
yielding the result.
\end{proof}
In fact, we have

\begin{thm} \label{thm:Slice} Let $\sA$ be a central simple algebra over a perfect field $k$.Then there is a natural isomorphism
\[
s_n(K^\sA)\cong \EM(\Z_\sA(n)[2n])
\]
\end{thm}

\begin{proof} We first show
\begin{claim} The presheaf $K^\sA$ is the 0-space of a $T$-spectrum $\sK(\sA)$.
\end{claim}
\begin{proof}[Proof of claim.] Take $Y\in \Sm/k$. For a $\sA\otimes_k\sO_Y$-module $M$, which is a locally free coherent $\sO_Y$ module, we associate the complex of sheaves on $Y\times\P^1$
\[
p_1^*(M)\otimes\sO(-1)\xrightarrow{\times X_0}p_1^*(M).
\]
This gives a functor from the category $\sP_{Y,\sA}$ to the category of perfect complexes of $\sA$-modules on $Y\times\P^1$ with support in $Y\times\infty$. Taking $K$-theory spectra, this gives a natural map
\[
\epsilon_Y:K(Y;\sA)\to K^{Y\times\infty}(Y\times\P^1;\sA)=(\Omega_TK^\sA)(Y)
\]
Using the maps $\epsilon_Y$, we thus build a $T$ $\Omega$ spectrum
\[
\sK(\sA):=((K^\sA, K^\sA,\ldots),\epsilon).
\]
Clearly 
\[
K^\sA=\Omega_T^\infty(\sK(\sA))
\]
proving the claim.
\end{proof}

By the claim and remark~\ref{rem:SliceSusp}, it suffices to prove the theorem for the case $n=0$.  We have the natural map
\[
K^\sA\to s_0(K^\sA)
\]
giving the map of Nisnevich sheaves
\[
\phi:\sK^\sA_0\to \pi_0(s_0(K^\sA)).
\]
But we have the canonical isomorphism (see \S\ref{subsec:0thSlice})
\[
s_0(K^\sA)(Y)\cong K(\hat\Delta^*_{k(Y)},\sA)
\]
Since $K^\sA$ is well-connected, we have
\[
\pi_n(K(\hat\Delta^*_{k(Y)},\sA))=\begin{cases}0&\text{ for }n\neq0\\ K_0(k(Y);\sA)&\text{ for }n=0
\end{cases}
\]
where the map $K_0(k(Y);\sA)\to \pi_0(K(\hat\Delta^*_{k(Y)},\sA))$ is induced by the canonical map
$K(k(Y);\sA)\to K(\hat\Delta^*_{k(Y)},\sA)$. In particular, $\phi$ is an isomorphism on function fields; since both  $\sK^\sA_0$ and  $\pi_0(s_0(K(-,\sA)))$ are birational sheaves, $\phi$ is an isomorphism. This gives us the isomorphism
\[
s_0(K^\sA)\cong \EM(\Z_A),
\]
as desired.
\end{proof}

\subsection{The reduced norm map} 
Let $A$ be a central simple algebra over $k$. We have already mentioned the reduced norm map
\[
\Nrd:K_0(A)\to K_0(k)
\]
in section~\ref{subsec:RedNorm1}; there are in fact reduced norm maps
\[
\Nrd:K_n(A)\to K_n(k)
\]
for $n=0,1,2$. For $n=0,1$, these may defined with the help of a splitting field $L\supset k$ for $A$ and Morita equivalence: Use the composition $A\subset A\otimes_kL\cong M_d(L)$
\[
K_n(A)\to K_n(A_L)\cong K_n(M_d(L))\cong K_n(L).
\]
For $n=0$, the map $K_0(k)\to K_0(L)$ is an isomorphism; one checks that the resulting map $K_0(A)\to K_0(k)$ is the reduced norm we have already defined. For $n=1$, one can take $L$ to be Galois over $k$ (with group say $G$) and use that fact that there is a 1-cocycle $\{\bar{g}_\sigma\}\in Z^1(G;\PGL_d(L))$ such that $A\subset M_d(L)$ is the invariant subalgebra under the $G$ action
\[
(\sigma,m)\mapsto \bar{g}_\sigma {}^\sigma m\bar{g}_\sigma^{-1}
\]
As $\det: K_1(M_d(L))\to K_1(L)=L^\times$ the  isomorphism given by Morita equivalence, one sees that the image of $K_1(A)$ in $L^\times$ lands in the $G$-invariants, i.e., in $k^\times=K_1(k)$. 

For $n=2$, the definition of the reduced norm map (due to Merkurjev-Suslin in the square-free degree case \cite[Th. 7.3]{K2Ref} and to Suslin in general \cite[Cor. 5.7]{suslin}) is more complicated; however, we do have the following result. Let 
$Spl_A$ be the set of field extensions $L/k$ that split $A$. 
\begin{prop} \label{prop:NrdProps}  Let $L\supset k$ be an extension field. \\
\\
1. For $n=0,1,2$, the diagram
\[
\xymatrix{
K_n(A_L)\ar[r]^{\Nrd}\ar[d]_{\Nm_{A_L/A}}&K_n(L)\ar[d]^{\Nm_{L/k}}\\
K_n(A)\ar[r]_{\Nrd}&K_n(k)
}
\]
commutes. Here $\Nm_{A_L/A}:K_n(A_L)\to K_n(A)$ is the map on the $K$-groups induced by the restriction of scalars functor, and similarly for $\Nm_{L/k}$.\\
\\
2. For $n=0,1$, the map
\[
\sum \Nm_{A_L/A}:\oplus_{L\in Spl_A}K_n(A_L)\to K_n(A)
\]
is surjective. If $A$ has square free index, $\sum \Nm_{A_L/A}$ is surjective for $n=2$ as well.
\end{prop}

For a proof of the last statement, see \cite[Th. 5.2]{K2Ref}.

Let $L\supset k$ be a field. Since $\CH^m(L,n;\sA)=0$ for $m>n$, due to dimensional reasons, we have the edge homomorphism
\[
p_{n,L;A}:\CH^n(L,n;A)\to K_n(A_L)
\]
coming from the spectral sequence of corollary~\ref{cor:HigherChowMotCoh}.

Let $L/k$ be a finite field extension. We let 
\[
\Nm_{L/k}:\CH^q(Y_L,p;A)\to \CH^q(Y,p;A)
\]
denote the push-forward map for the finite morphism $Y_L\to Y$.

\begin{lem}\label{lem:TopBottom} Let $L/k$ be a finite field extension, $f:\Spec L\to \Spec k$ the corresponding morphism. Then the diagram
\[
\xymatrix{
\CH^n(L,n;A)\ar[r]^-{p_{n,L;A}}\ar[d]_{\Nm_{L/k}}&K_n(A_L)\ar[d]^{\Nm_{A_L/A}}\\
\CH^n(k,n;A)\ar[r]_-{p_{n,k;A}}&K_n(A)}
\]
commutes.
\end{lem}

\begin{proof} Let $w$ be a closed point of $\Delta^n_L$, not contained in any face.  We have the composition
\begin{multline*}
K_0(L(w);A)\cong K_0^w(\Delta^n_L;A)\cong K_0^w(\Delta^n_L,\partial\Delta^n_L;A)\\
\xrightarrow{\alpha} K_0(\Delta^n_L,\partial\Delta^n_L;A)\cong K_n(A_L)
\end{multline*}
The first isomorphism is via the localization sequence for $K(-;A)$. We have the canonical map
\[
K^w(\Delta^n_L,\partial\Delta^n_L;A)\to K^w(\Delta^n_L;A)
\]
which is a weak equivalence since $w\cap \partial\Delta^n_L=\0$, giving us the second isomorphism. The map $\alpha$ is ``forget supports" and the last isomorphism follows from the homotopy property of $K(-;A)$. Denote this composition by
\[
\beta^w_{n,L;A}:K_0(k(w);A)\to K_n(A_L).
\]
Since $z^n(L,n;A)=\oplus_w K_0(k(w);A)$, where the sum is over all closed points $w\in \Delta^n_L\setminus\partial\Delta^n_L$, the maps $\beta^w_{n,L;A}$ induce
\[
\beta_{n,L;A}:z^n(L,n;A)\to K_n(A_L);
\]
we have as well the canonical surjection 
\[
\gamma_{n,L;A}:z^n(L,n;A)\to \CH^n(L,n;A).
\]
It follows easily from the definitions that the diagram
\[
\xymatrix{
z^n(L,n;A)\ar[dr]_{\beta_{n,L;A}}\ar[r]^{\gamma_{n,L;A}}& \CH^n(L,n;A)\ar[d]^{p_{n,L;A}}\\
&K_n(A_L)
}
\]
commutes. 

On the other hand, it is also a direct consequence of the definitions that, for $x\in \Delta^n_k$ the image of $w$ under $\Delta^n_L\to \Delta^n_k$, we have
\begin{align*}
\Nm_{L/k}\circ\gamma_{n,L;A}=\gamma_{n,k;A}\circ \Nm_{A_{L(w)}/A_{k(x)}}\\
\Nm_{L/k}\circ\beta_{n,L;A}=\beta_{n,k;A}\circ \Nm_{A_{L(w)}/A_{k(x)}}
\end{align*}
whence the result.
\end{proof}

\begin{lem}\label{lem:Gen} For all $n\ge0$, the map
\[
\sum_L\Nm_{L/k}:\oplus_{L\in Spl_A}\CH^n(L,n;A)\to\CH^n(k,n;A)
\]
is surjective.
\end{lem}

\begin{proof} In fact, the map 
\[
\sum_L\Nm_{L/k}:\oplus_{L\in Spl_A}z^n(L,n;A)\to z^n(k,n;A)
\]
is surjective. Indeed, let $x$ be a closed point of $\Delta^n_k\setminus\partial\Delta^n_k$. Then
\[
A_{k(x)}=M_n(D)
\]
for some division algebra $D$ over $k(x)$. Letting $L\subset D$ be a maximal subfield of $D$ containing $k(x)$, $L$ splits $D$, hence $L/k$ splits $A$. Since $L\supset k(x)$, there is a closed point $w\in \Delta^n_L\setminus\partial\Delta^n_L$ lying over $x$ with $L(w)=w$, i.e., $w$ is an $L$-point. 

Since $L$ is a maximal subfield of $D$, the degree of $L$ over $k(x)$ is exactly the index of $\Nrd(K_0(D))\subset K_0(k(x))$. Thus the norm map
\[
\Nm_{L/k(x)}:K_0(A_L)\to K_0(A_{k(x)})
\]
is surjective, i.e. $K_0(A_{k(x)})\cdot x$ is contained in the image of $\Nm_{L/k}(z^n(L,n;A))$. As
\[
z^n(k,n;A)=\oplus_xK_0(A_{k(x)})
\]
with the sum over all closed points $x\in \Delta^n_k\setminus\partial\Delta^n_k$, this proves the lemma.
\end{proof}

Recall from \S\ref{subsec:RedNorm2} the reduced norm map
\[
\Nrd_{X ;\sA}:z^q(Y,*;A)\to z^q(Y,*).
\]

\begin{lem}\label{lem:LH} Let $j:k\hookrightarrow L$ be a finite extension field, $Y\in\Sm/k$. Then the diagram
\[
\xymatrix{
z^q(Y_L,-;A)\ar[r]^-{\Nrd_{Y_L,A}}\ar[d]_{\Nm_{L/k}}&z^q(Y_L,-)\ar[d]^{\Nm_{L/k}}\\
z^q(Y,-;A)\ar[r]_-{\Nrd_{Y,A}}&z^q(Y,-)}
\]
commutes.
\end{lem}

\begin{proof} Take $w\in Y_L^{(q)}(n)$ and let $x\in Y^{(q)}(n)$ be the image of $w$ under $Y_L\times\Delta^n\to Y\times\Delta^n$. It is easy to check that the diagram
\[
\xymatrix{
K_0(A_{k(w)})\ar[r]^{\Nrd}\ar[d]_{\Nm_{A_{k(w)}/A_{k(x)}}}&K_0(k(w))\ar[d]^{\Nm_{k(w)/k(x)}}\\
K_0(A_{k(x)})\ar[r]_{\Nrd}&K_0(k(x))}
\]
commutes, from which the lemma follows.
\end{proof}

\begin{prop}\label{prop:NrdCommute}  For $n=0,1,2$ the diagram
\[
\xymatrix{
\CH^n(k,n;A)\ar[r]^-{p_{n,k;A}}\ar[d]_\Nrd&K_n(A)\ar[d]^{\Nrd}\\
\CH^n(k,n)\ar[r]_-{p_{n,k}}&K_n(k)}
\]
commutes.
\end{prop}

\begin{proof} Let $j:k\hookrightarrow L$ be a finite extension field of $k$. We have the diagram
\[
\xymatrix{
\CH^n(L,n;A_L)\ar[dr]_{\Nm_{L/k}}\ar[rr]^-{p_{n,L;A}}\ar[dd]_\Nrd&&K_n(A_L)\ar'[d][dd]^(.25){\Nrd}\ar[dr]^{\Nm_{L/k}}\\
&\CH^n(k,n;A)\ar[rr]^(.4){p_{n,k;A}}\ar[dd]_(.4)\Nrd&&K_n(A)\ar[dd]^{\Nrd}\\
\CH^n(L,n)\ar'[r][rr]_(.3){p_{n,L}}\ar[dr]_{\Nm_{L/k}}&&K_n(L)\ar[dr]_{\Nm_{L/k}}\\
&\CH^n(k,n)\ar[rr]_-{p_{n,k}}&&K_n(k)}
\]
The left hand square commutes by lemma~\ref{lem:LH}, the right hand square commutes by proposition~\ref{prop:NrdProps}, the top and bottom squares commute by lemma~\ref{lem:TopBottom}. 

Now suppose that $L$ splits $A$. Then, after using the Morita equivalence, the maps $\Nrd$ are identity maps, hence the back square commutes. Thus for  $b\in \CH^n(L,n;A)$,  $a=\Nm_{L/k}(b)\in \CH^n(k,n;A)$, we have
\[
\Nrd(p_{n,k,A}(a))=p_{n,k}(\Nrd(a)).
\]
But by lemma~\ref{lem:Gen},  $\CH^n(k,n;A)$ is generated by elements  $a$ of this form, as $L$ runs over all splitting fields of $A$, proving the result. 
\end{proof}

\subsection{Computations}

\begin{thm}[see also theorem~\protect{\ref{thm:Codim1Van} }]\label{thm:Comp} Let $A$ be a central simple algebra over $k$.\\
\\
1. For $n=0,1$, the edge homomorphism
\[
\CH^n(k,n;A)\xrightarrow{p_{n,k;A}}K_n(A)
\]
is an isomorphism.\\
\\
2. The sequence
\[
0\to \CH^1(k,3;A)\xrightarrow{d_2^{-2,-1}}\CH^2(k,2;A)\xrightarrow{p_{2,k;A}}K_2(A)\to 
 \CH^1(k,2;A)\to0
\]
is exact.  
\end{thm}

\begin{proof} We first note that $\CH^m(k,n;A)=0$ for $m>n$ for dimensional reasons. In addition $z^0(k,-;A)$ is the constant simplicial abelian group $n\mapsto K_0(A)$, hence $\CH^0(k,n;A)=0$ for $n\neq0$.  (1) follows thus from the spectral sequence of corollary~\ref{cor:SS}. 

For (2), the same argument gives the exact sequence
\[
0\to \CH^1(k,3;A)\xrightarrow{d_2^{-2,-1}}\CH^2(k,2;A)\xrightarrow{p_{2,k;A}}K_2(A)\to \CH^1(k,2;A)\to0.
\]
\end{proof}

\subsection{Codimension one} We recall the computation of the codimension one higher Chow groups due to Bloch:

\begin{prop}[Bloch \hbox{\cite[theorem 6.1]{AlgCyc}}] Let $F$ be a field. Then
\[
\CH^1(F,n)=\begin{cases} F^\times&\text{ for }n=1\\0&\text{ for }n\neq1\end{cases}
\]
\end{prop}

Note that $\CH^1(F,0)=0$ for dimensional reasons. To show that $\CH^1(F,n)=0$ for $n>1$, let
$D=\sum_in_iD_i$ be a divisor on $\Delta^n_F$, intersecting each face properly, i.e.,
containing no vertex of $\Delta^n_F$ in its support. Suppose that  $D$ represents an element
$[D]\in \CH^1(F,n)$, that is, $d_n(D)=0$. Using the degeneracy maps to add ``trivial"
components, we may assume that $D\cdot \Delta^{n-1}_j=0$ for all $j$, where $\Delta^{n-1}_j$ is the face
$t_j=0$. 

As $\Delta^n_F\cong\A^n_F$, the divisor $D$ is the divisor of a rational function $f$ on
$\Delta^n_F$. Since $D$ intersects each $\Delta^{n-1}_j$ properly, the restriction $f_j$ of $f$
to $\Delta^{n-1}_j$ is a well-defined rational function on $\Delta^{n-1}_j$; as $D\cdot
\Delta^{n-1}_j=0$, $\Div(f_j)=0$, so $f_j$ is a unit on $\Delta^{n-1}_j$, that is, $f_j=a_j$
for some $a_j\in k^\times$. Since $\Delta^{n-1}_j\cap \Delta^{n-1}_l\neq\0$ for all $j,l$\footnote{This is where we use the hypothesis $n>1$}, all
the $a_j$ are equal, thus $f_j=a\in k^\times$ for all $j$. Dividing $f$ by $a$ we may assume
that $f_j\equiv 1$ for all $j$.

Now let $\sD$ be the divisor of $g:=tf-(1-t)$ on $\Delta^n_F\times\A^1_F$, where
$\A^1_F:=\Spec F[t]$. As the restriction of $g$ to $\Delta^{n-1}_j\times\A^1$ is 1, $\sD$
defines an element $[\sD]\in\CH^1(\A^1_F,n)$ with $i_0^*([\sD])=[D]$, $i_1^*([\sD])=0$. By the
homotopy property, $[D]=0$.

We use essentially the same argument plus Wang's theorem \cite{Wang} to complete theorem~\ref{thm:Comp} as follows:

\begin{thm} \label{thm:Codim1Van}  Let $A$ be a central simple algebra over a field $F$. Suppose $A$ has square-free index. Then $\CH^1(F,n;A)=0$ for $n\neq1$, and the edge homomorphism
\[
\CH^2(k,2;A)\xrightarrow{p_{2,k;A}}K_2(A)
\]
is an isomorphism.
\end{thm}

\begin{proof} We reduce as usual to the case where $\deg A=p$ is prime. As above, the case $n=0$ is trivially true. We mimic the proof for $\CH^1(F,n)$ in case $n>1$.

If $A=M_p(k)$, then $\CH^1(F,n;A)=\CH^1(F,n)$, so there is nothing to prove; we therefore assume that $A$ is a degree $p$ division algebra over $k$. Then $A$ admits a splitting field $k'$ of degree $p$ over $k$; since $\CH^1(F\otimes_kk',n;A)=\CH^1(F\otimes_kk',n)=0$ for $n>1$, a norm argument shows that $\CH^1(F,n;A)$ is $p$-torsion.

We have seen in lemma~\ref{lem:Sherman} that the argument of Sherman \cite[theorem 2.4]{Sherman} for the degeneration of the Quillen spectral sequence for $K(\A^n_F)$ goes through word for word to give the degeneration of the Quillen spectral sequence for $K(\A^n_F;A)$. We will use this fact throughout the remainder of the proof.

Take a divisor $D$ representing a class $[D]\in\CH^1(F,n;A)$. Let $\sM^{(1)}(X;A)$ denote the category of $A\otimes_k\sO_X$-modules which are coherent as $\sO_X$-modules and are supported in codimension at least one on $X$. Since $K_0(\Delta^n_F;A)=K_0(A)$ by the homotopy property, the localization sequence
\[
 K_1(A\otimes_FF(\Delta^n))\xrightarrow{\partial}
 K_0(\sM^{(1)}(\Delta^n;A))\to K_0(\Delta^n;A)\to K_0(A\otimes_FF(\Delta^n))
 \]
 for $K(\Delta^n_F;A)$ gives us an element $f\in K_1(A\otimes_FF(\Delta^n))$ with 
\[
\partial f=D.
\]

In fact, since $D$ intersects each $\Delta^{n-1}_j$ properly, we may find such an $f\in
(A\otimes\sO_{\Delta^n,\partial\Delta^n})^\times$, where
\[
\partial\Delta^n:=\cup_j\Delta^{n-1}_j.
\]

We have the localization sequence
\[
0\to K_1(\Delta^{n-1}_j;A)\to K_1(k(\Delta^{n-1}_j);A)\xrightarrow{\partial}K_0(\sM^{(1)}(\Delta^{n-1}_j;A))\to
\]
By the degeneration of the Quillen spectral sequence on $\Delta^{n-1}_j$, it follows that
\[
K_0(\sM^{(1)}(\Delta^{n-1}_j;A))=\oplus_{w\in (\Delta^{n-1}_j)^{(1)}}K_0(A\otimes_kk(w))
\]
Thus,   the fact that $\Delta^{n-1}_j\cdot D=0$ implies that 
 restriction of $f$ to $f_j\in (A\otimes k(\Delta^{n-1}_j)^\times$ lifts uniquely to $K_1(\Delta^{n-1}_j;A)=K_1(A)$. Note that we have the surjection $A^\times\to K_1(A)$; similarly, we may lift $f$ to an element $\tilde f\in (A\otimes \sO_{\Delta^n\{v\}})^\times$ with restriction $\tilde f_j$ to $\Delta^{n-1}_j$.

Let $\partial^+\Delta^n=\cup_{j=1}^n\Delta^{n-1}_j$. The degeneracy maps give compatible splittings to the inclusions $\Delta^{n-1}_j\to\Delta^n$, so we can modify $f$ so that $\tilde{f}_j=1\in (A\otimes k(\Delta^{n-1}_j)^\times$.

Now let $L:=k(\Delta^{n-1}_0)$ and consider $\tilde f_0\in (A_L)^\times$. Restricting to $\Delta^{n-1}_0\cap\Delta^{n-1}_1$ shows that $f_0=1\in K_1(A_L)$. Furthermore, the reduced norm map 
\[
\Nrd:K_1(A_L)\to K_1(L)=L^\times
\]
is injective \cite{Wang}, and finally, for $a\in (A_L)^\times$, we have
\[
\Nrd(a)=\begin{cases} a^p&\text{ for } a\in L^\times\\ \Nm_{L(a)/L}(a)&\text{ for }a\in A_L^\times\setminus L^\times\end{cases}
\]

Now, $L(\tilde f_0)$ is a subfield of $A_L$ of degree $\le p$ over $L$. But since $A$ is a division algebra and 
$L$ is a pure transcendental extension of $k$, $ A_L$ is still a division algebra, and hence either $L(\tilde{f}_0)=L$ or $L(\tilde{f}_0)$ has degree exactly $p$ over 
$L$. In the former case, $1=\Nrd(\tilde{f_0})=f_0^p$, and since $\tilde{f}_j\equiv1$, it follows that $f_0\equiv 1$ as well.

In case $L(\tilde{f}_0)$ has degree exactly $p$ over $L$, 
then $\Nm_{L(\tilde{f}_0)/L}(\tilde{f}_0)=1$. Let $M\supset L(\tilde{f}_0)$ be the Galois closure of $L(\tilde{f}_0)$ over $L$,  let $M_0\subset M$ be the unique subfield of index $p$, and let $\sigma\in \Gal(M/L)$ be the generator for $\Gal(M/M_0)$. Then $\Nm_{M/M_0}(\tilde{f}_0)=1$, so by Hilbert's theorem 90, there is a $\tilde{g}\in M^\times$ with
\[
\tilde{f}_0=\frac{g^\sigma}{g}.
\]
Looking at the proof of Hilbert's theorem 90, we see that we may take $g$ in the integral closure of
$\sO_{\Delta^{n-1}_0,\partial^+\Delta^{n-1}}$, with $g\equiv1$ over all generic points of $\partial^+\Delta^{n-1}$.

By the Skolem-Noether theorem, there is an element $a_g\in A_{M_0}^\times$ with $g^\sigma=a_g^{-1}ga_g$, i.e.
\[
\tilde{f}_0=a_g^{-1}ga_gg^{-1}.
\]
As above, we may take $a_g$ to be a unit in the integral closure of $A\otimes \sO_{\Delta^{n-1}_0,\partial^+\Delta^{n-1}}$.

Let $\hat{L}:=k(\Delta^n)$, and let $\hat{M}\supset \hat{L}(\tilde{f})$ be the Galois closure of $\hat{L}(\tilde{f})$ over $\hat{L}$. Lift $g$ to $\hat{g}\in \hat{M}$ (or rather, in the integral closure $\hat{R}$ of $\sO_{\Delta^n,\Delta^{n-1}_0}$ in $\hat{M}$), with $\hat{g}\equiv1$ over all generic points of $\partial^+\Delta^n$.  Lift $a_g$ similarly to $\hat{a}_g$. Let $d=[\hat{M}_0: \hat{L}]$. We may replace $\tilde{f}^d$ with
\[
\hat{f}:=\Nm_{A_{\hat{M}_0}/A_{\hat{L}}}(\tilde{f}\hat{a}_g^{-1}\hat{g}\hat{a}_g\hat{g}^{-1})
\]
Then $\hat{f}$ restricts to 1 in $A\otimes k(\Delta^{n-1}_j)$ for all $j$, giving a trivialization of $d\cdot[D]$ in $\CH^1(F,n)$. Since $d$ is prime to $p$, it follows that $[D]=0$ in $\CH^1(F,n;A)$.
\end{proof}

\begin{rem}\label{rem:EdgeHomIso} 
We shall give in corollary~\ref{cor:CohComp} below a second proof of theorem~\ref{thm:Codim1Van} , relying on the Merkurjev-Suslin theorem, by proving that $H^p(k;\Z_A(1))=0$ for $A$ of square-free index over a perfect field $k$ and $p\neq1$.  Via the isomorphism of corollary~\ref{cor:HigherChowMotCoh}
\[
\CH^1(k,n;A)\cong H^{2-n}(k,\Z_A(1))
\]
this shows a second time that $\CH^1(k,n;A)=0$ for $n\neq 1$ in the square-free index case. We do not know if this holds for $A$ of arbitrary index.
\end{rem}

\subsection{A map from $SK_i(A)$ to \'etale cohomology}\label{Sect:4.9}

In this section, we use the \'etale version of
the spectral sequence in the previous section to construct homomorphisms from $SK_i(A)$ to quotients of $H^{i+2}_\et(k,\Q/\Z(i+1))$ for $i=1,2$.


The motivic Postnikov tower for $K^A$
\[
f_{n+1}K^A\to f_nK^A\to \ldots\to f_0K^A=K^A
\]
induces by \'etale sheafification the \'etale version
\[
[f_{n+1}K^A]^\et \to [f_nK^A]^\et\to \ldots\to [f_0K^A]^\et=[K^A]^\et
\]
with layers the \'etale sheafications  $s^\et_nK^A$ of the layers $s_nK^A$  of the original tower. Since $s_nK^A=\EM(\Z_A(n)[2n])$ (theorem~\ref{thm:Slice}), and $\Z_A(n)^\et=\Z(n)^\et$, we have
\[
s_n^\et K^A=\EM(\Z(n)^\et[2n]).
\]
Evaluating at $\Spec k$ and taking the spectral sequence of this tower gives the \'etale motivic Atiyah-Hirzebruch spectral sequence for $A$, with
Bloch-Lichtenbaum numbering:
\[E_2^{p,q}=H^{p-q}_\et(k,\Z(-q))\Rightarrow K_{-p-q}^\et(A).\]
 
Here is part of the corresponding $E_2$-plane:
\[
\begin{matrix}
-2&-1&0&1&2&3&4\\
\hline
0&0&H^0_\et(k,\Z)&0&H^2_\et(k,\Z)\\
0&0&H^1_\et(k,\Z(1))&0&H^3_\et(k,\Z(1))\\
H^0_\et(k,\Z(2))&H^1_\et(k,\Z(2))&H^2_\et(k,\Z(2))&0&H^4_\et(k,\Z(2))\\
H^1_\et(k,\Z(3))&H^2_\et(k,\Z(3))&H^3_\et(k,\Z(3))&0&H^5_\et(k,\Z(3))\\
H^2_\et(k,\Z(4))&H^3_\et(k,\Z(4))&H^4_\et(k,\Z(4))&0&H^6_\et(k,\Z(4))&H^7_\et(k,\Z(4))&H^8_\et(k,\Z(4))
\end{matrix}\]

For $i=1,2$, the composition
\[K_i(A)\to K_i^\et(A)\by{\varepsilon} H^i_\et(k,\Z(i))=K_i(k)\]
coincides with the reduced norm, where $\varepsilon$ is the edge homomorphism of the spectral
sequence and the isomorphism follows from the Beilinson-Lichtenbaum conjecture in weight $i$ (that is, Kummer theory for $i=1$ and the Merkurjev-Suslin theorem for $i=2$). Hence we get an induced map
\[SK_1(A) \to \coker(K_2(k)\simeq H^2_\et(k,\Z(2))\by{d_2} H^5_\et(k,\Z(3))).\]

Note that the map $H^4_\et(k,\Q/\Z(3))\to H^5_\et(k,\Z(3)))$ is an isomorphism, independently of the Beilinson-Lichtenbaum conjecture.

For $SK_2(A)$, we get {\it a priori} a map to the quotient of 
\[\coker(K_3^M(k)\simeq H^3_\et(k,\Z(3))\by{d_2} H^6_\et(k,\Z(4)))\]
by the image of a $d_3$ differential starting from $H^1_\et(k,\Z(2))\simeq K_3(k)_{ind}$. If $k$
contains a separably closed field, this group is divisible, hence its image by the torsion
differential $d_3$ is $0$. Note that we also have an isomorphism
\[H^5_\et(k,\Q/\Z(4))\iso H^6_\et(k,\Z(4))).\]
Here, the isomorphism $K_3^M(k)\simeq H^3_\et(k,\Z(3))$ follows from the Beilinson-Lichtenbaum conjecture in weight $3$; if one does not want to assume it, one gets a slightly more obscure quotient.

To compute $d_2$, we use the fact that this spectral sequence is a module on the
corresponding spectral sequence for $K^\et F$ \cite{Pablo}. The latter is
multiplicative \cite{Pablo} and $d_2$ is obviously $0$ on $K_0(F)$ and $K_1(F)$, hence
on all $K_i^M(F)$. For the one above, we then have
\[d_2(x)=x\cdot d_2(1)\]
where $d_2(1)$ is the image of $1\in K_0(F)$ in $Br(F)$. 

When we pass to the function field $K$ of the Severi-Brauer variety of
$A$, $A$ gets split so $d_2(1)_K=0$. By Amitsur's theorem, $d_2(1)$ is a multiple
$\delta[A]$ of $[A]$. 

In fact, we have $\delta=1$. The computation is very similar to our compution of a related boundary map for the motive of a Severi-Brauer variety (see proposition~\ref{prop:Boundary}) so we will be a little sketchy in our discussion here.

\begin{prop}\label{prop:BoundaryBis} $d_2(1)=[A]$.
\end{prop}

\begin{proof} We begin by noting that by naturality, it suffices to restrict the presheaf $Y\mapsto K(Y;A)$ to the small \'etale site over $k$. Fix a Galois splitting field $L$ over $k$ of $A$ with group $G$. As the field extensions of $L$ are cofinal in $k_\et$, it suffices to consider the functor
\[
F\mapsto K(F;A)
\]
on finite extensions $F$ of $k$ containing $L$; denote this subcategory of $k_\et$ by $k_\et(L)$. 

For such an $F$, $A_F$ is isomorphic to a matrix algebra, say $A_F\cong M_n(F)$, so by Morita equivalence, $K(F;A)$ is weakly equivalent to $K(F)$. Similarly,   $\Z_A=\Z$ on $k_\et(L)$. Since
\[
H^p(F,\Z(n))=0
\]
for $p<n$,  and since $\Z(1)\cong \G_m[-1]$, it follows from our identification of the slices (theorem~\ref{thm:Slice}) 
\[
s_nK^A\cong \EM(\Z_A(n)[2n]),
\]
 that the cofiber $f_{0/2}K^A$ of $f_2K^A\to f_0K^A$ is the same as the presheaf of cofibers of $K^A$ by its 1-connected cover 
 \[
 \tau_{\le 1}K^A:=\cofib[\tau_{\ge 2}K^A\to K^A].
 \]
 Thus, to compute $d_2(1)$, we just need to apply the usual machinery of $G$-cohomology to the fiber sequence
 \[
 \Sigma\EM(\sK_1^A)\to \tau_{\le 1}K^A\to \EM(\sK_0^A),
 \]
 similar to our computation in the proof of proposition~\ref{prop:Boundary}.
 
 Let us choose a cocycle $\sigma\mapsto \bar{g}_\sigma\in \PGL_n(L)$ representing the class of $A$ in $H^1(G,\PGL_n(L))$. Thus, if $g_\sigma\in \GL_n(L)$ is a lifting of $\bar{g}_\sigma$,  we have the action of $G$ on $M_n(L)$
 \[
 \phi_\sigma(m):=g_\sigma\cdot {}^\sigma m \cdot g_\sigma^{-1}
 \]
 where ${}^\sigma m$ is the usual action of $G$ by conjugation of the matrix coefficients.  $A$ is isomorphic to the $G$-invariant $k$-subalgebra of $M_n(L)$. Also,  the coboundary in $H^2_\et(k,\G_m)$ of the class  of $A$  in $H^1(k, \PGL_n)$  is represented by the  2-cocycle $\{c_{\tau,\sigma}\}\in Z^2(G,L^\times)$ defined by
 \[
 c_{\tau,\sigma}\cdot \id_{L^n}=g_\tau{}^\tau g_\sigma g_{\tau\sigma}^{-1}.
 \]

 The ring homomorphism $\phi_\sigma:M_n(L)\to M_n(L)$ induces an exact functor
 \[
 \phi_{\sigma*}:\Mod_{M_n(L)}\to \Mod_{M_n(L)}
 \]
 sending projectives to projectives, hence a natural map $\phi_{\sigma_*}:K(L;A)\to K(L;A)$ and thereby a map $\phi_{\sigma_*}:\tau_{\le1}K(L;A)\to \tau_{\le1}K(L;A)$. To compute $d_2(1)$, we apply the following procedure: lift $1\in K_0(L;A)$ to a representing $M_n(L)$-module $F$. For each $\sigma\in G$, choose an isomorphism $\psi_\sigma:\phi_{\sigma*}(F)\to F$, which gives us a path $\gamma_\sigma$ in the 0-space of $K(L;A)$. The path
 \[
 \gamma(\tau,\sigma):=\gamma_\tau\cdot\phi_{\tau_*}[\gamma_\sigma]\cdot\gamma^{-1}_{\tau\sigma}
 \]
 is a loop in $K(L;A)$, giving an element $c'_{\sigma,\tau}\in K_1(L;A)=L^\times$. This gives us a cocycle $\{c'_{\tau,\sigma}\}\in Z^2(G;L^\times)$, which represents $d_2(1)\in H^3_\et(k,\Z(1))=H^2_\et(k,\G_m)$.
 
 To make the computation concrete, let $F$ be a  left  $M_n(L)$-module. Then the isomorphism of abelian groups $F\to M_n(L)\otimes_{M_n(L)}F$ sending $v$ to $1\otimes v$ identifies $\phi_{\sigma*}(F)$ with the $M_n(L)$ module with underlying abelian group $F$, and with mutliplication
 \[
 m\cdot_\sigma v:={}^{\sigma^{-1}}[g_\sigma^{-1}m  g_\sigma]\cdot v.
 \]
 Under this identification, $\phi_{\sigma*}$ acts by the identity on morphisms.
 
 Take $F=L^n$ with the standard $M_n(L)$-module structure.
 One sees immediately that sending $v$ to $g_\sigma \cdot {}^\sigma v$ gives an $M_n(L)$-module isomorphism $\psi_\sigma: \phi_{\sigma_*}(F)\to F$. The loop $\gamma(\tau,\sigma)$ is thus represented by the automorphism $\psi_\tau\circ\phi_{\tau*}(\psi_\sigma)\circ\psi_{\tau\sigma}^{-1}$:
 \begin{align*}
 \psi_\tau\circ\phi_{\tau*}(\psi_\sigma)\circ\psi_{\tau\sigma}^{-1}(v)&=\psi_\tau\circ\phi_{\sigma*}(\psi_\sigma)({}^{(\tau\sigma)^{-1}}[g_{\tau\sigma}^{-1}v])\\
 &=\psi_\tau(g_\sigma\cdot{}^{\tau^{-1}}[g_{\tau\sigma}^{-1}\cdot v])\\&=g_\tau{}^\tau g_\sigma g_{\tau\sigma}^{-1}
 \end{align*}
Since the Morita equivalence $\Mod_{M_n(L)}\to \Mod_L$ sends multiplication by $c\in L$ on $F$ to multiplication by $c$ on  $L$, we have the explicit representation of $d_2(1)$ by the cocycle 
$\{c_{\tau,\sigma}\}$, completing the computation.
\end{proof} 

\section{The motivic Postnikov tower for a Severi-Brauer variety}\label{sec:SeveriBrauer}
Results of Huber-Kahn \cite{HuberKahn} give a computation of the sheaf $H^0$ of the slices of
$M(X)$ for $X$ any smooth projective variety and show that $H^n$ vanishes for $n>0$. For the
motive of a Severi-Brauer variety $X=\SB(A)$, we are able to show (in case $A$ has prime degree
$\ell$ over $k$) that the negative cohomology vanishes as well. We do this  by comparing with the slices of
the $K$-theory of $X$ and using Adams operations to split the appropriate spectral sequence,
proving our second main result theorem~\ref{Thm:Main2} (see theorem~\ref{thm:SBSlice}).

\subsection{The motivic Postnikov tower for a smooth variety}
\begin{lem}\label{lem:EasyVan} Let $X$ be a smooth projective variety, $M(X)\in \DM^\eff(k)$ the motive of $X$.\\
\\
1. $f_n^\mot M(X)=0$ for $n>\dim_kX$.\\
\\
2. For $\dim_k X\le n\le 0$, $\Omega^n_Tf_n^\mot X$ is represented by $C^\Sus(z_\equi(X,n))$.
\end{lem}
\begin{proof} (1) Since the objects $M(Z)[p]$ are dense in $\DM^\eff(k)$, it suffices to show that 
\[
\Hom_{\DM(k)}(M(Z)(n)[p],M(X))=0
\]
for all $p$ and all $n>\dim_k$. Since $RC_*$ is full, it suffices to show the same vanishing for the morphisms in $\DM^\eff_\gm(k)$; since $\DM^\eff_\gm(k)\to \DM_\gm(k)$ is faithful, it suffices to show the vanishing for the morphisms in $\DM_\gm(k)$. 

 $X$ is smooth and projective, so we have
\[
\Hom_{\DM_\gm(k)}(M(Z)(n)[p],M(X))=\Hom_{\DM_\gm(k)}(M(Z\times X),\Z(d-n)[-p])
\]
where $d=\dim_kX$. But
\[
\Hom_{\DM_\gm(k)}(M(Z\times X),\Z(d-n)[-p])=H^{-p}(Z\times X,\Z(d-n))
\]
which is zero for $d-n<0$.

For (2), it follows from \eqref{eqn:MotLoopIso} that 
\[
\Omega^n_Tf_n^\mot M(X)=f_0^\mot\Omega^n_TM(X)=\Omega^n_TM(X)
\]
By \cite[V, theorem 4.2.2]{FSV}, we have  the natural isomorphism
\[
\Hom_{\DM^\eff_-(k)}(M(Y\times\P^n),M(X)[m]) \cong H^m(C_*^\Sus(z_\equi(X\times\P^n,n))(Y))
\]
and one checks that the projection  
\[
\Hom_{\DM^\eff_-(k)}(M(Y\times\P^n),M(X)[m])\to  \Hom_{\DM^\eff_-(k)}(M(Y)(n)[2n],M(X)[m])
\]
corresponding to the summand $M(Y)(n)[2n]\subset M(Y\times\P^n)$ corresponds to the map induced by the projection
\[
z_\equi(X\times\P^n,n))\to z_\equi(X,n)
\]
This gives us the isomorphism
\[
\Omega^n_TM(X)=R\sHom(\Z(n)[2n], M(X))\cong C_*^\Sus(z_\equi(X,n)).
\]
\end{proof}

\begin{rem} Suppose $k$ admits resolution of singularities. Then for $X\in\Sm/k$, we have
\[
f_n^\mot M(X)=0
\]
for $n>\dim_kX$. In fact, let $U\subset Y$ be an open subscheme of some $Y\in \Sm/k$ such that the complement $Z:=Y\setminus U$ is smooth and of pure codimension $n$ on $Y$. We have the Gysin distinguished triangle
\[
M(U)\to M(Y)\to M(Z)(n)[2n]\to M(U)[1]
\]
Also, if $Y$ is in addition projective, then so is $Z$; since $M(Z)(n)[2n]$ is a  summand of $M(Z\times\P^n)$, we see that $M(U)$ is in the thick subcategory generated by smooth projective varieties of dimension $\le\dim_kY$. Using resolution of singularities, the existence of a compactification of $X$ with strict normal crossing divisor at infinity implies that every $X\in\Sm/k$ of dimension $\le d$ is in the thick subcategory of $\DM_\gm(k)$ generated by smooth projective varieties of dimension $\le d$. Thus lemma~\ref{lem:EasyVan} shows that
\[
f_n^\mot M(X)=0
\]
for $n>\dim_kX$.
\end{rem}

For later use, we make the following explicit computation:

\begin{lem} \label{lem:SliceSeq} Let $Y$ be in $\Sm/k$. Let $X$ be smooth, irreducible and projective of dimension $d$ over $k$. The canonical map $f^\mot_dM(X)\to f^\mot_{d-1}M(X)$ induces  the map (in $D(\Ab)$)
\[
[\Omega^{d-1}_Tf_d^\mot M(X)](Y)\xrightarrow{\alpha} [\Omega^{d-1}_Tf_{d-1}^\mot M(X)](Y)
\]
Then $\alpha$ is isomorphic to the map on Bloch's cycle complexes
\[
p_2^*:z^1(Y,*)\to z^1(X\times Y,*)
\]
induced by the projection $p_2:X\times Y\to Y$.
\end{lem}

\begin{proof} By \eqref{eqn:MotLoopIso}, we have
\[
\Omega^{d-1}_Tf_d^\mot M(X)=f^\mot_1 \Omega^{d-1}_TM(X) =f^\mot_1(\Omega^{d-1}_Tf^\mot_{d-1}M(X))
\]

By lemma~\ref{lem:EasyVan}(2), we have 
\[
\Omega^{d-1}_Tf_{d-1}^\mot M(X)=C_*^\Sus(z_\equi(X,d-1))
\]
hence 
\[
\Omega^{d-1}_Tf^\mot_dM(X)\cong f^\mot_1C_*^\Sus(z_\equi(X,d-1)))
\]
and the  map $\Omega^{d-1}_Tf^\mot_dM(X)\to \Omega^{d-1}_Tf^\mot_{d-1}M(X)$ is just the canonical map
\[
f^\mot_1C_*^\Sus(z_\equi(X,d-1)))\to C_*^\Sus(z_\equi(X,d-1)))
\]

Applying corollary~\ref{cor:HC}, we have isomorphisms in $D(\Ab)$
\[
f^\mot_1C_*^\Sus(z_\equi(X,d-1)))(Y)\cong f^1_\mot(Y,*;C_*^\Sus(z_\equi(X,d-1))))
\]
and (by remark~\ref{rem:HC}(2,3)) the canonical map 
\[
f^\mot_1C_*^\Sus(z_\equi(X,d-1)))\to C_*^\Sus(z_\equi(X,d-1)))
\]
 is isomorphic to
\[
\xymatrix{
f^1_\mot(Y,*;C_*^\Sus(z_\equi(X,d-1)))\ar[r]& f^0_\mot(Y,*;C_*^\Sus(z_\equi(X,d-1)))\ar@{=}[d]\\
&C_*^\Sus(z_\equi(X,d-1))(Y\times\Delta^*).}
\]

Now, for any $T\in\Sm/k$, the inclusion
\[
C_*^\Sus(z_\equi(X,d-1))(T)\subset z^1(T\times X,*)
\]
is a quasi-isomorphism. Thus, if $W\subset T$ is a closed subset, we have the quasi-isomorphism
\[
C_*^\Sus(z_\equi(X,d-1))^W(T)\to \cone(z^1(T\times X,*)\to z^1(T\times X\setminus W\times X,*))[-1]
\]
Now suppose that $W$ has pure codimension 1. By Bloch's localization theorem, we have the quasi-isomorphism
\[
z^0(W,*)\to \cone(z^1(T\times X,*)\to z^1(T\times X\setminus W\times X,*))[-1]
\]
also, $z^0(W)=z^0(W,0)\to z^0(W,*)$ is a quasi-isomorphism. If $\codim_XW>1$, a similar computation shows that $C_*^\Sus(z_\equi(X,d-1))^W(T)$ is acyclic. Applying this to the computation of 
$f^1_\mot(Y,*;C_*^\Sus(z_\equi(X,d-1))))$, we have the  isomorphism in $D(\Ab)$
\[
\phi:z^1(Y,*)\to f^1_\mot(Y,*;C_*^\Sus(z_\equi(X,d-1)))).
\]
Furthermore, the composition
\[
z^1(Y,*)\xrightarrow{\phi} f^1_\mot(Y,*;C_*^\Sus(z_\equi(X,d-1))))\to
f^1_\mot(Y,*;z^1(X\times -,*)) 
\]
is the map
\[
W\subset Y\times\Delta^n\mapsto
X\times W\times\Delta^0\subset X\times Y\times\Delta^n\times\Delta^0
\]

It is then easy to see that the composition
\[
z^1(Y,*)\xrightarrow{\phi} f^1_\mot(Y,*;C_*^\Sus(z_\equi(X,d-1))))\to 
 f^0_\mot(Y,*;C_*^\Sus(z_\equi(X,d-1))))
 \]
 combined with the isomorphism in $D(\Ab)$
 \[
  f^0_\mot(Y,*;C_*^\Sus(z_\equi(X,d-1))))\cong C_*^\Sus(z_\equi(X,d-1))(Y)\cong z^1(X\times Y,*)
  \]
  is just the pull-back
  \[
  p_2^*:z^1(Y,*)\to z^1(X\times Y,*).
  \]
\end{proof}

Let $X$ be in $\Sm/k$. For a presheaf of spectra $E$ on $\Sm/k$, we have the associated presheaf
$\sHom(X,E)$, defined by 
\[
\sHom(X,E)(Y):=E(X\times Y).
\]
Applying $\sHom(X,-)$ to a fibrant model defines the functor
\[
R\sHom(X,-):\SH_{S^1}(k)\to \SH_{S^1}(k).
\]
We use the notation $\sHom^\mot$ and $R\sHom^\mot$ for the analogous operations on $C(\PST(k))$ and on $\DM^\eff(k)$.

The operation $R\sHom(X,-)$ does not  in general commute with the truncation functors $f_n$. However, we do have
\begin{lem} \label{lem:0SliceVan} Take     $m> \dim_kX$. Then for all $E\in \SH_{S^1}(k)$, 
\[
s_0 R\sHom(X,f_mE)\cong 0.
\]
\end{lem}

\begin{proof}
Let $F$ be a presheaf of spectra on $\Sm/k$ which is $\A^1$-homotopy invariant and satisfies Nisnevich excision.  As we have seen in \S\ref{subsec:0thSlice}, we have a natural isomorphism in $\SH$
\[
(s_0F)(X)\cong F(\hat\Delta^*_{k(Y)}).
\]
 
Similarly, for $E$ homotopy invariant and satisfying Nisnevich excision, the spectrum
$\sHom(X,f_mE)(Y):=f_mE(X\times Y)$ is weakly equivalent to the simplicial spectrum $q\mapsto f_mE(X\times Y)(q)$ with
\[
f_mE(X\times Y)(q)=\colim_{W\in\sS_{X\times Y}^{(m)}(q)}E^W(X\times Y\times\Delta^q)
\]
Using the moving lemma described in \cite{LevineChowMov}, we thus have the natural weak equivalence
\[
f_mE(X\times \hat\Delta^p_{k(Y)})(q)\cong 
\colim_{W\in\sS_{X\times \hat\Delta^p_{k(Y)}}^{(m)}(q)_{\sC(p)}}E^W(X\times Y\times\Delta^q)
\]
where $\sC(p)$ is the set $X\times F$, with $F$ a face of $\hat\Delta^p_{k(Y)}$. 

Thus  $s_0\sHom(X,f_mE)(Y)$ is weakly equivalent to the total spectral of the bi-simplicial spectra
\[
(p,q)\mapsto s_0\sHom(X,f_mE)(Y)(p,q)=
 \colim_{W\in\sS_{X\times \hat\Delta^p_{k(Y)}}^{(m)}(q)_{\sC(p)}}E^W(X\times  \hat\Delta^p_{k(Y)}\times\Delta^q).
 \]
 We denote the total spectrum by $s_0\sHom(X,f_mE)(Y)(-,-)$.
 
 Let $s_0\sHom(X,f_mE)(Y)(-,q)$ be the total spectrum of the simplicial spectrum
 \[
p\mapsto s_0\sHom(X,f_mE)(Y)(p,q)=
 \colim_{W\in\sS_{X\times \hat\Delta^p_{k(Y)}}^{(m)}(q)_{\sC(p)}}E^W(X\times  \hat\Delta^p_{k(Y)}\times\Delta^q).
 \]
The face maps
 \[
 \delta^{q*}_i:s_0\sHom(X,f_mE)(Y)(-,q)\to
s_0 \sHom(X,f_mE)(Y)(-,q-1)
 \]
 are weak equivalences for all $i=0,\ldots, q$, $q\ge1$ (see \cite[claim, lemma 5.2.1]{LevineHC}). Thus the canonical map
 \[
 s_0\sHom(X,f_mE)(Y)(-,0)\to
 s_0sHom(X,f_mE)(Y)(-,-)
 \]
 is a weak equivalence.
 
Take $W\in\sS_{X\times \hat\Delta^p_{k(Y)}}^{(m)}(0)_{\sC(p)}$, so $W$ is a closed subset of $X\times\hat\Delta^p_{k(Y)}$ of codimension $\ge m>\dim_kX$, and $W\cap X\times F$ has codimension $\ge m$ on $X\times F$ for all faces $F$ of $\hat\Delta^p_{k(Y)}$. In particular, for each vertex $v$ of $\hat\Delta^p_{k(Y)}$,
\[
\codim_{X\times v}W\cap X\times v>\dim_kX.
\]
Thus $W\cap X\times v=\0$. Since $X$ is proper, the projection of $W$, $p_2(W)\subset \hat\Delta_{k(Y)}^p$, is a closed subset disjoint from all vertices $v$. Since $\hat\Delta_{k(Y)}^p$ is semi-local with closed points the set of vertices, this implies that $p_2(W)=\0$, hence $W=\0$, i.e. 
\[
\sS_{X\times \hat\Delta^p_{k(Y)}}^{(m)}(0)_{\sC(p)}=\{\0\}
\]
and thus $s_0\sHom(X,f_mE)(Y)(-,0)\sim0$. Our description of $s_0\sHom(X,f_mE)(Y)$ as a simplicial spectrum thus yields
\[
s_0\sHom(X,f_mE)(Y)\sim0
\]
for all $Y\in\Sm/k$, completing the proof.
\end{proof}

Thus, for $X\in\Sm/k$, smooth and projective of dimension $d$ over $k$, and for $E\in\SH_{S^1}(k)$, we have the tower in $\SH_{S^1}(k)$
\begin{multline}\label{eqn:Slice2}
0=s_0 R\sHom(X,f_{d+1}E)\to s_0 R\sHom(X,f_dE)\to\ldots\\
\to s_0 R\sHom(X,f_0E)=
s_0 R\sHom(X,E)
\end{multline}
gotten by applying $s_0 R\sHom(X,-)$ to the $T$-Postnikov tower of $E$. Since the functors $s_0$ and $R\sHom(X,-)$  are exact, the $m$th layer in the tower \eqref{eqn:Slice2} is isomorphic to  $s_0R\sHom(X,s_mE)$, $m=0,\ldots, \dim_kX$. Evaluating at some $Y\in\Sm/k$, we have the strongly convergent spectral sequence
\begin{equation}\label{eqn:DoubleSliceSSeq}
E^1_{a,b}=\pi_{a+b}(s_0R\sHom(X,s_aE)(Y))\Longrightarrow \pi_{a+b}(s_0R\sHom(X,E)(Y)).
\end{equation}

\subsection{The case of $K$-theory} We take $E=K$, where $K(Y)$ is the Quillen $K$-theory spectrum of the smooth $k$-scheme $Y$. By \cite[theorem 6.4.2]{LevineHC} we have the natural isomorphism
\[
(s_mK)(Y)\cong \EM(z^m(Y,*))\cong \EM(\Z(m)[2m])(Y)
\]
In addition, we have natural Adams operations $\psi_k$, $k=2, 3,\ldots$ acting on $K$ and on the $T$-Postnikov tower of $K$, with  $\psi_k$ acting on $\pi_*(s_mK)(Y)$ by multiplication by $k^m$ for all $Y\in\Sm/k$ (see \cite{LevineKThyMotCoh}).

Thus we have
\begin{lem} \label{lem:main}Suppose $X$ has dimension $p-1$ over $k$ for some prime $p$. Suppose that for $Y\in\Sm/k$ the differentials in the spectral sequence \eqref{eqn:DoubleSliceSSeq}, for $E=K$, are all zero after localising at $p$. Then  the spectral sequence \eqref{eqn:DoubleSliceSSeq} degenerates at $E_1$.
\end{lem}

\begin{proof} The Adams operations act on the spectral sequence and $\psi_k$ acts by multiplication by 
$k^a$ on $E^r_{a,b}$. Thus the differential $d_r:E^r_{a,b}\to E^r_{a-r,q+r-1}$ is killed by inverting $k^a-k^{a-r}$. Since $d^r_{a,b}=0$ if $a+r\ge p$, we need only invert the numbers
\[
\tau_{a,r}:=\text{g.c.d.}\{k^a(1-k^r)\ |\ k=2,3,\ldots\}
\]
for   $a=0,1,2,\ldots, p-2$, $r=1,\ldots, p-a-1$, which only involve primes $q<p$.
\end{proof}

\subsection{The Chow sheaf}
For a smooth projective variety $X$, we have the Nisnevich sheaf with transfers $\sCH^n(X)$ on $\Sm/k$, this being the sheaf associated to the presheaf
\[
Y\mapsto \CH^n(X\times Y).
\]
It is shown in \cite[theorem 2.2]{HuberKahn} that $\sCH^n(X)$ is a  birational motivic sheaf. We can also label with the relative dimension, defining
\[
\sCH_n(X):=\sCH^{\dim_kX-n}(X).
\]

For our next computation, we need:

\begin{lem} \label{lem:computation} Take  $\sF\in C(\Sh^{tr}_\Nis(k))$ which is homotopy invariant and satisfies Nisnevich excision. Suppose in addition that  $\sF$ is connected.   Then the sheaf   $\sH_0^\Nis(s_0^\mot R\sHom(X,s_n^\mot\sF))$ is the Nisnevich sheaf associated to the presheaf $H_0(s^\mot_0R\sHom(X,s^\mot_n\sF))$ with value at $Y\in\Sm/k$  given by the exactness of
\begin{multline*}
\colim_{\substack{W'\in \sS^{n+1}_{X\times Y}(1)\\ W\in \sS^{n}_{X\times Y}(1)}}H_0(\sF^{W\setminus W'}(X\times Y\times\Delta^1\setminus W'))\\\xrightarrow{i_1^*-i_0^*}
\colim_{\substack{W'\in \sS^{n+1}_{X\times Y}(0)\\ W\in \sS^{n}_{X\times Y}(0)}}H_0(\sF^{W\setminus W'}(X\times Y\setminus W'))\\\to H_0(s^\mot_0R\sHom(X,s^\mot_n\sF))(Y)\to0
\end{multline*}
\end{lem}

\begin{proof}  From proposition~\ref{prop:HC} , $R\sHom(X,s^\mot_n\sF)(Y)=(s^\mot_n\sF)(X\times Y)$ is isomorphic in $D(\Ab)$ to $s^n_\mot(X\times Y,-;\sF)$, the total complex of the simplicial complex
\[
m\mapsto s^n_\mot(X\times Y,m;\sF):=\colim_{\substack{W\in\sS_{X\times Y}^{(n)}(m)\\W'\in\sS_{X\times Y^{(n+1)}(n)}}}\sF^{W\setminus W'}(X\times Y\times\Delta^m\setminus W').
\]
By lemma~\ref{lem:SliceConn}, the spectra $ s^n_\mot(X\times Y,m;\sF)$ are all -1 connected. Thus we have the exact sequence
\[
H_0(s^n(X\times Y,1;\sF))\xrightarrow{i_0^*-i_1^*}H_0(s^n_\mot(X\times Y,0;\sF))\to
H_0(s^n_\mot(X\times Y,-;\sF)).
\]

In any case $R\sHom(X,s^\mot_n\sF)$ is in $\DM^\eff(k)$, hence the homology presheaf 
\[
Y\mapsto 
H_0(R\sHom(X,s^\mot_n\sF)(Y))= H_0(s_n^\mot(X\times Y,-;\sF))
\]
is a homotopy invariant presheaf with transfers. Thus, by \cite[III, corollary 4.18]{FSV}, if $Y$ is local, the restriction map
\begin{equation}\label{eqn:ResMap}
H_0(s_n^\mot(X\times Y,-;\sF))
\to H_0(s_n^\mot(X_{k(Y)},-;\sF))
\end{equation}
is injective. In addition,  $R\sHom(X,s^\mot_n\sF)$ is connected. Indeed, $s^\mot_n\sF$ is connected by proposition~\ref{prop:SliceConn}, and this implies that $R\sHom(X,s^\mot_n\sF)$ is connected.
Thus the restriction map \eqref{eqn:ResMap} is also surjective, hence an isomorphism. 

By theorem~\ref{thm:BiratS0}, $s^\mot_0R\sHom(X,s^\mot_n\sF)$ is also birational, and is connected by proposition~\ref{prop:SliceConn}, hence the same argument shows that
\[
H_0(s^\mot_0R\sHom(X,s^\mot_n\sF)(Y))\to H_0(s_0^\mot R\sHom(X,s^\mot_n\sF)(k(Y)))
\]
is an isomorphism.

We now return to the situation $Y\in\Sm/k$. As in the proof of lemma~\ref{lem:0SliceVan}, $s^\mot_0R\sHom(X,s^\mot_n\sF)(Y)$ is given by evaluating
$R\sHom(X,s^\mot_n\sF)$ on $\hat\Delta^*_{k(Y)}$.  Since $R\sHom(X,s^\mot_n\sF)$ is connected by proposition~\ref{prop:SliceConn}, it follows that we have the exact sequence
\begin{multline*}
H_0(R\sHom(X,s^\mot_n\sF))(\hat\Delta^1_{k(Y)})\xrightarrow{i_0^*-i_1*} H_0(R\sHom(X,s^\mot_n\sF))(\hat\Delta^0_{k(Y)})\\
\to H_0(s^\mot_0R\sHom(X,s^\mot_n\sF)(Y))\to 0.
\end{multline*}
But since $R\sHom(X,s^\mot_n\sF)$ is connected, the restriction map 
\[
H_0(R\sHom(X,s_n^\mot\sF)(\Delta^1_{k(Y)}))\to H_0(R\sHom(X,s^\mot_n\sF)(\hat\Delta^1_{k(Y)}))
\]
is surjective, which shows that 
\[
H_0(R\sHom(X,s^\mot_n\sF)(k(Y))\cong H_0(s_0^\mot R\sHom(X,s^\mot_n\sF)(Y)).
\]
Since the restriction map \eqref{eqn:ResMap} is an isomorphism for $Y$ local, it follows that the canonical map
\[
H_0(R\sHom(X,s^\mot_n\sF)(Y))\to H_0(s_0^\mot R\sHom(X,s^\mot_n\sF)(Y))
\]
is an isomorphism for $Y$ local.

Putting this together with our description above of $H_0(R\sHom(X,s^\mot_n\sF)(Y))$ proves the result.
\end{proof}

\begin{lem} \label{lem:CycSliceComp} Let $X$ be a smooth projective variety of dimension $d$. There is a natural isomorphism
\[
\sH_0^\Nis(s^\mot_0R\sHom(X,\Z(n)[2n]))\cong \sCH^n(X)
\]
\end{lem}

\begin{proof}  Since $\Z$ is a birational motive, we have (remark~\ref{rem:BiratSlice})
\[
\Z(n)[2n]\cong s_n^\mot(\Z(n)[2n]).
\]
We can now use lemma~\ref{lem:computation} to compute $H_0^\Nis(s^\mot_0R\sHom(X,s_n^\mot(\Z(n)[2n])))$. 

By lemma~\ref{lem:CycComp},  for $W\subset Y$ a  closed subvariety of codimension $n$, $Y\in\Sm/k$, there is a natural isomorphism
\[
H_0((\Z(n)[2n])^W(T))=H^{2n}_W(Y,\Z(n))\xrightarrow{\rho_{Y,W,n}} z^n_W(Y).
\]
From this, it  follows from  lemma~\ref{lem:computation}  that  $H^\Nis_0(s^\mot_0R\sHom(X,s_n^\mot(\Z(n)[2n])))$ is just the sheafification of
\[
Y\mapsto \CH^n(X\times Y), 
\]
i.e.,
\[
\sH^\Nis_0(s^\mot_0R\sHom(X,s_n^\mot(\Z(n)[2n])))\cong\sCH^n(X).
\]
\end{proof}

\subsection{The slices of $M(X)$} Before proving our main theorem on the slices of the motive of a Severi-Brauer variety, we first use duality to shift the computation of the $n$th slice to a 0th slice of a related motive. 0th slices are easier to handle, because their cohomology sheaves are birational sheaves.

\begin{lem}\label{lem:duality} Let $X$ be smooth and projective of dimension $d$ over $k$. Then for $0\le n\le d$ there is a natural isomorphism
\[
s_n^\mot M(X)\cong s_0^\mot \left(R\sHom^\mot(X,\Z(d-n))\right)(n)[2d]
\]
\end{lem}

\begin{proof} By \cite{HuberKahn}
\begin{align*}
f_n^\mot M(X)&=\sHom_{\DM^\eff(k)}(\Z(n),M(X))(n)\\
&=\sHom_{\DM^\eff(k)}(\Z(d)[2d],M(X)(d-n)[2d])(n)\\
&=\sHom_{\DM^\eff(k)}(M(X),\Z(d-n))(n)[2d]\\
\end{align*}
In addition, using the isomorphism \eqref{eqn:S1LoopIso}, we have
\[
f_{n-1}^\mot\circ\sHom_{\DM^\eff}(\Z(1),-)=\sHom_{\DM^\eff}(\Z(1),-)\circ f^\mot_n;
\]
this plus Voevodsky's cancellation theorem \cite{voecan} implies 
\[
f^\mot_n(F(1))\cong f^\mot_{n-1}(F)(1)
\]
and similarly for the slices $s_n^\mot$.
Thus
\begin{align*}
s_n^\mot M(X)&=s_n^\mot(f_n^\mot(M(X))\\
&=s_n^\mot\left(\sHom_{\DM^\eff(k)}(M(X),\Z(d-n))(n)[2d]\right)\\
&=s_0^\mot\left(\sHom_{\DM^\eff(k)}(M(X),\Z(d-n))\right)(n)[2d]\\
&=s_0^\mot \left(R\sHom^\mot(X,\Z(d-n))\right)(n)[2d]
\end{align*}
\end{proof}

\begin{thm}\label{thm:SBSlice} Let $X$ be a Severi-Brauer variety of dimension $p-1$, $p$ a
prime, associated to a central simple algebra $\sA$ of degree $p$ over $k$. Then
\begin{enumerate}
\item 
\[
s_n^\mot M(X)\cong \sCH_n(X)(n)[2n]
\]
for $n=0,\ldots, p-1$, $s_n^\mot M(X)=0$ for $n\ge p$.
\item There is a canonical isomorphism
\[
\oplus_{n=0}^{p-1}\sCH^n(X)\cong \oplus_{n=0}^{p-1}\Z_{\sA^{\otimes n}}
\]
\item For $n=0,\ldots, p-1$, we have
\[
\sCH^n(X)\cong \Z_{\sA^{\otimes n}}\cong \begin{cases} \Z_{\sA}&\text{ for } n=1,\ldots, p-1\\
\Z&\text{ for }n=0\end{cases}
\]
\end{enumerate}
\end{thm}

\begin{proof} We first note that $X$ satisfies the conditions of lemma~\ref{lem:main}. If $X=\P^{p-1}$, then the projective bundle formula gives the  weak equivalence
\[
R\sHom(\P^{p-1},f_mK)\cong \oplus_{i=0}^{p-1}f_{m-i}K
\]
from which the degeneration of the spectral sequence at $E_1$ for all $Y\in\Sm/k$ easily follows.  In general,  there is a splitting field $L$ for $\sA$ of degree $p$ over $k$, so $X_L\cong \P^{p-1}_L$, and thus the differentials are all killed by $\times p$.

We recall that $s_n K\cong \EM(\Z(n)[2n])$ \cite[theorem 6.4.2]{LevineHC}. By Quillen's computation of the $K$-theory of Severi-Brauer varieties, 
\[
R\sHom(X,K)\cong \oplus_{n=0}^{p-1}K(-;\sA^{\otimes n}).
\]
Finally, the fact that $K(-;\sA^{\otimes n})$ is well-connected (remark~\ref{rem:WC}) implies
\[
s_0(K(-;\sA^{\otimes n}))=\EM(\Z_{\sA^{\otimes n}}).
\]
Since our spectral sequence degenerates at $E^1$, we therefore have the isomorphism
\[
\oplus_{n=0}^{p-1} \pi_*^\Nis s_0(R\sHom(X, \EM(\Z(n)[2n]))\cong \oplus_{n=0}^{p-1}\pi_*^\Nis\EM(\Z_{\sA^{\otimes n}})
\]

Also, we have $s_0\circ \EM=\EM\circ s_0^\mot$,  and
\begin{align*}
&R\sHom(X, \EM(F))=\EM(R\sHom^\mot(X,F))\\
&\pi_m^\Nis(\EM(F))=\sH^{-m}_\Nis(F)
\end{align*}
for $F\in \DM^\eff(k)$. Thus we see that 
\[
\sH^m_\Nis(s_0^\mot(R\sHom^\mot(X,\Z(n)[2n])))=0
\]
 for $m\neq0$ and 
\begin{equation}\label{eqn:Iso1}
\oplus_{n=0}^{p-1} \sH^0_\Nis\left(s_0^\mot(R\sHom^\mot(X, \Z(n)[2n])\right)\cong \oplus_{n=0}^{p-1}\Z_{\sA^{\otimes n}}.
\end{equation}
In particular, $s_0^\mot(R\sHom^\mot(X, \Z(n)[2n]))$ is concentrated in degree 0. Thus, it follows from lemma~\ref{lem:CycSliceComp} that
\[
s_0^\mot(R\sHom^\mot(X, \Z(n)[2n]))\cong \sCH^n(X)
\]
for $n=0,\ldots, p-1$, which together with \eqref{eqn:Iso1} proves (2).

Together with lemma~\ref{lem:duality}, this gives
\begin{align*}
s_n^\mot M(X)&\cong s_0^\mot(R\sHom^\mot(X, \Z(p-1-n))(n)[2p-2]\\
&\cong \sCH_n(X)(n)[2n]
\end{align*}
proving (1). 

For (3), take a finite Galois splitting field $L/k$ for $\sA$ with Galois group $G$. We have the natural map
\[
\pi^*:\sCH^n(X)\to \sCH^n(X_L)^G\cong \Z
\]
with kernel and cokernel killed by $p$. By (2), $\sCH^n(X)$ is torsion-free. Similarly, we have the inclusion
\[
\pi^*:\Z_{\sA^{\otimes n}}\to (\Z_{\sA^{\otimes n}_L})^G\cong \Z.
\]
We thus have compatible inclusions
\[
\xymatrixcolsep{5pt}
\xymatrix{
\oplus_{n=0}^{p-1}\sCH^n(X)\ar@{}[r]|-{\displaystyle\cong}\ar@{^{(}->}[d]&
\oplus_{n=0}^{p-1}\Z_{\sA^{\otimes n}}\ar@{^{(}->}[d]\\
\oplus_{n=0}^{p-1}\sCH^n(X_L)^G\ar@{}[r]|-{\displaystyle\cong}&
\oplus_{n=0}^{p-1}(\Z_{\sA^{\otimes n}_L})^G
}
\]
Clearly $\sCH^0(X)\cong \Z$. For $y\in Y\in\Sm/k$ the quotient  
\[
(\oplus_{n=0}^{p-1}(\Z_{\sA^{\otimes n}_L})^G)_y/(\oplus_{n=0}^{p-1}\Z_{\sA^{\otimes n}})_y
\]
has order $p^{p-1}$ if $\sA_y$ is not split, and 1 otherwise. Thus, for $n=1,\ldots, p-1$,  $\sCH^n(X)_y\subset \sCH^n(X_L)^G_y=\Z$ has index $p$ if $\sA_y$ is not split and index 1 if $\sA_y$ is split. Thus we can write
\[
\sCH^n(X)\cong \Z_{\sA^{\otimes n}}
\]
for $n=0,\ldots, p-1$, completing the proof.
\end{proof}

\section{Applications}\label{sec:Applications}
In this section, we let $X$ be the Severi-Brauer variety $X:=\SB(A)$ associated to a central simple algebra $A$ of prime degree $\ell$ over $k$. We use our computations of the layers for $M(X)$, together with a duality argument and the 
Beilinson-Lichtenbaum conjecture, to study the reduced norm map
\[
\Nrd:H^p(k,\Z_A(q))\to H^p(k,\Z(q))
\]
and prove the first of our main applications corollary~\ref{Cor:Main1} (see theorem~\ref{thm:CohComp}). Combining these results with our identification of the low-degree $K$-theory of $A$ with the twisted Milnor $K$-theory of $k$ gives us our main result on the vanishing of $SK_2(A)$ for $A$ of square-free degree (corollary~\ref{Cor:Main2}, see also theorem~\ref{thm:SK0}).

\subsection{A spectral sequence for motivic homology}  We have the motivic Postnikov tower for $M(X)$
\begin{multline}\label{eqn:MotPost1}
0=f^\mot_\ell M(X)\to f_{\ell-1}^\mot M(X)\to\ldots\\\to f^\mot_1M(X)\to f^\mot_0M(X)=M(X)
\end{multline}
with slices
\[
s_b^\mot M(X)\cong \Z_{A^{\otimes b+1}}(b)[2b]; \ b=0,\ldots, \ell-1.
\]

Let $\alpha^*:\DM^\eff(k)\to \DM^\eff(k)^\et$ be the change of topologies functor, with right adjoint $\alpha_*:\DM^\eff(k)^\et\to \DM^\eff(k)$. The functors $\alpha^*$ and $\alpha_*$ are exact, and applying $\alpha^*$ to $\Z_A(n)\to \Z(n)$ gives an isomorphism $\alpha^*\Z_A(n)\to \alpha^*\Z(n))$. Thus,
we have the tower
\begin{multline}\label{eqn:MotPost2}
0=\alpha_*\alpha^*f^\mot_\ell M(X)\to \alpha_*\alpha^*f_{p-1}^\mot M(X)\to\ldots\\\to \alpha_*\alpha^*f^\mot_1M(X)\to \alpha_*\alpha^*f^\mot_0M(X)=\alpha_*\alpha^*M(X)
\end{multline}
with slices
\[
\alpha_*\alpha^*s_b^\mot M(X)\cong \alpha_*\alpha^*\Z(b)[2b]; \ b=0,\ldots, \ell-1.
\]

Since $\alpha_*$ is right adjoint to $\alpha^*$, the unit $\eta$ of the adjunction gives the natural transformation of towers $\eta:\eqref{eqn:MotPost1}\to \eqref{eqn:MotPost2}$. Defining $\bar{M}(X)$, $f_n^\mot\bar{M}(X)$ and $\bar{\Z}_{A^{\otimes b+1}}(a)$ by the distinguished triangles

\begin{gather*}
M(X)\to \alpha_*\alpha^*M(X)\to \bar{M}(X)\to M(X)[1]\\
f_n^\mot M(X)\to \alpha_*\alpha^*f_n^\mot M(X)\to f_n^\mot\bar{M}(X)\to f_n^\mot M(X)[1]\\
\Z_{A^{\otimes b+1}}(a)\to \alpha_*\alpha^*\Z(a)\to \bar{\Z}_{A^{\otimes b+1}}(a)\to
\bar{\Z}_{A^{\otimes b+1}}(a)[1]
\end{gather*}
we have the tower 
\begin{multline}\label{eqn:MotPost3}
0=f^\mot_p\bar{M}(X)\to f_{p-1}^\mot \bar{M}(X)\to\ldots\\\to f^\mot_1\bar{M}(X)\to f^\mot_0\bar{M}(X)=\bar{M}(X)
\end{multline}
with slices
\[
s_b^\mot\bar{M}(X)\cong \bar{\Z}_{A^{\otimes b+1}}(b)[2b]; \ b=0,\ldots, p-1.
\]

This last  tower thus gives rise to the strongly convergent spectral sequence
\begin{equation}\label{eqn:SliceSS}
E^{p,q}_2\Longrightarrow \Hom_{\DM^\eff(k)}(\Z(a)[b],\bar{M}(X)(a')[p+q])
\end{equation}
with
\[
E^{p,q}_2=\begin{cases}
\Hom_{\DM^\eff(k)}(\Z(a)[b], \bar{\Z}_{A^{\otimes -q+1}}(a'-q)[p-q])&\text{ for }0\le -q\le \ell-1\\
0&\text{ else.}\end{cases}
\]
\begin{lem} \label{lem:Elem} $\Hom_{\DM_-^\eff(k)}(\Z(r'), \Z_A(r)[q])=0$ for\\
\\
1. $r=0$, $r'>0$ and all $q$\\
2. $r=0=r'$ and $q\neq 0$\\
3.  $r'=0$ and $1\le r<q$
\end{lem}

\begin{proof}   $\Z_A$ is a homotopy invariant Nisnevich sheaf with transfers, so 
\[
\Hom_{\DM_-^\eff(k)}(\Z(r'), \Z_A[q-2r'])=\ker[H^q_\Zar(\P^{r'},\Z_A)\to H^q_\Zar(\P^{r'-1},\Z_A)]
\]
Since $\Z_A$ is a constant sheaf in the Zariski topology 
\[
H^q_\Zar(\P^{r'},\Z_A)=\begin{cases} 0&\text{ for }q\neq0\\ \Z_A(k)&\text{ for }q=0,\end{cases}
\]
the last identity following from the homotopy invariance of $\Z_A$. This proves (1) and (2).

For (3), we have seen that  
\[
\Hom_{\DM^\eff_-(k)}(M(X),\Z_A(r)[2r+n])\cong \CH^r(X;A,n)
\]
for all $n$. Taking $X=\Spec k$, (3)   follows from the fact that $z^r(\Spec k;A,n)=0$ for $n<r$ for dimensional reasons.
\end{proof}

\begin{lem} \label{lem:BL} The Beilinson-Lichtenbaum conjecture for weight $n+1$ implies
\[
\Hom_{\DM^\eff(k)}(\Z(d)[2d],\bar{M}(X)(n+1)[m])=0\text{ for }m\le n+2
\]
and that the sequence
\[
0\to H^{n+3}(X,\Z(n+1))\to H^{n+3}_\et(X,\Z^\et(n+1))\to H^{n+3}_\et(k(X),\Z(n+1))
\]
 is exact.
\end{lem}

\begin{proof} The Beilinson-Lichtenbaum conjecture for weight $n+1$ says that the cohomology
sheaves of $\bar \Z(n+1)$ are $0$ in degree $\le n+2$, hence the natural map
\[
H^m(X,\Z(n+1))\to H^m_\et(X,\Z^\et(n+1))
\]
is an isomorphism for $m\le n+2$ and there is an exact sequence
\[0\to H^{n+3}(X,\Z(n+1))\to H^{n+3}_\et(X,\Z^\et(n+1))\to H^0_\Zar(X,\sH^{n+3}_\et(\Z(n+1)))\]
since the cohomology sheaves of $\Z(n+1)$ vanish in degrees $>n+1$. By the Gersten conjecture for 
$\sH^{n+3}_\et(\Z(n+1)))$, the map
\[
\sH^{n+3}_\et(\Z(n+1)))\to H^{n+3}_\et(k(X),\Z(n+1))
\]
is injective, which gives the exact sequence in the statement of the lemma.

In terms of morphisms in
$\DM^\eff(k)$ and
$\DM^\eff(k)^\et$, this says that the change of topologies map
\[
\Hom_{\DM^\eff(k)}(M(X),\Z(n+1)[m])\to \Hom_{\DM^\eff(k)^\et}(\alpha^*M(X),\alpha^*\Z(n+1)[m])
\]
is an  isomorphism for $m\le n+2$ and an injection for $m=n+3$.

$M(X)$ has dual $M(X)^\vee=M(X)(-d)[-2d]$, hence the \'etale version $\alpha^*M(X)$ has dual 
$\alpha^*M(X)(-d)[-2d]$. In other words, we have natural isomorphisms
\[
\Hom_{\DM^\eff(k)}(M(X),\Z(n+1)[m])\cong \Hom_{\DM^\eff(k)}(\Z(d)[2d], M(X)(n+1)[m])
\]
and
\begin{align*}
\Hom_{\DM^\eff(k)^\et}(\alpha^*M(X),&\alpha^*\Z(n+1)[m])\\&\cong \Hom_{\DM^\eff(k)^\et}(\alpha^*\Z(d)[2d], \alpha^*M(X)(n+1)[m])\\
&\cong \Hom_{\DM^\eff(k)}(\Z(d)[2d], \alpha_*\alpha^*M(X)(n+1)[m]).
\end{align*}
Thus, the natural map $M(X)\to \alpha_*\alpha^*M(X)$ induces an isomorphism
\begin{multline*}
\Hom_{\DM^\eff(k)}(\Z(d)[2d], M(X)(n+1)[m])\\
\to  \Hom_{\DM^\eff(k)}(\Z(d)[2d], \alpha_*\alpha^*M(X)(n+1)[m])
\end{multline*}
for $m\le n+2$ and an injection for $m=n+3$,  hence the lemma.
\end{proof}

\begin{thm}\label{thm:CohComp} Let $A$ be a central simple algebra over
$k$  of prime degree $\ell$. Let $n\ge 0$, and assume that the Beilinson-Lichtenbaum conjecture holds in weights $\le n+1$and  for the prime   $\ell$.\\
\\
1. For  $m<n$, the reduced norm
\[
\Nrd: H^m(k,\Z_A(n))\to H^m(k,\Z(n))
\]
is an isomorphism. \\
\\
2. There is an exact sequence
\begin{multline*}
0\to H^n(k,\Z_A(n))\xrightarrow{\Nrd}H^n(k,\Z(n))\xrightarrow{\partial_n}\\ 
H^{n+3}_\et(k,\Z(n+1))\to H^{n+3}_\et(k(X),\Z(n+1))
\end{multline*}
where $X$ is the Severi-Brauer variety of $A$.
\end{thm}

\begin{proof}  
We proceed by induction on $n$. (1) in case $n=0$ follows from lemma~\ref{lem:Elem}(2). For (2), the reduced norm identifies $H^0(k,\Z_A)$ with $\ell H^0(k,\Z)$ in case $[A]\neq0$, and is an isomorphism if $[A]=0$; since $\ell [A]=0$, (2) follows.

Now assume the result for all $n'<n$. By the Beilinson-Lichtenbaum conjecture in weight $n'$,
the map $\Z(n')\to \alpha_*\alpha^*\Z(n')$ induces an isomorphism on
$\Hom_{\DM^\eff(k)}(\Z,-[m])$ for all $m\le n'+1$ and an injection for $m=n'+2$. Thus
\[
\Hom_{\DM^\eff(k)}(\Z,\bar{\Z}(n')[m])=0\text{ for }m\le n'+1
\]
Similarly, applying (1) and (2) to the distinguished triangle defining $\bar{\Z}_A(n')$, our induction assumption gives
\begin{equation}\label{eqn:Vanishing1}
\Hom_{\DM^\eff}(\Z, \bar{\Z}_A(n')[m])=0\text{ for } m<n'.
\end{equation}
Finally, by lemma~\ref{lem:BL},  the Beilinson-Lichtenbaum conjecture for weight $n+1$ gives
\begin{equation}\label{eqn:Vanishing2}
\Hom_{\DM^\eff(k)}(\Z(d)[2d], \bar{M}(X)(n+1)[m])=0\text{ for }m\le n+2.
\end{equation}

Now consider our spectral sequence \eqref{eqn:SliceSS} with $a=d$, $b=2d-2n-2$, $a'=n+1$, where $d=\dim_kX=\ell-1$.  We have
\begin{multline*}
\Hom(\Z(d)[2d-n-2],\bar{M}(X)(n+1)[p+q])
\\=\Hom(\Z(d)[2d],\bar{M}(X)(n+1)[n+2+p+q])
\end{multline*}
so by \eqref{eqn:Vanishing2}  the spectral sequence converges to 0 for $p+q\le 0$. 

The $E^{p,q}_2$ term is 
\[
E_2^{p,q}=\Hom(\Z, \bar{\Z}_{A^{\otimes -q+1}}(n+1-d-q)[n+2-2d+p-q])
\]
for $0\le -q\le d$ and 0 otherwise. For $0\le -q<d-1$ and $p+q\le0$, we have
\begin{align*}
&n':=n+1-d-q<n\\
&n+2-2d+p-q<n'
\end{align*}
For $-q=d$, $A^{\otimes -q+1}$ is a matrix algebra, hence $\bar{\Z}_{A^{\otimes -q+1}}(N)=\bar{\Z}(N)$. Thus
\[
E_2^{p,-d}=\Hom(\Z, \bar{\Z}(n+1)[n+2-d+p])
\]
If $p+q\le 0$, then $p\le d$, so $n+2-d+p\le n+2$.
Thus our induction hypothesis plus Hilbert's theorem 90 in weight $n+1$ yields
\[
E^{p,q}_2=0\text{ for }  0\le-q\le d,\ -q\neq d-1,\ p+q\le 0.
\]
hence there
is exactly one $E_2$ term that is possibly non-zero, namely
\[
E^{p,1-d}_2=\Hom(\Z, \bar{\Z}_{A^{\otimes d}}(n)[n+1-d+p])
\]
The $d_2$ differential is
\[
E^{p,1-d}_2\xrightarrow{d_2}E^{p+2,-d}_2
\]
and since $p+2-d\le0$, $E^{p+2,-d}_2=0$. Since $E^{*,q}_2=0$ for $q<-d$, there are no higher
differentials coming out of $E^{p,1-d}_2$. Similarly, our induction hypothesis implies that
there are no $d_r$ differentials going to $E^{p,1-d}_r$. Thus $E^{p,1-d}_2=E^{p,1-d}_\infty=0$.

Now take $p+q=0$. The abutment of the spectral sequence is still 0 and there is still only
one possibly non-zero $E_2$ term, 
\[
E^{d-1,1-d}_2=\Hom(\Z, \bar{\Z}_{A}(n)[n]).
\]
 
There is a single possibly non-trivial $d_2$ differential, since 
\[
E^{d+1,-d}_2=\Hom(\Z, \bar{\Z}(n+1)[n+3]).
\]

As above, there are no other non-zero  $d_r$ differentials, hence 
\[
d_2^{d-1,1-d}:E^{d-1,1-d}_2\to E^{d+1,-d}_2
\]
is an injection. Moreover, all $d_r$ differentials abutting to $E^{d+1,-d}_r$ have a source
equal to $0$, hence $E^{d+1,-d}_3=E^{d+1,-d}_\infty$.

Let us collect the information obtained so far:

\begin{itemize}
\item $E_2^{p,q}=0$ for $p+q\le 0$, except possibly $(p,q)=(d-1,1-d)$.
\item The differential $d_2^{d-1,1-d}$ induces an exact sequence
\begin{equation}\label{eqn:Exact}
0\to E^{d-1,1-d}_2\to E^{d+1,-d}_2\to \Hom(\Z(d)[2d],\bar{M}(X)(n+1)[n+3]).
\end{equation}
\end{itemize}

Since $E^{p,1-d}_2=0$, we get that the map
\[
\Hom(\Z, \Z_{A^{\otimes d}}(n)[n+1-d+p])\to \Hom(\Z, \alpha_*\alpha^*\Z_{A^{\otimes d}}(n)[n+1-d+p])
\]
is an isomorphism for $p<d-1$ and an injection for $p=d-1$. Since $\Z_A\cong \Z_{A^{\otimes \ell-1}}$, we have
\begin{align*}
&\Hom(\Z, \Z_{A^{\otimes d}}(n)[n+1-d+p])\cong H^{n+1+p-d}(k,\Z_A(n))\\
&\Hom(\Z, \alpha_*\alpha^*\Z_{A^{\otimes d}}(n)[n+1-d+p])\cong H_\et^{n+1+p-d}(k,\Z(n))
\end{align*}
hence the canonical map
\[
\alpha_A:H^m(k,\Z_A(n))\to H^m_\et(k,\Z(n))
\]
is an isomorphism for $m<n$ and an injection for $m=n$. Since $\alpha_A$ factors as
\[
\xymatrix{
 H^m(k,\Z_A(n))\ar[r]^{\alpha_A}\ar[d]_\Nrd& H^m_\et(k,\Z(n))\\
 H^m(k,\Z(n))\ar[ru]_\alpha
 }
 \]
and $\alpha:H^m(k,\Z(n))\to H^m_\et(k,\Z(n)^\et)$ is an isomorphism for $m\le n$ by the
Beilinson-Lichtenbaum conjecture in weight $n$, it follows that $\Nrd$ is an isomorphism for
$m<n$ and an injection for $m=n$, proving (1) and the injectivity of $\Nrd$ in (2). 

From the distinguished triangles defining $\bar{\Z}_A$ and $\bar{\Z}$, we have exact sequences
\begin{multline*}
\to H^{n-1}(k,\Z_A(n))\to H^{n-1}_\et(k,\Z(n))\to E^{d-2,1-d}_2\\
\to H^n(k,\Z_A(n))\to H^n_\et(k,\Z(n))\to E^{d-1,1-d}_2\to H^{n+1}(k,\Z_A(n))\to
\end{multline*}
and
\[
\to H^{n+3}(k,\Z(n+1))\to H^{n+3}_\et(k,\Z(n+1))\to E^{d+1,-d}_2\to H^{n+4}(k,\Z(n+1))\to
\]
But we have already shown $E^{d-2,1-d}_2=0$. Also, 
\[
 H^{n+1}(k,\Z_A(n))=H^{n+3}(k,\Z(n+1))=H^{n+4}(k,\Z(n+1))=0
 \]
for dimensional reasons and $H^n(k,\Z(n))=H^n_\et(k,\Z(n))$ by Bloch-Kato in weight $n$.
Thus we get an exact sequence
\[0\to H^n(k,\Z_A(n))\by{\Nrd} H^n(k,\Z(n))\to E^{d-1,1-d}_2\to 0\]
and an isomorphism
\[H^{n+3}_\et(k,\Z(n+1))\iso E^{d+1,-d}_2.\]

Putting this together with \eqref{eqn:Exact}, we get the exact sequence
 \begin{multline*}
 0\to H^n(k,\Z_A(n))\xrightarrow{\Nrd}
H^n(k,\Z(n))\xrightarrow{\partial_n}H^{n+3}_\et(k,\Z(n+1))\\
\to\Hom(\Z(d)[2d],\bar{M}(X)(n+1)[n+3]),
 \end{multline*}
 where $\partial_n$ is the map induced by $d_2^{d-1,1-d}$. By comparing the spectral sequence for 
 \[\Hom(\Z(d)[2d], M(X)((n+1)[*]),  \Hom(\Z(d)[2d], \alpha_*\alpha^*M(X)((n+1)[*])
 \]
  and
\[\Hom(\Z(d)[2d], \bar{M}(X)((n+1)[*]),\] we see that $H^{n+3}_\et(k,\Z(n+1))
\to\Hom(\Z(d)[2d],\bar{M}(X)(n+1)[n+3])$ factors through the image of 
\[\Hom(\Z(d)[2d],\alpha_*\alpha^*M(X)(n+1)[n+3])\to \Hom(\Z(d)[2d],\bar{M}(X)(n+1)[n+3]).\]

By the exact sequence of lemma~\ref{lem:BL}, we thus have the exact sequence
 \begin{multline*}
 0\to H^n(k,\Z_A(n))\xrightarrow{\Nrd}
H^n(k,\Z(n))\xrightarrow{\partial_n}H^{n+3}_\et(k,\Z(n+1))\\
\to H^{n+3}_\et(k(X),\Z(n+1)).
 \end{multline*}
The resulting map
\[H^{n+3}_\et(k,\Z(n+1))\to H^{n+3}_\et(k(X),\Z(n+1))\]
is induced by an edge homomorphism of our spectral sequence, hence equals the extension of
scalars map. This completes the proof.
 \end{proof}

 \begin{cor} \label{cor:CohComp} Let $A$ be a central simple algebra of square-free index over $k$.  For $n\neq 1$, $H^n(k,\Z_A(1))=0$.
 \end{cor}
 
Of course, we have already proved this by a direct argument (theorem~\ref{thm:Codim1Van} ).
This second argument uses our main result on the reduced norm, theorem~\ref{thm:CohComp},
which, in the weight one case, relies on the Merkurjev-Suslin theorem to prove 
Beilinson-Lichtenbaum in weight two (using in turn \cite{SuslinVoev} or \cite{GeisserLev}).
 
 \begin{proof}  We first reduce to the case of $A$ of prime degree $\ell$. Write
\[
\Deg(A)=\prod \ell_i=d.
\]
where the $\ell_i$ are distinct primes. Write $A=M_n(D)$ for some division algebra $D$ of
degree $d$ over $k$, and let $F\subset D$ be a maximal subfield. Then $F$ has degree $d$ over
$k$ and splits $D$. Let $\ell=\ell_i$ for some $i$, let $k(\ell)\supset k$ be the maximal prime
to $\ell$ extension of $k$ and let $F(\ell):=Fk(\ell)$. Then clearly $F(\ell)$ has degree
$\ell$ over $k(\ell)$ and splits $A_{k(\ell)}$; since $k(\ell)$ has no prime to  $\ell$
extensions, $F(\ell)$ is Galois over $k(\ell)$.  Passing
from $k$ to the $\Gal(k(\ell)/k)$ invariants  alters only the prime to $\ell$ torsion. Thus we
may replace $k$ with $k(\ell)$ and assume that $A$ is split by a degree $\ell$ Galois extension
of $k$. But then $A$ is Morita equivalent to an algebra of degree $\ell$, which achieves the
reduction.
 
 It follows from \cite[theorem 6.1]{AlgCyc} that
 \[
 0=\CH^1(k,2-n)\cong  H^n(k,\Z(1))
 \]
 for $n\neq1$. By theorem~\ref{thm:CohComp}(1), this implies that $H^n(k,\Z_A(1))=0$ for
$n<1$. Additionally,  we have
 \[
 H^n(k,\Z_A(1))\cong \CH^1(k,2-n;A)
 \]
 by corollary~\ref{cor:HigherChowMotCoh}. Since $\CH^1(k,m;A)=0$ for $m<0$ and  $\CH^1(k,0;A)=0$ for dimensional reasons, the proof is complete.
 \end{proof}
 
  \begin{cor} \label{cor:K2Iso} Let $A$ be a central simple algebra of square-free index over $k$. 
  The the edge homomorphism
  \[
  p_{2,k;A}:\CH^2(k,2;A)\to K_2(A)
  \]
  is an isomorphism
  \end{cor}
 
 \begin{proof} From corollary~\ref{cor:CohComp}, $\CH^1(k,n;A)=0$ for $n\neq1$. From
theorem~\ref{thm:Comp}(2), we have the exact sequence
\[
0\to \CH^1(k,3;A)\xrightarrow{d_2^{-2,-1}}\CH^2(k,2;A)\xrightarrow{p_{2,k;A}}K_2(A)\to
\CH^1(k,2;A)\to  0,
\]
hence the edge-homomorphism $p_{2,k;A}:\CH^2(k,2;A)\to K_2(A)$ is an isomorphism. 
\end{proof}

Finally, here is a global version of theorem~\ref{thm:CohComp}:

\begin{cor}\label{cor:CohCompGlobal} Let $\tilde \Z_A$ denote the cokernel of the reduced norm
map $\Nrd:\Z_A\to \Z$. Suppose that $A$ has square-free index and assume the
Beilinson-Lichtenbaum conjecture. Then,
\begin{enumerate}
\item For all $n\ge 0$, the complex $\tilde \Z_A(n)\in
DM^\eff(k)$ is concentrated in degree $n$.
\item Let $\sF_n=\sH^n(\tilde \Z_A(n))$.
Then the stalk of $\sF_n$ at a function field $K$ is isomorphic to
\[\ker(H^{n+3}_\et(K,\Z(n+1)\to H^{n+3}_\et(K(X),\Z(n+1)))\]
where $X$ is the Severi-Brauer variety of $A$. 
\item For any smooth scheme $U$ we have a Gersten
resolution
\[0\to \sF_n\to \bigoplus_{x\in U^{(0)}}(i_x)_*(\sF_n)\to\bigoplus_{x\in
U^{(1)}}(i_x)_*(\sF_{n-1})\to\dots\to\bigoplus_{x\in U^{(p)}}(i_x)_*(\sF_{n-p})\to \dots
\]
\end{enumerate}
\end{cor}

\begin{proof} As in the proof of corollary~\ref{cor:CohComp}, it suffices to handle the case of $A$ of prime degree over $k$.

Clearly, $\Z_A(n)$ has no cohomology in degrees $>n$; by Voevodsky's form of
Gersten's conjecture, the vanishing of $\sH^i(\tilde \Z_A(n))$ for $i<n$ reduces to theorem~\ref{thm:CohComp}. The computation of the stalks of $\sH^n(\tilde \Z_A(n))$ also follows from
theorem~\ref{thm:CohComp}.

For (3), we first show (with the notation of
\cite[3.1]{Voetrans}, that the Zariski sheaf associated to the presheaf $(\sF_n)_{-1}$ is
$\sF_{n-1}$. This follows immediately from Voevodsky's cancellation
theorem \cite{voecan}: by definition
\begin{multline*}
(\sF_n)_{-1}(U) = \coker\left(\sF_n(U\times \A^1)\to \sF_n(U\times (\A^1-\{0\})\right)\\
=\coker\left(H^n(U\times \A^1,\tilde \Z_A(n))\to H^n(U\times (\A^1-\{0\}),\tilde
\Z_A(n))\right). 
\end{multline*}
By purity, the localization sequence for $U\times (\A^1-\{0\})\subset U\times\A^1$, and part (1) of the corollary,   the latter cokernel is isomorphic to
\[\ker(H^{n-1}(U,\tilde \Z_A(n-1))\to H^{n+1}(U,\tilde \Z_A(n))\simeq H^1_\Zar(U,\sF_n)\]
hence the Zariski sheaf associated to $(\sF_n)_{-1}$ is the sheaf associated to
\[U\mapsto H^{n-1}(U,\tilde\Z_A(n-1))\simeq\sF_{n-1}(U).\]
The statement on the Gersten complex follows from this and {\it loc. cit.} theorem 4.37.
\end{proof}

 \subsection{Computing the boundary map} To finish our study of $H^n(k,\Z_A(n))$, we need to compute the boundary map $\partial_n$ in theorem~\ref{thm:CohComp}. As above, we fix a central simple algebra $A$ over $k$ of prime degree $\ell$, let $d=\ell-1$ and let $X$ be the Severi-Brauer variety $\SB(A)$. We let $[A]\in H^2_\et(k,\G_m)$ denote the class of $A$ in the (cohomological) Brauer group of $k$.
 
 Concentrating on $f^\mot_{d-1}M(X)$ gives us the distinguished triangle
 \[
 s^\mot_dM(X)\to f^\mot_{d-1}M(X)\to s^\mot_{d-1}M(X)\to  s^\mot_dM(X)[1]
 \]
 which by theorem~\ref{thm:SBSlice}  is
 \[
 \Z(d)[2d]\to f^\mot_{d-1}M(X)\to \Z_A(d-1)[2d-2]\to \Z(d)[2d+1]
 \]
 Applying $\Omega^{d-1}_T$ gives
  \[
 \Z(1)[2]\to \Omega^{d-1}_Tf^\mot_{d-1}M(X)\to \Z_A \to \Z(1)[3]
 \]
 
 Applying the \'etale sheafification $\alpha^*$ and noting that $\Z_A^\et\cong \Z^\et$ gives the distinguished triangle
  \begin{equation}\label{eqn:EtSeq}
 \Z(1)^\et[2]\to \alpha^*\Omega^{d-1}_Tf^\mot_{d-1}M(X)\to \Z^\et \xrightarrow{\partial} \Z(1)^\et[3]
 \end{equation}
 Thus $\partial:\Z^\et\to \Z(1)^\et[3]$ gives us the element
 \[
 \beta_A\in H^3_\et(k,\Z(1)^\et)=H^2_\et(k,\G_m).
 \]
 \begin{prop}\label{prop:Boundary} $\beta_A=[A]$.
 \end{prop}
 
 \begin{proof} To calculate $\beta_A$, it suffices to restrict \eqref{eqn:EtSeq} to the small \'etale site on $k$.    By lemma~\ref{lem:SliceSeq}, \eqref{eqn:EtSeq} on $k_\et$ is  isomorphic (in $D(\Sh_\et(k))$) to the sheafification of the sequence of presheaves
 \begin{equation}\label{eqn:ConeSeq}
 L\mapsto  \left(z^1(L,*)\xrightarrow{p^*} z^1(X_L,*)\to \cone(p^*)\to z^1(L,*)[1]\right).
\end{equation}
Here, and in the remainder of this proof, we consider the cycle complexes as cohomological complexes: 
\[
z^1(Y,*)^n:=z^1(Y,-n).
\]

We recall that $z^1(X_L,*)$ has non-zero cohomology only in degrees 0 and -1, and that
\begin{align*}
&H^{-1}(z^1(X_L,*))=\Gamma(X_L,\sO_{X_L}^\times),\\
&H^0(z^1(X_L,*))=\CH^1(X_L).
\end{align*}
Similarly, $H^{-1}(z^1(L,*))=L^\times$ and all other cohomology of $z^1(L,*)$ vanishes. Since $X$ is geometrically irreducible and projective,  
\[
p^*:L^\times\to \Gamma(X_L,\sO_{X_L}^\times)
\]
is an isomorphism, and thus the cone of $z^1(L,*)\xrightarrow{p^*} z^1(X_L,*)$ has only cohomology in degree 0, namely 
\[
H^0(\cone(p^*))= \CH^1(X_{L}).
 \]
Thus the  sequence \eqref{eqn:ConeSeq}   is naturally isomorphic  (in $D(\PSh_\et(k))$) to the canonical sequence
\begin{multline}\label{eqn:CanSeq}
L\mapsto\\
\left( H^{-1}(z^1(X_{L},*))[1]\to \tau_{\ge -1}z^1(X_{L},*)\to H^0(z^1(X_{L},*))\to H^{-1}(z^1(X_{L},*))[2]\right).
\end{multline}

We can explicitly calculate a co-cycle representing $\beta_A$ as follows: Take $L/k$ to be a Galois extension with group $G$ such that $A_{L}$ is a matrix algebra over $L$. Then \eqref{eqn:CanSeq} gives a distinguished triangle in the derived category of $G$-modules, so we have in particular the connecting homomorphism
\[
\partial_L:H^0(G,H^0(z^1(X_{L},*)))\to H^2(G; H^{-1}(z^1(X_{L},*)))=H^2(G;L^\times)
\]
Also $X_L\cong \P^d_L$. As $H^0(z^1(X_{L},*))=\CH^1(X_L)$, $H^0(z^1(X_{L},*))$ has a canonical $G$-invariant generator $1$, namely the element corresponding to $c_1(\sO(1))$. We can apply $\partial_L$ to $1$, giving the element $\partial_L(1)\in H^2(G;L^\times)$ which maps to $\beta_A$ under the canonical map
\[
H^2(G,L^\times)\to H^2_\et(k,\G_m).
\]

Since $A_L$ is a matrix algebra over $L$,  $A$ is given by a 1-cocycle 
 \[
\{ \bar{g}_\sigma\ |\ \sigma\in G\}\in Z^1(G,\PGL_\ell(L))
\]
and $X$ is the form of $\P^d$ defined by $\{ \bar{g}_{\sigma}\}$. This mean that there is an $L$ isomorphism $\psi:X_{L}\to\P^d_{L}$  such that, for each $\sigma\in G$,  we have
\[
\bar{g}_\sigma:=\psi\circ{}^\sigma\psi^{-1},
\]
under the usual identification $\Aut_L(\P^d_L)=\PGL_{d+1}(L)$.

Lifting $\bar{g}_\sigma$ to $g_\sigma\in\GL_\ell(L)$ and defining $c_{\tau,\sigma}\in L^\times$ by
\[
c_{\tau,\sigma}\id:= g_\tau{}^\tau g_\sigma g^{-1}_{\tau\sigma}
\]
we have the co-cycle $\{c_{\tau,\sigma}\}\in Z^2(G,L^\times)$ representing $[A]$.

For a $G$-module $M$, let $(C^*(G;M),\hat{d})$ denote the standard co-chain complex computing $H^*(G;M)$, i.e., $C^n(G;M)$ is a group of $n$ co-chains of $G$ with values in $M$.  We will show that  
 $\partial_L(1)=\{c_{\tau,\sigma}\}$ in $H^2(G,L^\times)$  by applying $C^*(G;-)$ to the sequence \eqref{eqn:CanSeq} and making an explicit computation of the boundary map.

Fix a hyperplane $H\subset \P^d_k$.  Then  $D:=\psi^*(H_{L})\in z^1(X_{L},*)^0$ represents the positive generator  $1\in \CH^1(X_{L})\cong \Z$. As the class of $D$ in $\CH^1(X_{L})$ is $G$-invariant, there is for  each $\sigma\in G$ a rational function $f_\sigma$ on $X_{L}$ such that
\[
\Div(f_\sigma)={}^\sigma D-D.
\]
Given $\tau,\sigma\in G$, we thus have
\[
\Div(f_\sigma^\tau f_{\tau\sigma}^{-1}f_\tau)={}^{\tau\sigma}D-{}^\tau D-({}^{\tau\sigma}D-D)+{}^\tau D-D=0
\]
so there is a $c'_{\tau,\sigma}\in \Gamma(X_L,\sO^\times_{X_L})= L^\times$ with
\[
c'_{\tau,\sigma}= f_\sigma^\tau f_{\tau\sigma}^{-1} f_\tau.
\]
Using the fact that 
\[
{}^\sigma D=\psi^*(\bar{g}_\sigma(H_{L}))
\]
one can easily calculate that
\[
c'_{\tau,\sigma}=c_{\tau,\sigma}.
\]
 Indeed,  take a $k$-linear form $L_0$ so that $H$ is the hyperplane defined by $L_0=0$. Let
\[
F_\sigma:=\frac{L_0\circ g_\sigma^{-1}}{L_0}
\]
so $\Div(F_\sigma)=\bar{g}_\sigma(H)-H$. Letting $f_\sigma:=\psi^*F_\sigma$, we have 
\[
\Div(f_\sigma)=\psi^*(\Div(F_\sigma))=\psi^*(\bar{g}_\sigma(H)-H)={}^\sigma D-D,
\]
and
\[
{}^\tau f_\sigma=\psi^*(\frac{L_0\circ{}^\tau g_\sigma^{-1}\circ g_\tau^{-1}}{L_0\circ g_\tau^{-1}}).
\]
Thus
\begin{align*}
c'_{\tau,\sigma}&=
\psi^*\left(\frac{L_0\circ{}^\tau g_\sigma^{-1}\circ g_\tau^{-1}}{L_0\circ g_\tau^{-1}}\right)\cdot\psi^*\left(\frac{L_0\circ g_{\tau\sigma}^{-1}}{L_0}\right)^{-1}\cdot  \psi^*\left(\frac{L_0\circ g_\tau^{-1}}{L_0}\right)\\
&=\psi^*\left(\frac{L_0\circ {}^\tau g_\sigma^{-1}g_\tau^{-1}}{L_0\circ g_{\tau\sigma}^{-1}}\right)\\
&=c_{\tau,\sigma}
\end{align*}

On the other hand, we can calculate the boundary $\partial_L(1)$ by lifting the generator $1=[D]\in \CH^1(X_L)^G$ to the element $D\in z^1(X_L,*)^0$ and taking  \v{C}ech co-boundaries. Explicitly, let $\Gamma_\sigma\subset X_L\times\Delta^1$ be the closure of graph of $f_\sigma$, after identifying $(\Delta^1,0,1)$ with $(\P^1\setminus\{1\},0,\infty)$. Define $\Gamma_{c_{\sigma,\tau}}\in z^1(L,*)^{-1}$ similarly as the point of $\Delta^1_L$ corresponding to $c_{\tau,\sigma}\in \A^1(k)\subset \P^1(k)$, and let $\delta$ denote the boundary in the complex $z^1(X_L,*)$. For $\sigma\in G$, we have
\[
\delta^{-1}(\Gamma_\sigma)={}^\sigma D-D=\hat{d}^0(D)_\sigma.
\]
Since $H^{-1}(z^1(X_L,*))=\Gamma(X_L,\sO_{X_L}^\times)=L^\times$, there is for each $\sigma,\tau\in G$, an element $B_{\sigma,\tau}\in z^1(X_L,2)$ with
\begin{align*}
p^*\Gamma_{c_{\sigma,\tau}}&=
{}^\tau\Gamma_\sigma-\Gamma_{\tau\sigma}+\Gamma_\tau+\delta^{-2}(B_{\sigma,\tau})\\
&=\hat{d}^1(\sigma\mapsto \Gamma_\sigma)_{\tau,\sigma}\in \tau_{\ge-1}z^1(X_L,*)^{-1}.
\end{align*}
Thus  
\[
\partial_L([D])=\{c_{\sigma,\tau}\}\in H^2(G, H^{-1}(z^1(L,*)))= H^2(G, L^\times).
\]
This completes the computation of $\partial_L(1)$ and the proof of the  proposition.
\end{proof}

\begin{thm}\label{thm:SK0} Let $A$ be a central simple algebra over $k$ of square-free index $e$. Let $n\ge 0$, and assume that the Beilinson-Lichtenbaum conjecture holds in weights $\le n+1$ at all primes dividing $e$.\\
\\
1.  For  $m<n$, the reduced norm
\[
\Nrd: H^m(k,\Z_A(n))\to H^m(k,\Z(n))
\]
is an isomorphism. \\
\\
2. We have an exact sequence
\begin{multline*}
0\to H^n(k,\Z_A(n))\xrightarrow{\Nrd}H^n(k,\Z(n))\simeq K_n^M(k)\\\xrightarrow{\cup[A]} H^{n+2}_\et(k,\Z/e(n+1))
\to H^{n+2}_\et(k(X),\Z/e(n+1)).
\end{multline*}
\\
3. ($n=1$) $SK_1(A)=0$. More precisely, we have an exact sequence
\[
0\to K_1(A)\xrightarrow{\Nrd}K_1(k)\xrightarrow{\cup[A]}H^3_\et(k,\Z/e(2))\to H^3_\et(k(X),\Z/e(2)).
\]
\\
4. ($n=2$) $SK_2(A)=0$. More precisely, we have an exact sequence
\[
0\to K_2(A)\xrightarrow{\Nrd}K_2(k)\xrightarrow{\cup[A]}H^4_\et(k,\Z/e(3))\to H^4_\et(k(X),\Z/e(3))
\]
\end{thm} 
To explain the map $\cup[A]$ in (2), (3) and (4): We have isomorphisms
\begin{align*}
&K_1(k)=k^\times\cong H^1(k,\Z(1))\\
&K_2(k)\cong H^2(k,\Z(2)) \\
&H^n_\et(k,\G_m)\cong H^{n+1}_\et(k,\Z(1)^\et).
\end{align*}
Thus we have $[A]\in H^3_\et(k,\Z(1)^\et)$ and cup product maps
\[
H^n(k,\Z(n))\to H^n_\et(k,\Z(n)^\et)\xrightarrow{\cup[A]}H^{n+3}_\et(k,\Z(n+1)^\et)
\]
which obviously land into ${}_e H^{n+3}_\et(k,\Z(n+1)^\et$.
On the other hand, the exact triangle
\[\Z(n+1)^\et\by{e}\Z(n+1)^\et\to\Z/e(n+1)\by{+1}\]
and the Beilinson-Lichtenbaum conjecture in weight $n+1$ give an isomorphism
\[H^{n+2}_\et(k,\Z/e(n+1))\iso {}_e H^{n+3}_\et(k,\Z(n+1)^\et).\]

\begin{proof} As in the proof of corollary~\ref{cor:CohComp}, it suffices to handle the case of $A$ of prime degree over $k$. Thus, (1) follows from theorem~\ref{thm:CohComp}(1).

For (2), applying $\alpha^*$ to the distinguished triangle
 \[
 \Z(1)[2]\to \Omega^{d-1}_Tf^\mot_{d-1}M(X)\to \Z_A \to \Z(1)[2]
 \]
we have 
 \[
 \Z(1)^\et[2]\to \alpha^* \Omega^{d-1}_Tf^\mot_{d-1}M(X)\to \Z^\et\xrightarrow{\partial} \Z(1)^\et[3]
 \]
 It follows from proposition~\ref{prop:Boundary} that the $\partial$ is given by cup product with $[A]\in H^3_\et(k,\Z(1)^\et)$.  Since the map $\partial_n$ in theorem~\ref{thm:CohComp} is just the map induced by $\partial$ after tensoring with $\Z(n)^\et[n]$, (2) is proven in the form of an exact sequence
\begin{multline*}
0\to H^n(k,\Z_A(n))\xrightarrow{\Nrd}H^n(k,\Z(n))\\\xrightarrow{\cup[A]} H^{n+3}_\et(k,\Z(n+1)^\et)
\to H^{n+3}_\et(k(X),\Z(n+1)^\et).
\end{multline*}

But the Beilinson-Lichtenbaum conjecture in weight $n+1$, applied both to $k$ and $k(X)$, shows that in the commutative diagram
\[\begin{CD}
H^{n+2}_\et(k,\Z/e(n+1))@>>> H^{n+2}_\et(k(X),\Z/e(n+1))\\
@V{\partial}VV @V{\partial}VV \\
H^{n+3}_\et(k,\Z(n+1)^\et)@>>> H^{n+3}_\et(k(X),\Z(n+1)^\et)
 \end{CD}\]
both horizontal maps have isomorphic kernels, hence the form of (2) appearing in Theorem \ref{thm:SK0}.

For (3) and (4),  we have the isomorphism (corollary~\ref{cor:HigherChowMotCoh})
\[
\psi_{p,q;A}:H^p(k,\Z_A(q))\to \CH^q(k,2q-p;A)
\]
compatible with the respective reduced norm maps. From
corollary~\ref{cor:K2Iso}, the edge-homomorphism $p_{2,k;A}:\CH^2(k,2;A)\to K_2(A)$ is an isomorphism. It follows from 
theorem~\ref{thm:Comp}(1) that the edge homomorphism $p_{1,k;A}:\CH^1(k,1;A)\to K_1(A)$ is an isomorphism as well. Together with proposition~\ref{prop:NrdCommute}, this gives us the commutative diagram  for $n=1,2$:
\[
\xymatrix{
H^n(k,\Z_A(n))\ar[r]^-{\psi_{n,n;A}}\ar[d]^\Nrd& \CH^n(k,n;A) \ar[r]^-{p_{n,k;A}}\ar[d]^\Nrd
&K_n(A)\ar[d]^\Nrd\\
H^n(k,\Z(n))\ar[r]_-{\psi_{n,n;k}}& \CH^n(k,n) \ar[r]_-{p_{n,k;k}}&K_n(k)}
\]
with all horizontal maps isomorphisms. Thus, in the sequence (1), we may replace 
$H^n(k,\Z_A(n))$ with $K_n(A)$ and  $H^n(k,\Z(n))$ with $K_n(k)$ for $n=1,2$, proving (3) and
(4).
\end{proof}

\appendix

\section{Modules over Azumaya algebras}\label{AppSec:Azumaya}
We collect some basic results for use throughout the paper.

Let $R$ be a commutative ring and $A$ an Azumaya $R$-algebra.

\begin{aplem}\label{lA1} If $R$ is Noetherian, $A$ is left and right Noetherian.
\end{aplem}

\begin{proof} Indeed, $A$ is a Noetherian $R$-module, hence a Noetherian $A$-module (on the left and on the right).
\end{proof}

\begin{aplem}\label{lA2} For an $A-A$-bimodule $M$, let
\[M^A=\{m\in M\mid am = ma.\}\]
Then the functor $M\mapsto M^A$ is exact and sends injective $A-A$-bimodules to injective $R$-modules.
\end{aplem}

\begin{proof} Let $A^e=A\otimes_R A^{op}$ be the enveloping algebra of $A$. We may view $M$ as a left $A^e$-module. A special $A-A$-bimodule is $A$ itself, and we clearly have
\[M^A =\Hom_{A^e}(A,M).\]

Since $A$ is an Azumaya algebra, the map $A^e\to \End_R(A)$ is an isomorphism of $R$-algebras; via this isomorphism, $\Hom_{A^e}(A,M)$ may be canonically identified with $A^*\otimes_{\End_R(A)} M$, with $A^*=\Hom_R(A,R)$. Hence $M^A$ is the transform of $M$ under the Morita functor from $\End_R(A)$-modules to $R$-modules; since this functor is an equivalence of categories, it is exact and preserves injectives.
\end{proof}

\begin{approp}\label{pA1} For any two left $A$-modules $M,N$ and any $q\ge 0$, we have
\[\Ext^q_A(M,N)\simeq \Ext^q_R(M,N)^A.\]
(Note that $\Ext^q_R(M,N)$ is naturally an $A-A$-bimodule, which gives a meaning to the statement.)
\end{approp}

\begin{proof} The bifunctor $(M,N)\mapsto \Hom_A(M,N)$ is clearly the composition of the two functors
\[(M,N)\mapsto \Hom_R(M,N)\]
(from left $A$-modules to $A-A$-bimodules) and
\[Q\mapsto Q^A\]
(from $A-A$-bimodules to $R$-modules). Note also that, if $P$ is $A$-projective and $I$ is $A$-injective, then $\Hom_R(P,I)$ is an injective $A-A$-bimodule. The conclusion therefore follows from lemma~\ref{lA1}.
\end{proof}

\begin{apcor} \label{cA1} Let $M$ be a left $A$-module. Then $M$ is $A$-projective if and only if it is $R$-projective.
\end{apcor}

\begin{proof} If $M$ is $A$-projective, it is $R$-projective since $A$ is a projective $R$-module. The converse follows from proposition~\ref{pA1}.
\end{proof}

\begin{apcor}\label{cA2} Suppose $R$ regular of dimension $d$. Then any finitely generated left $A$-module $M$ has a left resolution of length $\le d$ by finitely generated projective $A$-modules. In particular, $A$ is regular.
\end{apcor}

\begin{proof} Since $R$ is regular, it is Noetherian and so is $A$ by lemma~\ref{lA1}. Proposition~\ref{pA1} also shows that $\Ext^{d+1}_A(M,N) =0$ for any $N$. The conclusion is now classical \cite[Ch. VI, Prop. 2.1 and Ch. V, Prop. 1.3]{CartanEilenberg}.
\end{proof}

\section{Regularity} \label{AppSec:Reg} We prove the main result on the regularity of the
functor $K(-;A)$ that we  need to compute the layers in the homotopy coniveau tower for
$G(X;\sA)$ in section~\ref{sec:SpecSeq}.

Fix a noetherian commutative ring $R$. We let $\Ralg$ denote the category of commutative $R$-algebras which are localizations of finitely generated commutative $R$-algebras.

Following Bass \cite[Ch. XII, \S 7, pp. 657--658]{bass}, for an additive functor $F:\Ralg\to \Ab$, we let $NF:\Ralg\to\Ab$ be the functor
\[
NF(A):=\ker\left(F(A[t])\to F(A[t]/(t))\right)
\]
where $A[t]$ is the polynomial algebra over $A$. We set $N^qF:=N(N^{q-1}F)$; $F$
is called \emph{regular} if $N^qF=0$ for all $q>0$ (equivalently, if $NF=0$).

For $f\in A$ the morphism $A[X]\rightarrow A[X], X\mapsto f\cdot X$
induces a group endomorphism $NF(A)\rightarrow NF(A)$. So $NF(A)$ becomes
a $\Z[T]$ module. We denote by $NF(A)_{[f]}$ the $\Z[T,T^{-1}]$ module
$\Z[T,T^{-1}]\otimes_{\Z[T]} NF(A)$. With these notations
Vorst proves the following theorem in \cite{Vorst}.

\begin{apthm}\label{theorem23}
Let $A\in \Ralg$ and let $a_1,\dots, a_n$ be elements of $A$ which generate
the unit ideal. Suppose further that the map
\[
NF(R[T]_{a_{i_0},\dots,\widehat{a_{i_j}},\dots,a_{i_p}})_{[a_{i_j}]}\rightarrow
                 NF(A[T]_{a_{i_0},\dots,a_{i_p}})
                 \]
is an isomorphism, for each set of indexes $1\leq i_0<\dots <i_p\leq n$. Then
the canonical morphism
\[
\epsilon: NF(A)\rightarrow \bigoplus_{j=1}^n NF(A_{a_j})
\]
is injective.
\end{apthm}

\begin{proof}
Compare \cite[Theorem 1.2]{Vorst} or \cite[Lemma 1.1]{LevineRevisited}.
\end{proof}

This is extended by van der Kallen, in the case of the functor $A\mapsto K_n(A)$,
to prove an \'etale descent result, namely,

\begin{apthm}\label{thm:vdK1} Let $A$ be a noetherian commutative ring such that
each zero-divisor of $A$ is contained in a minimal prime ideal of $A$. Let $A\to
B$ be an \'etale and faithfully flat extension of $A$. Then the Amitsur complex
\[
0\to N^qK_n(A)\to N^qK_n(B)\to N^qK_n(B\otimes_AB)\to \ldots
\]
is exact for each $q$ and $n$.
\end{apthm}

In fact, one can abstract van der Kallen's argument to give conditions on a
functor $F:\Ralg\to \Ab$ as above so that the conclusion of theorem~\ref{thm:vdK1}
holds for the Amitsur complex for $NF$. For this, we recall the \emph{big Witt
vectors} $W(A)$ of a commutative ring $A$, with the canonical surjection $W(A)\to
A$ and the multiplicative \emph{Teichm\"uller lifting} $A\to W(A)$ sending $a\in
A$ to $[a]\in W(A)$. We have as well the Witt vectors of length $n$, with
surjection $W(A)\to W_n(A)$; we let $F^nW(A)\subset W(A)$ be the kernel. If $M$ is
a $W(A)$-module, we say $M$ is a \emph{continuous} $W(A)$ module if $M$ is a union
of the submodules $M_n$ killed by $F^nW(A)$. Then one has

\begin{apthm}\label{vdK2} Let $F:\Ralg\to \Ab$ be a functor. Suppose that $F$
satisfies:
\begin{enumerate}
\item Given $a\in A\in\Ralg$, the natural map $F(A_a)\to F(A)_{[a]}$ is an
isomorphism.
\item Sending $a\in A$ to the endomorphism $[a]:NF(A)\to NF(A)$ extends to a
continuous $W(A)$-module structure on $NF(A)$, natural in $A$, with the
Teichm\"uller lifting $[a]\in W(A)$ acting by $[a]:NF(A)\to NF(A)$. 
\item $F$ commutes with filtered direct limits
\end{enumerate}
Let $A\in \Ralg$ be such that each zero-divisor of $A$ is contained in a minimal
prime ideal of $A$. Let $A\to B$ be an \'etale and faithfully flat extension of
$A$. Then the Amitsur complex
\[
0\to NF(A)\to NF(B)\to NF(B\otimes_AB)\to NF(B\otimes_AB\otimes_AB)\to \ldots
\]
is exact.
\end{apthm}

The main example of interest for us is the following: Let $\sA$ be a noetherian
central $R$-algebra, and let $K_n(\sA)$ be the $n$th $K$-group of the category of
finitely generated projective (left) $\sA$-modules. 

\begin{apcor}\label{cor:EtDescent}  Let $F:\Ralg\to \Ab$ be the functor
\[
F(A):=N^qK_n(\sA\otimes_RA).
\]
Then $F$ satisfies the conditions of theorem~\ref{vdK2}, hence (assuming $A$
satisfies the hypothesis on zero-divisors) if $A\to B$ is an \'etale and
faithfully flat extension of $A$, then the Amitsur complex
\[
0\to N^qK_n(\sA\otimes_RA)\to N^qK_n(\sA\otimes_RB)\to
N^qK_n(\sA\otimes_RB\otimes_AB)\to \ldots
\]
is exact.
\end{apcor}

\begin{proof}
Weibel \cite{Weibel} has shown that $N^qK_n(\sA)$ admits a $W(A)$-module
structure, satisfying the conditions  (1) and (2) of theorem~\ref{vdK2}. Since
$K$-theory commutes with filtered direct limits, this proves that the given $F$
satisfies the conditions of theorem~\ref{vdK2}, whence the result.
\end{proof}

Now let $X$ be an $R$-scheme and let $\sA$ be a sheaf of Azumaya algebras over
$\sO_X$. We have the category $\sP_{X;\sA}$ of left $\sA$-Modules $\sE$ which are
locally free as $\sO_X$-Modules. We let $K(X;\sA)$ denote the
$K$-theory spectrum of $\sP_{X;\sA}$. We extend  $K(X;\sA)$ to a spectrum which
is (possibly) not $(-1)$-connected by taking the Bass delooping, and denote this
spectrum by $KB(X;\sA)$. For $f:Y\to X$ an $X$-scheme, we write $K(Y;\sA)$ for
$K(Y;f^*\sA)$, and similarly for $KB$.

The spectra $KB(X;\sA)$ have the following properties:
\begin{enumerate}
\item There is a canonical map $K(X;\sA)\to KB(X;\sA)$, identifying $K(X;\sA)$
with is the -1-connected cover of $KB(X;\sA)$.
\item There is the natural exact sequence
\begin{multline*}
0\to KB_p(X;\sA)\to KB_p(X\times\A^1;\sA)\oplus KB_p(X\times\A^1;\sA)\\
\to KB_p(X\times\G_m;\sA)\to
KB_{p-1}(X;\sA)\to 0
\end{multline*}
called the \emph{fundamental exact sequence}.
\item If $X$ is regular, then $K(X;\sA)\to KB(X;\sA)$ is a weak equivalence.
\end{enumerate}
From now on, we will drop the notation $KB(X;\sA)$ and write $K(X;\sA)$ for the
(possibly) non-connected version.

\begin{approp}\label{prop:descent} Let $X$ be a noetherian affine $R$-scheme such
that $\sO_X$ has no nilpotent elements, and let $p:Y\to X$ be an \'etale cover.
Let $\tilde\sA$ be a sheaf of Azumaya algebras over $\sO_X$. For each point $y\in
Y$, let $Y_y:=\Spec\sO_{Y,y}$ and let $p_y:Y_y\to X$ be the map induced by $p$.
Fix and integer $q\ge1$. Suppose there is an $M$ such that,  for each smooth
affine $k$-scheme $T$,  $N^qK_n(T\times_kY_y,(p_y\circ p_2)^*\sA)=0$ for each
$y\in Y$ and each $n\le M$. Then $N^qK_n(T\times_k X;\sA)=0$ for each smooth
affine $T$ and each $n\le M$.
\end{approp}

\begin{proof} Write $X=\Spec A$. Then $\prod p_y^*:A\to B:=\prod_y\sO_{Y,y}$ is
faithfully flat and \'etale. Since $X$ is affine, $\tilde\sA$ is the sheaf
associated to a central $A$-algebra $\sA$ and since $\tilde\sA$ is a sheaf of
Azumaya algebras, each finitely generated projective left $\sA$ module is finitely
generated and projective as an $A$-module. Thus $N^qK_n(X,\tilde\sA)=N^qK_n(\sA)$.
Similarly, $N^qK_n(Y_y, p_y^*\tilde\sA)=N^qK_n(p_y^*\sA)$. By
corollary~\ref{cor:EtDescent}, $N^qK_n(\sA)=0$ for $n\ge0$. The same argument,
with $T\times X$ replacing $X$ and $T\times Y_y$ replacing $Y_y$, proves the
result for $M\ge n\ge 0$ and all $T$. To handle the cases $n<0$, use the Bass
fundamental sequence and descending induction starting with $n=0$. 
\end{proof}

\end{document}